\documentclass[11pt,a4paper]{article}

\usepackage{amsmath}
\usepackage{amssymb}
\usepackage{mathrsfs}
\usepackage{color}
\usepackage{graphicx}
\usepackage[percent]{overpic}
\usepackage{url}
\usepackage{subfig}
\usepackage{enumerate}
\usepackage{overpic}
\usepackage{amsthm}
\usepackage{moreverb}
\usepackage[pdftex,colorlinks,bookmarksopen,bookmarksnumbered,citecolor=red,urlcolor=red]{hyperref}
\usepackage{authblk}

\def\b{{\bm b}}

\def\f{{\bm f}}
\def\g{{\bm g}}
\def\n{{\bm n}}
\def\u{{\bm u}}
\def\v{{\bm v}}

\def\0{\boldsymbol{0}}
\def\ss{\boldsymbol{\sigma}}

\def\mubar{\overline{\mu}}
\def\ubar{\overline{\u}}
\def\vbar{\overline{\v}}

\def\qbar{\overline{q}}
\def\dt{\partial_t}

\def\vh{\v_h}

\def\vbarh{\vbar_h}

\def\cl {\nonumber \\}
\def\el {\nonumber }

\newtheorem{rem}{Remark}[section]

\newcommand{\bm}[1]{\mbox{\boldmath{$#1$}}}
\def\div{\nabla\cdot}

\usepackage{verbatim}
\usepackage{graphicx}
\usepackage{epstopdf}

\begin{document}
\date{}
\title{A Finite Volume approximation of the Navier-Stokes equations with nonlinear filtering stabilization}
\author[1]{Michele Girfoglio\thanks{mgirfogl@sissa.it}}
\author[2]{Annalisa Quaini\thanks{quaini@math.uh.edu}}
\author[1]{Gianluigi Rozza\thanks{grozza@sissa.it}}
\affil[1]{SISSA, International School for Advanced Studies, Mathematics Area, mathLab, via Bonomea, Trieste 265 34136, Italy}
\affil[2]{Department of Mathematics, University of Houston, Houston TX 77204, USA}
\maketitle

\begin{abstract}
We consider a Leray model with a nonlinear differential low-pass filter 
for the simulation of incompressible fluid flow at moderately large Reynolds number (in the range of a few thousands) with under-refined meshes. For the implementation of the model, we adopt the three-step algorithm Evolve-Filter-Relax (EFR). The Leray model has been extensively applied within a Finite Element (FE) framework. 
Here, we propose to combine the EFR algorithm with 
a computationally efficient Finite Volume (FV) method. 
Our approach is validated against numerical data available 
in the literature for the 2D flow past a cylinder and against experimental measurements 
for the 3D fluid flow in an idealized medical device, as recommended by the U.S. Food and Drug Administration. 
We will show that for similar levels of mesh refinement FV and FE methods provide significantly different results. 
Through our numerical experiments, we are able to provide practical directions to tune the parameters involved in the model. Furthermore, we are able to investigate the impact of mesh features (element type, 
non-orthogonality, local refinement, and element aspect ratio) and the discretization method for 
the convective term on the agreement between numerical solutions and experimental data.
\end{abstract}

\vspace*{0.5cm}

\section{Introduction}

The Direct Numerical Simulation (DNS) of the Navier-Stokes Equations (NSE) computes the evolution
of all the significant flow structures by resolving them with a properly refined mesh.
Unfortunately, when the convection dominates the dynamics, i.e. the Reynolds number 
(the dimensionless number that weights the importance of inertial forces versus viscous ones) is ``large'', 
this requires very fine meshes, making DNS computationally unaffordable for practical purposes.
Therefore, for many realistic problems the simulation of turbulent flows is performed by using different models. 

The NSE can be properly averaged (in different ways, quite often in time), 
leading to the so-called Reynolds-Averaged Navier-Stokes equations (RANS), 
or filtered (usually in space), leading to Large Eddy Simulation (LES) techniques (see, e.g., \cite{B-Ferziger}). 
In this work, we focus on the latter approach. We consider a variant of the so-called Leray model \cite{Leray1934}, 
where the small-scale effects are described by a set of equations to be added to the discrete NSE formulated 
on the under-refined mesh. The extra-problem can be devised in different ways, for instance by 
a functional splitting of the solved and unresolved scales \cite{Bazilevs2007173} or by 
resorting to the concept of ``suitability'' of weak solutions \cite{guermond2011suitable}.
We consider a variation of the Leray model for which the extra problem acts as a differential low-pass filter \cite{abigail_CMAME}.
For its  actual implementation,  we use the Evolve-Filter-Relax (EFR) algorithm proposed in \cite{layton_CMAME}. 
For a formal rewriting of the EFR algorithm as a LES technique we refer to \cite{Olshanskii2013},
while we refer to \cite{BQV} for a reformulation 
in an operator-splitting framework. 
One of the advantages of this approach is that it is easily implemented in a legacy Navier-Stokes solver,
since the filtering step requires the solution of a Stokes-like problem.

The Leray model has been extensively applied within a Finite Element (FE) framework,
although obviously other space approximations are possible. Here, we focus
on Finite Volume (FV) methods, which have been widely used in the LES context. 
So far, the application of the Leray model in a FV framework has been unexplored.
In this manuscript, we intend to fill this gap by proposing a computationally efficient
FV method for the EFR algorithm. We target applications involving incompressible fluid flow in both two and three
dimensions, with the 3D applications featuring moderately large Reynolds numbers (in
the range of a few thousands). We will show that that for similar levels of mesh refinement 
FV and FE methods provide significantly different results. Thus, changing the 
space discretization method for the EFR algorithm is not a trivial exercise.

We choose two benchmarks to showcase the features of our approach. The first is 
a classical academic benchmark: 2D flow past a cylinder \cite{John2004, Turek1996}
at Reynolds number 100. Although this flow field is not turbulent, it turns out to be very 
challenging for most numerical models and methods, especially on coarser computational grids.
This first test serves the purpose of getting a preliminary understanding of how the EFR algorithm
works in combination with a FV method, in terms of element type, level of mesh refinement, 
and choice of critical parameter values. The second benchmark has been 
issued by the U.S. Food and Drug Administration (FDA) and it involves
3D flow at Reynolds numbers up to 6500 through a nozzle.
Three independent  laboratories were requested by FDA  to perform flow visualization experiments 
on fabricated nozzle models for different flow rates \cite{hariharang}. 
This resulted in benchmark data available online to the scientific community for the validation of Computational Fluid
Dynamics (CFD) simulations \cite{fdacfd}. Available experimental measurements  enable us to check 
the effectiveness of our FV-based EFR algorithm in simulating average macroscopic quantities.
Through our numerical experiments, we are able to provide practical directions to tune the parameters 
involved in the model. Moreover, we are able to investigate the impact of mesh features 
(element type, non-orthogonality, local refinement, and element aspect ratio) 
and the discretization method for the convective term on the agreement between numerical 
solutions and experimental data. 


All the computational results presented in this article have been performed with 
OpenFOAM\textsuperscript{\textregistered} \cite{Weller1998}, an open source 
finite volume C++ library widely used by commercial and academic organizations.
See, e.g., \cite{Feymark2012, Mao2017} for numerical results obtained with LES
techniques implemented in OpenFOAM\textsuperscript{\textregistered}. 
An important outcome of this work is that the code created for it is incorporated 
in an open-source library\footnote{\url{https://mathlab.sissa.it/cse-software}}
and therefore is readily shared with the community.

This work is outlined as follows. 
In Sec.~\ref{sec:pbd}, we introduce the continuous Leray model as well as the numerical approximation proposed in \cite{layton_CMAME}. 
In Sec.~\ref{subsec:space-discrete}, 
we detail our strategy for space discretization, which combines an operator splitting
and a Finite Volume method. 
Additionally, we show how to tune certain model parameters using physical and discretization quantities. 
The numerical results for the flow around the cylinder are reported in Sec.~\ref{sec:cylinder}, 
while the comparison between numerical results and experimental data provided by 
the FDA is presented in Sec.~\ref{sec:FDA}. Conclusions are drawn in Sec.~\ref{sec:conclusions}.

\section{Problem definition}
\label{sec:pbd}

\subsection{The Navier-Stokes equations}
\label{sec:NS Equations}

We consider the motion of an incompressible viscous fluid in a time-independent domain $\Omega$ 
over a time interval of interest $(t_0, T)$. The flow is described by the incompressible Navier-Stokes equations
\begin{align}
\rho\,\dt \u + \rho\,\div \left(\u \otimes \u\right) - \div \ss & = \f\quad \mbox{ in }\Omega \times (t_0,T),\label{eq:ns-mom}\\
\div \u & = 0\quad\, \mbox{ in }\Omega \times(t_0,T),\label{eq:ns-mass}
\end{align}
endowed with the boundary conditions 

\begin{align}
\u & = \u_D\quad\ \mbox{ on } \partial\Omega_D \times(t_0,T),\label{eq:bc-ns-d}\\
\ss\cdot\n & = \g\quad \quad\mbox{ on } \partial\Omega_N \times(t_0,T),\label{eq:bc-ns-n}
\end{align}
and the initial data $\u = \u_0$ in $\Omega \times\{t_0\}$.
Here $\overline{\partial\Omega_D}\cup\overline{\partial\Omega_N}=\overline{\partial\Omega}$ and $\partial\Omega_D \cap\partial\Omega_N=\emptyset$. In addition 
$\rho$ is the fluid density, $\u$ is the fluid velocity, $\dt$ denotes the time derivative, $\ss$ is the Cauchy stress tensor, $\f$ accounts for possible body forces (such as, e.g., gravity), $\u_D,\g$ and $\u_0$ are 
given.
 Equation (\ref{eq:ns-mom}) represents the conservation of the linear momentum, while eq. (\ref{eq:ns-mass}) represents the conservation of the mass. For Newtonian fluids $\ss$ can be written as
\begin{equation}\label{eq:newtonian}
\ss (\u, p) = -p \mathbf{I} +\mu (\nabla\u + \nabla\u^T),
\end{equation}
where $p$ is the pressure and $\mu$ is the constant \emph{dynamic} viscosity.
Notice that by plugging \eqref{eq:newtonian} into eq.~\eqref{eq:ns-mom}, eq.~\eqref{eq:ns-mom}
can be rewritten as
\begin{align}
\rho\, \dt \u + \rho\,\div \left(\u \otimes \u\right) + \nabla p - 2\mu\,\Delta\u  = \f\mbox{ in }\Omega \times (t_0,T).\label{eq:ns-lapls-1}
\end{align}

In order to characterize the flow regime under consideration, we define the Reynolds number as
\begin{equation}\label{eq:re}
Re = \frac{U L}{\nu},
\end{equation}
where $\nu=\mu/\rho$ is the \emph{kinematic} viscosity of the fluid, and $U$ and $L$ are characteristic macroscopic velocity and length, respectively. For an internal flow in a cylindrical pipe, $U$ is the mean sectional velocity and $L$ is the diameter. 
For  moderately large Reynolds numbers  the effects of flow disturbances cannot be neglected, and yet Reynolds-averaged Navier-Stokes (RANS) models \cite{pope} are generally too crude.

\subsection{Leray model}\label{sec:leray_model}

In the framework of the Kolmogorov 1941 theory \cite{Kolmogorov41-1,Kolmogorov41-2}, the 
turbulent kinetic energy, which is the kinetic energy associated with eddies in the turbulent flow,
 is injected in the system at the large scales (low wave numbers). Since the large scale eddies are unstable, they break down, transferring the energy to smaller eddies. Finally, the turbulent kinetic energy 
 is dissipated by the viscous forces at the small scales (high wave numbers). This process is usually referred to as \emph{energy cascade}. The scale at which the viscous forces dissipate energy is referred to as \emph{Kolmogorov scale}. For a flow in developed turbulent regime and at statistical equilibrium, the Kolmogorov scale can be expressed as
\begin{equation}\label{eq:eta-re}
\eta=Re^{-3/4}L.
\end{equation}

In order to correctly capture the dissipated energy, 
DNS needs a mesh with spacing $h\sim\eta$. As the Reynolds number increases, DNS leads to a huge number of unknowns and prohibitive computational costs.
On the other hand, when the mesh size $h$ fails to resolve the Kolmogorov scale, the under-diffusion in the simulation leads to nonphysical computed velocities.
In some cases, this is detectable simply looking at the velocity field, which features nonphysical oscillations eventually leading to a simulation break down. 
However, in some cases the velocity field does not display oscillations, 
yet it does not correspond to the physical solution. A possible remedy to this issue is to introduce a model which filters the 
nonphysical oscillations in the velocity field and conveys the energy lost to resolved scales. 

The so called \emph{Leray model} couples the Navier-Stokes equations \eqref{eq:ns-lapls-1},\eqref{eq:ns-mass} with a differential filter. The resulting system reads

\begin{align}
\rho\, \dt \u + \rho\,\div \left(\ubar \otimes \u\right) - 2\mu \Delta\u + \nabla p = \f & \quad {\rm in}~\Omega \times
(t_0,T), \label{eq:filter-ns1}\\
\div \u = 0 & \quad {\rm in}~\Omega \times(t_0,T), \label{eq:filter-ns2} \\
-2 \alpha^2\div \left(a(\u) \nabla\ubar\right) +\ubar +\nabla \lambda = \u & \quad {\rm in}~\Omega \times
(t_0,T),\label{eq:filter-mom}\\
\div \ubar = 0 & \quad {\rm in}~\Omega \times(t_0,T).\label{eq:filter-mass}
\end{align}
Here, $\ubar$ is the \emph{filtered velocity}, $\alpha$ can be interpreted as the \emph{filtering radius} (that is, the radius of the neighborhood where the filter extracts information from the unresolved scales), 
the variable $\lambda$ is a Lagrange multiplier to enforce the incompressibility constraint for $\ubar$
and $a(\cdot)$ is a scalar function such that:
\begin{align*}
a(\u)\simeq 0 & \mbox{ where the velocity $\u$ does not need regularization;}\\
a(\u)\simeq 1 & \mbox{ where the velocity $\u$ does need regularization.}
\end{align*}
This function, called \emph{indicator function}, is crucial for the success of the Leray model. 
In Sec.~\ref{sec:indicator}, we will discuss our choice of $a(\u)$. 
Here, we mention that the choice $a(\u)\equiv 1$ corresponds to the classic Leray-$\alpha$ model \cite{Leray1934}. 
This model has the advantage of making the operator in the filter equations linear and constant in time, 
but its effectivity is rather limited, since it introduces the same amount of regularization everywhere in the domain, hence causing overdiffusion.

Equations (\ref{eq:filter-mom})-(\ref{eq:filter-mass}) require suitable boundary conditions. These are chosen to be
\begin{align}
\ubar & = \u_D \mbox{ on } \partial\Omega_D \times(t_0,T),\label{eq:bc-filter-d}\\
(2\alpha^2 a(\u)\nabla\ubar - \lambda\mathbf{I})\n & = \0 \quad\mbox{ on } \partial\Omega_N \times(t_0,T).\label{eq:bc-filter-n}
\end{align}


Even though  (\ref{eq:filter-ns1})-(\ref{eq:filter-ns2}) are linear in $(\u,p)$ and the filter problem is linear in $(\ubar,\lambda)$, the coupling is non-linear, due to the term $\div \left(\ubar \otimes \u\right)$ in eq.~(\ref{eq:filter-ns1}) and the term $a(\u)\nabla\ubar$ in eq.~(\ref{eq:filter-mom}) (when $a(\cdot)$ is not constant). 

\subsection{Time discrete problem}
\label{subsec:time-discrete}

To discretize in time problem (\ref{eq:filter-ns1})-(\ref{eq:filter-mass}), let $\Delta t \in \mathbb{R}$, $t^n = t_0 + n \Delta t$, with $n = 0, ..., N_T$ and $T = t_0 + N_T \Delta t$. Moreover, we denote by $y^n$ the approximation of a generic quantity $y$ at the time $t^n$. 

For the time discretization of system (\ref{eq:filter-ns1})-(\ref{eq:filter-mass}), we adopt a Backward Differentiation Formula of order 2 (BDF2), 
see e.g. \cite{quarteroni2007numerical}. The Leray system discretized in time reads: given $\u^0$, for $n \geq 0$ find the solution $(\u^{n+1}, p^{n+1},\ubar^{n+1},\lambda^{n+1})$ of system:
\begin{align}
\rho\, \frac{3}{2\Delta t}\, \u^{n+1} + \rho\, \div \left(\ubar^{n+1} \otimes \u^{n+1}\right) - 2\mu\Delta\u^{n+1} +\nabla p^{n+1} & = \b^{n+1},\label{eq:disc_filter_ns-1}\\
\div \u^{n+1} & = 0, \label{eq:disc_filter_ns-2}\\
-\alpha^2\div \left(a(\u^{n+1}) \nabla\ubar^{n+1}\right) +\ubar^{n+1} +\nabla \lambda^{n+1} & = \u^{n+1},\label{eq:disc_filter_mom}\\
\div \ubar^{n+1} & = 0, \label{eq:disc_filter_mass}
\end{align}
where $\b^{n+1} = \f^{n+1} + (4\u^n - \u^{n-1})/(2\Delta t)$. 
Obviously, other discretization schemes are possible. However, for clarity of exposition we will restrict the description of the approach to the case of BDF2.

A monolithic approach for problem (\ref{eq:disc_filter_ns-1})-(\ref{eq:disc_filter_mass}) would lead to high computational costs, making the advantage compared to DNS questionable. 
To decouple the Navier-Stokes system (\ref{eq:disc_filter_ns-1})-(\ref{eq:disc_filter_ns-2}) from the filter system (\ref{eq:disc_filter_mom})-(\ref{eq:disc_filter_mass}) at the time $t^{n+1}$, we have two options:
\begin{enumerate}[i)]
\item[1.] {\it Filter-then-solve}: Solve the filter equations (\ref{eq:disc_filter_mom})-(\ref{eq:disc_filter_mass}) first, with $\u^{n+1}$ replaced by a 
suitable extrapolation $\u^*$ and $a(\u^{n+1})$ replaced by $a(\u^*)$, and then solve equations (\ref{eq:disc_filter_ns-1})-(\ref{eq:disc_filter_ns-2}) with advection field given by the filtered velocity previously computed.
\item[2.] {\it Solve-then-filter}: Solve equations (\ref{eq:disc_filter_ns-1})-(\ref{eq:disc_filter_ns-2}) first, replacing  the advection field $\ubar^{n+1}$ with  a suitable extrapolation $\u^*$, and then solve the 
filter problem (\ref{eq:disc_filter_mom})-(\ref{eq:disc_filter_mass}).
\end{enumerate}
In this work, we will focus on approach 2. 
To keep the computational costs low, we adopt a semi-implicit approach, i.e.~we perform only one iteration per time step.
In particular, we consider a modified version of approach 2 called Evolve-Filter-Relax (EFR), which was 
first proposed in \cite{layton_CMAME}. This algorithm reads as follows: given the velocities $\u^{n-1}$ and $\u^{n}$, at $t^{n+1}$:

\begin{enumerate}[i)]
\item \textit{evolve}: find intermediate velocity and pressure $(\v^{n+1},q^{n+1})$ such that
\begin{align}
\rho\, \frac{3}{2\Delta t}\, \v^{n+1} + \rho\, \div \left(\u^* \otimes \v^{n+1}\right) - 2\mu\Delta\v^{n+1} +\nabla q^{n+1} & = \b^{n+1},\label{eq:evolve-1.1}\\
\div \v^{n+1} & = 0\label{eq:evolve-1.2},
\end{align}
where $\u^* = 2 \u^n-\u^{n-1}$.

\item \textit{filter}: find $(\vbar^{n+1},\lambda^{n+1})$ such that
\begin{align}
-\alpha^2\div \left(a(\v^{n+1}) \nabla\vbar^{n+1}\right) +\vbar^{n+1} +\nabla \lambda^{n+1} & = \v^{n+1}, \label{eq:evolve-2.1}\\
\div \vbar^{n+1} & = 0 \label{eq:filter-1.2}.
\end{align}
\item \textit{relax}: set 
\begin{align}
\u^{n+1}&=(1-\chi)\v^{n+1} + \chi\vbar^{n+1}, \label{eq:relax-1} \\
p^{n+1}&= q^{n+1},  \label{eq:relax-2}
\end{align}
where $\chi\in(0,1]$ is a relaxation parameter.
\end{enumerate}

\begin{rem}\label{rem:gen_Stokes}
Filter problem \eqref{eq:evolve-2.1}-\eqref{eq:filter-1.2} can be considered 
a generalized Stokes problem. In fact, by dividing eq.~\eqref{eq:evolve-2.1}
by $\Delta t$ and rearranging the terms we obtain:
\begin{align}
\frac{\rho}{\Delta t} \vbar^{n+1}  - \div \left( \mubar \nabla\vbar^{n+1}\right) + \nabla \qbar^{n+1} & = \frac{\rho}{\Delta t} \v^{n+1}, \quad \mubar = \rho \frac{\alpha^2}{\Delta t} a(\v^{n+1}), \label{eq:filter-1.1}
\end{align}
where $\qbar^{n+1} = \rho \lambda^{n+1}/\Delta t$. Problem \eqref{eq:filter-1.1},\eqref{eq:filter-1.2} can be
seen as a time dependent Stokes problem with a non-constant viscosity $\mubar$, discretized
by the Backward Euler (or BDF1) scheme. 
A solver for problem \eqref{eq:filter-1.1},\eqref{eq:filter-1.2} can then be obtained 
by adapting a standard linearized Navier-Stokes solver. 
Notice that for $\alpha^2/\Delta t \simeq h^2/\Delta t \to 0$ we have that $\mubar \rightarrow 0$ and $\vbarh^{n+1} \to \vh^{n+1}$.
\end{rem}

\begin{rem}\label{rem:p_relax}
In the EFR algorithm proposed in \cite{layton_CMAME} there is no relaxation for the pressure, i.e.~the end-of-step
pressure is set equal to the pressure of the Evolve step. 
In \cite{BQV}, two relaxations for the pressure were considered: $p^{n+1} = q^{n+1} + \frac{3}{2} \chi \qbar^{n+1}$ or
$p^{n+1} =(1-\chi)q^{n+1} + \chi\qbar^{n+1}$. Notice while $\qbar^{n+1}$ has the same dimensional 
units as $q^{n+1}$, $\lambda^{n+1}$ does not.
\end{rem}

In the rest of the paper, we will call EFR algorithm \eqref{eq:evolve-1.1}, \eqref{eq:evolve-1.2}, 
\eqref{eq:filter-1.1}, \eqref{eq:filter-1.2}-\eqref{eq:relax-2}.
We will call the algorithm simply Evolve-Filter (EF) when $\chi = 1$, since there is no actual relaxation step.
In the next subsection, we will consider an indicator function that leads to a nonlinear alternative to the Leray-$\alpha$ model.

 \subsection{The indicator function}\label{sec:indicator}

Different choices of $a(\cdot)$ have been proposed and compared in \cite{Borggaard2009,layton_CMAME,O-hunt1988,Vreman2004,Bowers2012}. 
Here, we focus on a class of deconvolution-based indicator functions:
\begin{equation}
a(\v) = a_{D}(\v) = \left|  \v - D (F(\v)) \right|^2, \label{eq:a_deconv}
\end{equation}
where $F$ is a linear filter (an invertible, self-adjoint, compact operator from a Hilbert space to itself)
and $D$ is a bounded regularized approximation of $F^{-1}$.
A popular choice for $D$ is the Van Cittert deconvolution operator $D_N$, defined as
\begin{equation}
D_N = \sum_{n = 0}^N (I - F)^n. \el
\end{equation}
The evaluation of $a_D$ with $D=D_N$ (deconvolution of order $N$) requires then to apply 
the filter $F$ a total of $N+1$ times.
Since $F^{-1}$ is not bounded, in practice $N$ is chosen to be small, 
as the result of a trade-off between accuracy (for a regular solution) and filtering (for a non-regular one).
In this paper  we consider $N = 0$, corresponding to $D_0=I$. For this choice of $N$, the indicator function (\ref{eq:a_deconv}) becomes
\begin{align}
a_{D_0}(\v) = \left|  \v - F(\v) \right|. \label{eq:a_D0_a_D1}
\end{align}

We select $F$ to be the linear Helmholtz filter operator $F_H$ defined by 
\begin{equation}
	F=F_H = \left(I - \alpha^2  \Delta \right)^{-1}. \el
\end{equation}
It is possible to prove \cite{Dunca2005} that
\begin{align}
\v - D_N (F_H(\v)) = (-1)^{N+1} \alpha^{2N+2} \Delta^{N+1} F_H^{N+1} \v. \el
\end{align}
Therefore, $a_{D_N}(\v)$ is close to zero in the regions of the domain where $\v$ is smooth.
We remark that finding $F_H(\v^{n+1}) = \tilde{\v} ^{n+1}$ is equivalent to finding
$\tilde{\v} ^{n+1}$ such that:
\begin{equation}\label{eq:vtilde}
\tilde{\v} ^{n+1} - \alpha^2  \Delta \tilde{\v} ^{n+1}= \v^{n+1}. 
\end{equation}

\def\resY{J_{\alpha^2}}
\def\regY{\mathcal{L}_{\alpha^2}}


We will refer to eq.~\eqref{eq:filter-ns1}-\eqref{eq:filter-mass} with indicator function given by 
\eqref{eq:a_D0_a_D1} as Leray-$\alpha$-NL model.

\section{Space discrete problem: the Finite Volume approximation}\label{subsec:space-discrete}

In this section we discuss the space discretization of problems \eqref{eq:evolve-1.1}-\eqref{eq:evolve-1.2}
and \eqref{eq:filter-1.1},\eqref{eq:filter-1.2} for both the classic Leray-$\alpha$ and the Leray-$\alpha$-NL model.

We adopt the Finite Volume (FV) approximation that is derived
directly from the integral form of the governing equations. 
We have chosen to implement the EFR algorithm within
the finite volume C++ library OpenFOAM\textsuperscript{\textregistered} \cite{Weller1998}.
We partition the computational domain $\Omega$ into cells or control volumes $\Omega_i$,
with $i = 1, \dots, N_{c}$, where $N_{c}$ is the total number of cells in the mesh. 
Let  \textbf{A}$_j$ be the surface vector of each face of the control volume, 
with $j = 1, \dots, M$. The value of $M$ depends on the dimension of the domain and the type of mesh 
(hexahedral vs prismatic). 

\subsection{Numerical discretization for the \emph{evolve} step of the problem} 
\label{subsec:space-discrete-evolve}


The integral form of eq.~(\ref{eq:evolve-1.1}) for each volume $\Omega_i$ is given by:

\begin{align}\label{eq:evolveFVtemp-1.1}
\rho\, \frac{3}{2\Delta t}\, \int_{\Omega_i} \v^{n+1} d\Omega &+ \rho\, \int_{\Omega_i} \div \left(\u^* \otimes \v^{n+1}\right) d\Omega - 2\mu \int_{\Omega_i} \Delta\v^{n+1} d\Omega \cl
&+ \int_{\Omega_i}\nabla q^{n+1} d\Omega  = \int_{\Omega_i}{\bm b}^{n+1} d\Omega.
\end{align}
By applying the Gauss-divergence theorem, eq.~\eqref{eq:evolveFVtemp-1.1} becomes:

\begin{align}\label{eq:evolveFV-1.1}
\rho\, \frac{3}{2\Delta t}\, \int_{\Omega_i} \v^{n+1} d\Omega &+ \rho\, \int_{\partial \Omega_i} \left(\u^* \otimes \v^{n+1}\right) \cdot d\textbf{A} - 2\mu \int_{\partial \Omega_i} \nabla\v^{n+1} \cdot d\textbf{A} \cl
&+ \int_{\partial \Omega_i}q^{n+1} d\textbf{A}  = \int_{\Omega_i}{\bm b}^{n+1} d\Omega.
\end{align}
Each term in eq.~\eqref{eq:evolveFV-1.1} is approximated as follows: 

\begin{itemize}
\item[-] \textit{Gradient term}: 

\begin{align}\label{eq:grad}
\int_{\partial \Omega_i}q^{n+1} d\textbf{A} \approx \sum_j^{} q^{n+1}_j \textbf{A}_j, 
\end{align} 
where $q_j$ is the value of the pressure relative to centroid of the $j^{\text{th}}$ face. In OpenFOAM\textsuperscript{\textregistered} solvers, the face center pressure values $q_j$ are typically obtained from 
the cell center values by means of a linear interpolation scheme. Such scheme is rigorously second-order accurate only on structured meshes \cite{Syrakos2017}.


\item[-] \textit{Convective term}: 
\begin{align}\label{eq:conv}
\int_{\partial \Omega_i} \left(\u^* \otimes \v^{n+1}\right) \cdot d\textbf{A} \approx \sum_j^{} \left(\u^{*}_j \otimes \v^{n+1}_j\right) \cdot \textbf{A}_j = \sum_j^{} \varphi^*_j \v^{n+1}_j, \quad \varphi^*_j = \u^{*}_j \cdot \textbf{A}_j,
\end{align} 
where $\u^{*}_j$ and $\v^{n+1}_j$ are respectively the extrapolated convective velocity and the fluid velocity relative 
to the centroid of each control volume face. In \eqref{eq:conv}, $\varphi^*_j$ is the convective flux associated to $\u^{*}$ through face $j$ of the control volume. In OpenFOAM\textsuperscript{\textregistered} solvers, 
the convective flux at the cell faces is typically a linear interpolation of the values from the adjacent cells. 
Also $\v^{n+1}$ needs to be approximated at cell face $j$ in order to get the face 
value $\v^{n+1}_j$. Different interpolation methods can be applied: central, upwind, second order upwind and blended differencing schemes \cite{jasakphd}.
\item[-] \textit{Diffusion term}: 
\begin{align}
\int_{\partial \Omega_i} \nabla\v^{n+1} \cdot d\textbf{A} \approx \sum_j^{} (\nabla\v^{n+1})_j \cdot \textbf{A}_j, \el
\end{align} 
where $(\nabla\v^{n+1})_j$ is the gradient of $\v^{n+1}$ at face $j$. 
We are going to briefly explain how $(\nabla\v^{n+1})_j$ is approximated with
second order accuracy on a structured, orthogonal mesh. Let $P$ and $Q$ be two neighboring control volumes.
The term $(\nabla\v^{n+1})_j$ is evaluated by subtracting
 the value of velocity at the cell centroid on the $P$-side of the face, denoted with $\v^{n+1}_P$,
 from the value of velocity at the centroid on the $Q$-side, denoted with $\v^{n+1}_Q$,
 and dividing by the magnitude of the distance vector $\textbf{d}_j$ connecting the two cell centroids:
\begin{align}
(\nabla\v^{n+1})_j \cdot \textbf{A}_j = \dfrac{\v^{n+1}_Q - \v^{n+1}_P}{|\textbf{d}_j|} |\textbf{A}_j|. \el
\end{align} 
For non-structured, non-orthogonal meshes
(see Fig.~\ref{fig:gradient_image}), an explicit non-orthogonal correction has to be added to the orthogonal component
in order to preserve second order accuracy. See \cite{jasakphd} for details.
\begin{figure}[h!]
\centering
\includegraphics[width=0.5\textwidth]{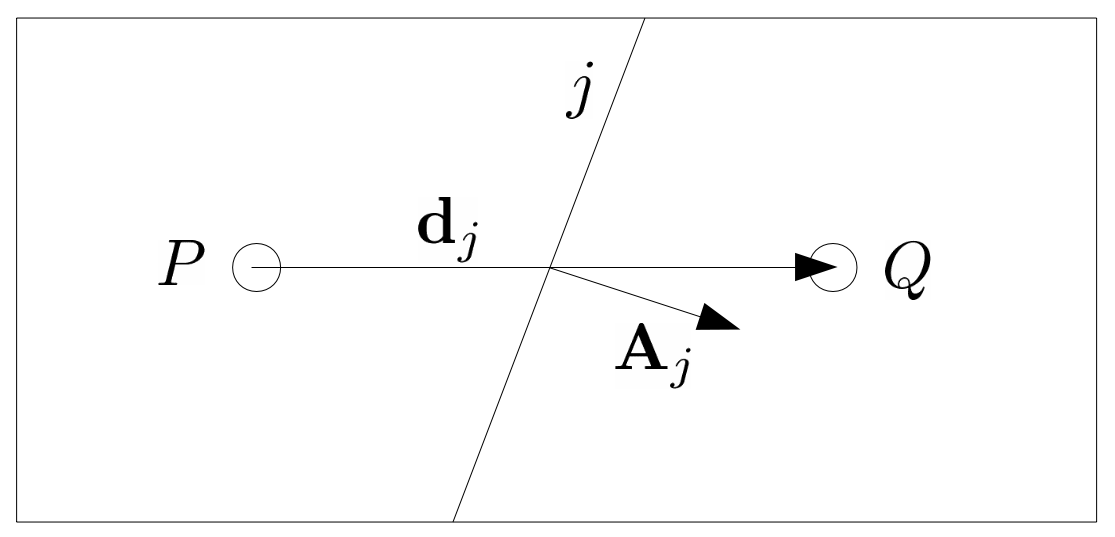}
\caption{Close-up view of two non-orthogonal control volumes in a 2D configuration.}
\label{fig:gradient_image}
\end{figure}
\end{itemize}

Let us denote with $\v^{n+1}_i$ and ${\bm b}^{n+1}_i$ the average velocity 
and source term in control volume $\Omega_i$, respectively.
Moreover, we denote with $\v^{n+1}_{i,j}$ and $q^{n+1}_{i,j}$ the velocity and pressure
associated to the centroid of face $j$ normalized by the volume of $\Omega_i$.
Then the discretized form of eq.~\eqref{eq:evolveFV-1.1}, divided by the control volume 
$\Omega_i$, can be written as:
\begin{align}\label{eq:evolveFV-1.1_disc}
\rho\, \frac{3}{2\Delta t}\, \v^{n+1}_i &+ \rho\, \sum_j^{} \varphi^*_j \v^{n+1}_{i,j} - 2\mu \sum_j^{} (\nabla\v^{n+1}_i)_j \cdot \textbf{A}_j + \sum_j^{} q^{n+1}_{i,j} \textbf{A}_j  = {\bm b}^{n+1}_i.
\end{align}
Following \cite{Rhie1983}, we now write eq.~\eqref{eq:evolveFV-1.1_disc} in semi-discretized form, i.e.
with the pressure term in continuous form while all the other terms
are in discrete form 
\begin{align}\label{eq:semidisc}
\v^{n+1} = \frac{2  \Delta t}{3} \left(\textbf{H}(\v^{n+1}) - \nabla q^{n+1} \right),
\end{align}
with
\begin{align}
\textbf{H}(\v^{n+1}_i) = -\rho \sum_j^{} \varphi^*_j \v^{n+1}_{i,j} + 2\mu \sum_j^{} (\nabla\v^{n+1}_i)_j \cdot \textbf{A}_j + {\bm b}^{n+1}_i. \el
\end{align}
Next, we take the divergence of the eq.~(\ref{eq:semidisc}) and make use of eq.~(\ref{eq:evolve-1.2}) to obtain
\begin{align}\label{eq:q_poissonL_tmp}
\Delta q^{n+1} = \nabla \cdot \textbf{H}(\v^{n+1})
\end{align}
By integrating eq.~\eqref{eq:q_poissonL_tmp} over the control volume $\Omega_i$ 
and applying the Gauss-divergence theorem, we get
\begin{align}
\int_{\partial \Omega_i} \nabla q^{n+1} d \textbf{A} = \int_{\partial \Omega_i} \textbf{H}(\v^{n+1}) \cdot d \textbf{A}. \el
\end{align}
So, the space discretized eq.~\eqref{eq:q_poissonL_tmp}, divided by the control volume 
$\Omega_i$, can be expressed as:
\begin{align}\label{eq:q_poissonL_iii}
\sum_j^{} (\nabla q^{n+1})_j \cdot \textbf{A}_j = \sum_j^{} (\textbf{H}(\v_i^{n+1}))_j \cdot \textbf{A}_j
\end{align}
In eq.~\eqref{eq:q_poissonL_iii}, $(\nabla q^{n+1})_j$ is the gradient of $q^{n+1}$ at faces $j$ 
and it is approximated in the same way as $(\nabla\v^{n+1})_j$.
Finally, the fully discretized form of problem \eqref{eq:evolve-1.1}-\eqref{eq:evolve-1.2} is given by system~\eqref{eq:evolveFV-1.1_disc},\eqref{eq:q_poissonL_iii}. In this work, we choose 
a partitioned approach to deal with the pressure-velocity coupling. The partitioned algorithms 
available in OpenFOAM\textsuperscript{\textregistered} are SIMPLE \cite{SIMPLE} for steady-state problems, and PISO \cite{PISO} and PIMPLE \cite{PIMPLE} for transient problems. For the results reported in 
Sec.~\ref{sec:num_res}, we used the PISO algorithm. The splitting of operations in the solution of the discretised momentum and pressure equations gives rise to a formal order of accuracy of the order of powers of $\Delta t$ depending on the number of operation-splittings used (see \cite{PISO} for more details).
We remark that in the OpenFOAM\textsuperscript{\textregistered} implementation of the PISO solver the mass flux is modified through an additional term than can cause artificial dissipation and is not a part of the original PISO algorithm in \cite{PISO}. 
See \cite{ddtPhiCorr} for more details.



\subsection{Numerical discretization for the \emph{filter} step of the problem} 
\label{subsec:space-discrete-filter}


In this subsection, we present the discretization of the filter problem \eqref{eq:filter-1.1},\eqref{eq:filter-1.2}.
The integral form of the eq.~(\ref{eq:filter-1.1}) for each volume $\Omega_i$ is given by:
\begin{align}
\frac{\rho}{\Delta t}\int_{\Omega_i} \vbar^{n+1} d \Omega - \int_{\Omega_i} \div \left(\mubar \nabla\vbar^{n+1}\right) d \Omega + \int_{\Omega_i} \nabla \qbar^{n+1} d \Omega & = \frac{\rho}{\Delta t} \int_{\Omega_i} \v^{n+1} d \Omega. \label{eq:disc_filter_momFV}
\end{align}
By making use of the Gauss-divergence theorem, eq.~\eqref{eq:disc_filter_momFV} becomes:
\begin{align}\label{eq:FV1}
\frac{\rho}{\Delta t}\int_{\Omega_i} \vbar^{n+1} d \Omega - \int_{\partial \Omega_i} \mubar \nabla\vbar^{n+1} \cdot d\textbf{A}  + \int_{\partial \Omega_i}\qbar^{n+1} d\textbf{A} & = \frac{\rho}{\Delta t} \int_{\Omega_i}{\v}^{n+1} d\Omega.
\end{align}
Like in Sec.~\ref{subsec:space-discrete-evolve}, we proceed with providing the 
approximation of each term in eq.~\eqref{eq:FV1}:

\begin{itemize}
\item[-] \textit{Gradient term}: 
\begin{align}
\int_{\partial \Omega_i}\qbar^{n+1} d\textbf{A} \approx \sum_j^{} \qbar^{n+1}_j \textbf{A}_j,
\end{align} 
where $\qbar_j^{n+1}$ is the value of the auxiliary pressure associated to the centroid of the $j^{\text{th}}$ face.
\item[-] \textit{Diffusion term}: 
\begin{align}
\int_{\partial \Omega_i} \mubar\nabla\vbar^{n+1} \cdot d\textbf{A} \approx \sum_j^{} \mubar_j(\nabla\vbar^{n+1})_j \cdot \textbf{A}_j,
\end{align} 
where $(\nabla\vbar^{n+1})_j$ is the gradient of $\vbar^{n+1}$ at face $j$. It is approximated in the same
way as $(\nabla\v^{n+1})_j$; see Sec.~\ref{subsec:space-discrete-evolve}. 
\end{itemize}

Upon division by the volume of $\Omega_i$, the discretized form of the eq.~\eqref{eq:FV1} can be written as 
\begin{align}\label{eq:evolveFV-2.1_disc}
\frac{\rho}{\Delta t} \vbar^{n+1}_i - \sum_j^{} \mubar_j(\nabla\vbar^{n+1}_i)_j \cdot \textbf{A}_j +  \sum_j^{} \qbar^{n+1}_{i,j} \textbf{A}_j = \frac{\rho}{\Delta t} \v^{n+1}_i.
\end{align}
In eq.~\eqref{eq:evolveFV-2.1_disc}, we denoted with $\vbar^{n+1}_i$ the average filtered velocity 
in control volume $\Omega_i$, while $\qbar^{n+1}_{i,j}$ is the auxiliary pressure
at the centroid of face $j$ normalized by the volume of $\Omega_i$.

Next, we rewrite eq.~\eqref{eq:evolveFV-2.1_disc} in semi-discretized form, i.e.~with 
the pressure term in continuous form while all the other terms
are in discrete form, take its divergence and use eq.~\eqref{eq:filter-1.2} to get:
\begin{align}\label{eq:q_poissonNL_p1}
\Delta \qbar^{n+1} = \nabla \cdot \overline{\textbf{H}}(\vbar^{n+1}), ~\text{with} \quad\overline{\textbf{H}}(\vbar^{n+1}) =  \sum_j^{} \mubar_j(\nabla\vbar^{n+1}_i)_j \cdot \textbf{A}_j + \dfrac{\rho}{\Delta t}\v^{n+1}.
\end{align}
%
By integrating eq.~\eqref{eq:q_poissonNL_p1} over the control volume $\Omega_i$ and by applying the Gauss-divergence theorem, we obtain
\begin{align}\label{eq:q_poissonNL_p2}
\int_{\partial \Omega_i} \nabla \qbar^{n+1} \cdot d\textbf{A} = \int_{\partial \Omega_i} \overline{\textbf{H}}(\vbar^{n+1}) \cdot d\textbf{A}.
\end{align}
Finally, we divide by the volume of $\Omega_i$ to get
\begin{align}\label{eq:q_poissonNL_iii}
\sum_j^{} (\nabla \qbar^{n+1}_i)_j \cdot \textbf{A}_j = \sum_j^{} (\overline{\textbf{H}}(\vbar^{n+1}_i))_j \cdot \textbf{A}_j. 
\end{align}


The fully discrete problem associated to the filter is given by \eqref{eq:evolveFV-2.1_disc}, \eqref{eq:q_poissonNL_iii}. 
Also for this problem, we choose a partitioned algorithm:
a slightly modified version of the SIMPLE algorithm, called SIMPLEC algorithm \cite{Doormaal1984},
that features improved and accelerated convergence towards a steady state solution.
We found that the SIMPLEC algorithm is a necessary choice for the filter problem in the Leray-$\alpha$ model,
because the standard SIMPLE algorithm does not converge. On the other hand, 
the SIMPLE algorithm does converge for the filter problem in the Leray-$\alpha$-NL model.  


For the Leray-$\alpha$-NL model, we also need to approximate 
the solution to problem \eqref{eq:vtilde}. 
After using Gauss-divergence theorem, the integral form of eq.~\eqref{eq:vtilde} reads:
\begin{align}\label{eq:Helmholz_int}
\int_{\Omega_i}  \tilde{\v} ^{n+1} d \Omega - \alpha^2\int_{\partial \Omega_i}\nabla \tilde{\v} ^{n+1} \cdot d\textbf{A} = \int_{\Omega_i}\v^{n+1} d \Omega.
\end{align} 
We approximate the diffusion term as:
\begin{align}
\int_{\partial \Omega_i}\nabla \tilde{\v} ^{n+1} \cdot d\textbf{A} \approx \sum_j^{} (\nabla \tilde{\v} ^{n+1})_j \cdot \textbf{A}_j, \el
\end{align} 
where $(\nabla \tilde{\v} ^{n+1})_j$ is treated in the same way as 
$(\nabla\v^{n+1})_j$ and $(\nabla\vbar^{n+1})_j$. Once we divide by the volume of $\Omega_i$, the space
discretized form of eq.~\eqref{eq:Helmholz_int} is 
\begin{align}\label{eq:disc_F}
\tilde{\v} ^{n+1}_i - \alpha^2 \sum_j^{} (\nabla \tilde{\v} ^{n+1}_i)_j \cdot \textbf{A}_j = \v^{n+1}_i,
\end{align} 
where $\tilde{\v} ^{n+1}_i$ is the average value of $\tilde{\v} ^{n+1}$ in control volume $\Omega_i$. 




\subsection{Setting of the relaxation parameter $\chi$}\label{subsec:chi}


In \cite{BQV}, a heuristic formula for estimating the value of $\chi$ is provided. The proposed idea is to 
set $\chi$ such that the total viscous stress on an under-refined mesh of size $h$ provides the 
same amount of dissipation of the viscous stress on a fully refined mesh of size $\eta$. This leads to:
\begin{equation}\label{eq:old_chi}
\chi_1 = \dfrac{ 2\mu }{3 \rho ||a||_\infty \eta \alpha^2} \text{max}\left(h - \eta, 0 \right) \Delta t.
\end{equation}

We propose here an alternative approach. We set $\chi$ such that the viscous contribution on a fully refined 
mesh of size $\eta$ is equivalent to the convex linear combination of the viscous contribution on an under-resolved mesh of size $h$ and
the viscous contribution due to the artificial viscosity introduced by the Leray model: 
\begin{equation}\label{eq:new_chi}
\mu \nabla_\eta \u_h^{n+1} \simeq \chi_2 \left(\mu \nabla_h \u_h^{n+1}\right) + \left(1-\chi_2\right)\left(\nabla_h\cdot\left( \mubar \nabla \u_h^{n+1} \right)\right).
\end{equation}
Eq. \eqref{eq:new_chi} can be interpreted as a relaxation step for the velocity gradient. With the approximation $\nabla_\xi \approx \xi^{-1}$, we obtain
\begin{equation}
\chi_2 = \dfrac{h - \eta}{||a||_\infty \dfrac{\rho \alpha^2}{\mu \Delta t} \eta - \eta}. \el
\end{equation}

In practice, we may set $||a||_\infty = 1$, which is a good approximation on coarse meshes, and obtain:
\begin{equation}\label{eq:chi_comp}
\chi_1 = \dfrac{2\mu}{3 \rho \eta \alpha^2} \text{max}\left(h - \eta, 0 \right) \Delta t, \quad
\chi_2 = \dfrac{h - \eta}{\dfrac{\rho \alpha^2}{\mu \Delta t} \eta - \eta}. 
\end{equation}
We observe that both estimates of $\chi$ defined in eq.~\eqref{eq:chi_comp} are clearly positive. In fact, they would be negative when $h<\eta$, which means that the mesh
is refined enough for a DNS and  the filtering step is not needed. 

\section{Numerical results}\label{sec:num_res}

In this section, we present several numerical results for the NSE (i.e.~no turbulence model), Leray-$\alpha$ and Leray-$\alpha$ NL models. 
Two test cases are investigated: 2D flow past a cylinder and a 3D benchmark from the FDA. 

The number of PISO loops and non-orthogonal correctors has been fixed to 2 for all the simulations.
The following solvers have provided a good compromise between stability, accuracy, and numerical cost. 
The linear algebraic system associated with eq.~\eqref{eq:evolveFV-1.1_disc} is
solved using an iterative solver with symmetric Gauss-Seidel smoother. 
For eq.~\eqref{eq:evolveFV-2.1_disc} and \eqref{eq:disc_F}, we use the
Diagonal Incomplete Cholesky Preconditioned Conjugate Gradient.
Finally, for Poisson problems \eqref{eq:q_poissonL_iii} and \eqref{eq:q_poissonNL_iii} we
use Geometric Agglomerated Algebraic Multigrid Solver GAMG with the Gauss-Seidel smoother. 
The required accuracy is 1e-6 at each time step. 

\subsection{2D benchmark: channel flow around a cylinder}\label{sec:cylinder}
\begin{figure}[h]
\centering
\subfloat[][mesh $16k_H$]{\includegraphics[width=0.45\textwidth]{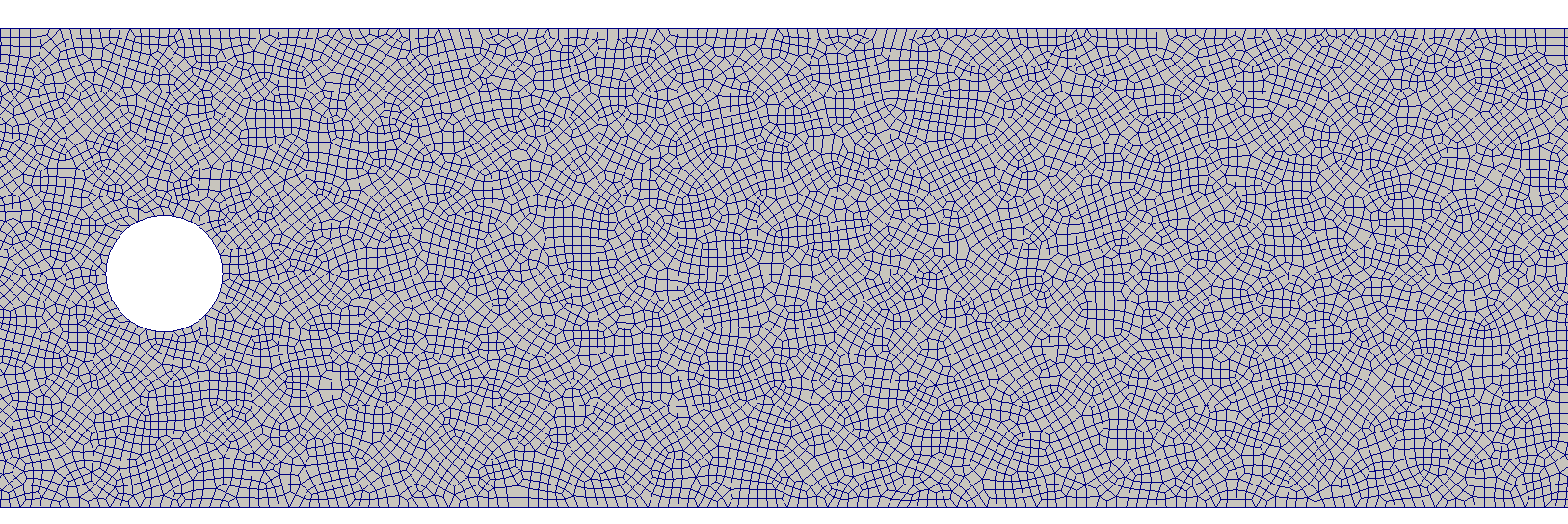}}~
\subfloat[][velocity magnitude at $t = 8$]{\includegraphics[width=0.45\textwidth]{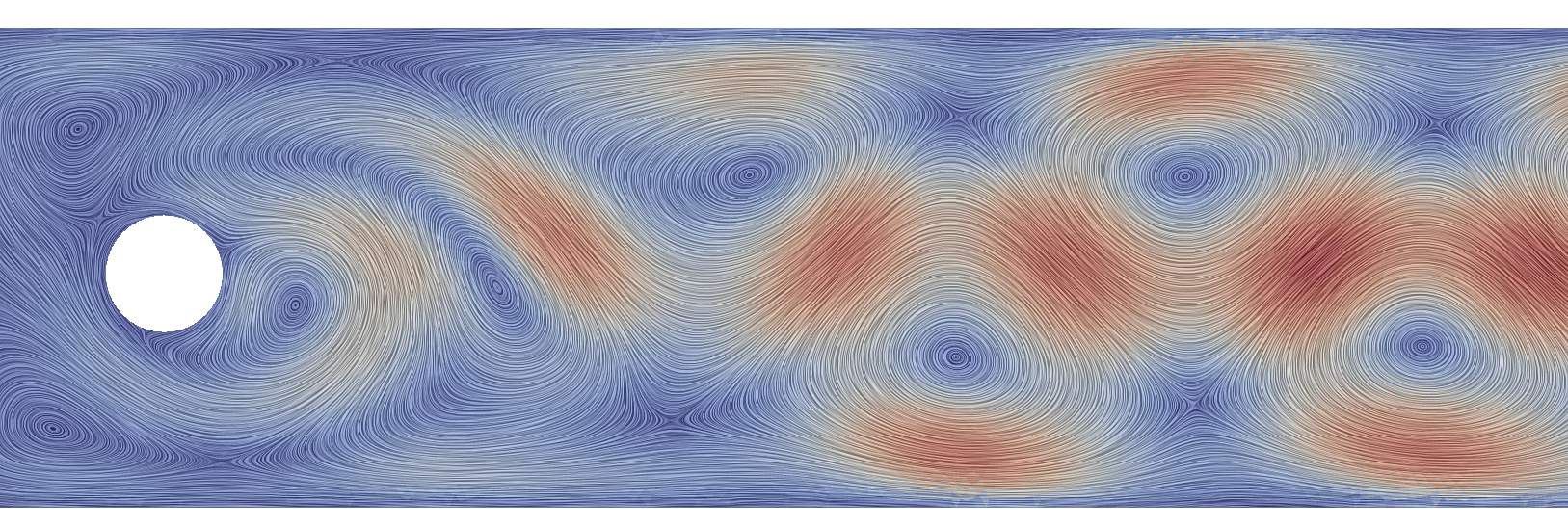}}
\caption{(a) Part of mesh $16k_H$  and (b) velocity magnitude at $t = 8$ computed by the NSE on mesh $16k_H$.
The velocity magnitude goes from 0 (blue) to 3.9e-1 (red).}
\label{fig:example_cyl}
\end{figure}

The first test we consider is a well-known benchmark \cite{John2004,Turek1996}. 
The computational domain is a 2.2 $\times$ 0.41 rectangular channel with a cylinder of radius 0.05 centered at (0.2, 0.2), 
when taking the bottom left corner of the channel as the origin of the axes. 
Fig.~\ref{fig:example_cyl} shows part of the computational domain. 
We impose a no slip boundary condition on the upper and lower wall and on the cylinder. At the inflow and the outflow we prescribe
the following velocity profile:
\begin{align}\label{eq:cyl_bc}
\v(0,y,t) = \left(\dfrac{6}{0.41^2} \sin\left(\pi t/8 \right) y \left(0.41 - y \right), 0\right), \quad y \in [0, 2.2], \quad t \in (0, 8].
\end{align}
The partitioned algorithms we use (see Sec.~\ref{subsec:space-discrete-evolve}
and \ref{subsec:space-discrete-filter}) require a boundary condition for the pressure too. We choose ${\partial q}/{\partial \n} = {\partial \qbar}/{\partial \n} = 0$ on each boundary where $\n$ is the outward normal. 
We set density $\rho = 1$ and viscosity $\mu = 10^{-3}$. 
There is no external force, so $\f = \mathbf{0}$. 
We start the simulations from fluid at rest.
Note that the Reynolds number is time dependent, with $0 \leq Re \leq 100$ \cite{Turek1996}. 
The quantities of interest for this benchmark are the drag and lift coefficients:
\begin{align}\label{eq:cd_cl}
c_d(t) = \dfrac{2}{\rho L_{r} {U}^2_{r}} \int_S \left(\boldsymbol{\sigma}
\cdot \boldsymbol{n}\right) \cdot \boldsymbol{t}~dS, \quad
c_l(t) = \dfrac{2}{\rho L_{r} {U}^2_{r}} \int_S \left(\boldsymbol{\sigma}
\cdot \boldsymbol{n}\right) \cdot \boldsymbol{n}~dS,
\end{align}
where $U_{r}= 1$ is the maximum velocity at the inlet/outlet, $L_r = 0.1$ is the diameter of the cylinder, 
$S$ is the surface of the cylinder, and $\boldsymbol{t}$ and $\boldsymbol{n}$ are the tangential and normal 
unit vectors, respectively.

We consider several meshes with prismatic and hexahedral elements
generated using \emph{Gmsh} \cite{Gmsh}.
Table \ref{tab:1} reports mesh name, minimum and maximum diameter, and number of
cells for each mesh. 
The name of each grid refers to the number of cells and the subscript denotes the kind of element,
prismatic ($P$) or hexahedral ($H$). 
See Fig.~\ref{fig:example_cyl} (a) shows part mesh $16k_H$. 
The quality of all the meshes under consideration is high.
Hexahedral and prismatic meshes have very low values of maximum non-orthogonality 
(33$^\circ$ to 39$^\circ$), average non-orthogonality (4$^\circ$ to 7$^\circ$), 
skewness (0.5 to 1.2), and maximum aspect ratio (up to 2). 
Time step used for the simulations was governed by the maximum Courant-Friedrichs-Lewy number ($CFL_{max}$), 
set to $CFL_{max} = 0.2$. For all the simulations, we use the PISO algorithm with
flux correction term \cite{ddtPhiCorr} and  second-order accurate Linear Upwind Differencing (LUD) \cite{Warming1976} scheme for the convective term.

\begin{table}[h]
\centering
\begin{tabular}{cccc|cccc}
\multicolumn{2}{c}{} \\
\cline{1-8}
mesh name    & $h_{min}$ & $h_{max}$ & No. of cells &  mesh name    & $h_{min}$ & $h_{max}$ & No. of cells \\
\hline
$16k_P$      & 7.4e-3  & 1.5e-2   & 1.6e4   & $16k_H$      &   4.2e-3 & 1.1e-2   & 1.59e4 \\
\hline
$25k_P$       &  5.9e-3  &  1.3e-2  &  2.5e4 & $25k_H$     &   3.5e-3  &   9.6e-3 &    2.45e4 \\
\hline
$63k_P$      &   3.7e-3  &   8.1e-3 & 6.33e4  & $63k_H$      &   1.9e-3 & 5.8e-3   &     6.28e4 \\
\hline
$120k_P$      &    2.7e-3 &  6.2e-3  &   1.21e5 & $120k_H$      &    1.4e-3 &  4.2e-3 & 1.22e5    \\
\hline
$200k_P$ &    2e-3  &  4.8e-3 &    1.98e5 & $200k_H$      &    1e-3   & 3.4e-3    & 1.98e5    \\
\hline
\end{tabular}
\caption{Name, minimum diameter $h_{min}$, maximum diameter $h_{max}$, and number of cells
for all the prismatic ($P$) or hexahedral ($H$) and meshes used for the 2D flow past a cylinder.}
\label{tab:1}
\end{table}

\begin{rem}\label{rem:Kolm_cyl}
As discussed in Sec. \ref{sec:leray_model}, DNS needs a mesh with spacing $h\sim\eta$. 
For the current test, the Reynolds number is time dependent. One could either 
calculate the Kolmogorov scale (\ref{eq:eta-re}) based on the maximum Reynolds number, 
which gives $\eta \approx 3.2e-3$, or introduce a mean Kolmogorov scale $\bar{\eta}$:
\begin{align}\label{eq:mean_Kolm}
\bar{\eta} = \dfrac{1}{T} \int_0^T \eta(t) dt \approx 3.7e-3
\end{align}
where $T = 8$ is the final computational time. Notice that $\eta \approx \bar{\eta}$,
thus  a DNS would roughly need a number of cells of the order of 100k.  
\end{rem}

\begin{figure}[h]
\centering
\subfloat[][Lift coefficient]{\includegraphics[width=0.5\textwidth]{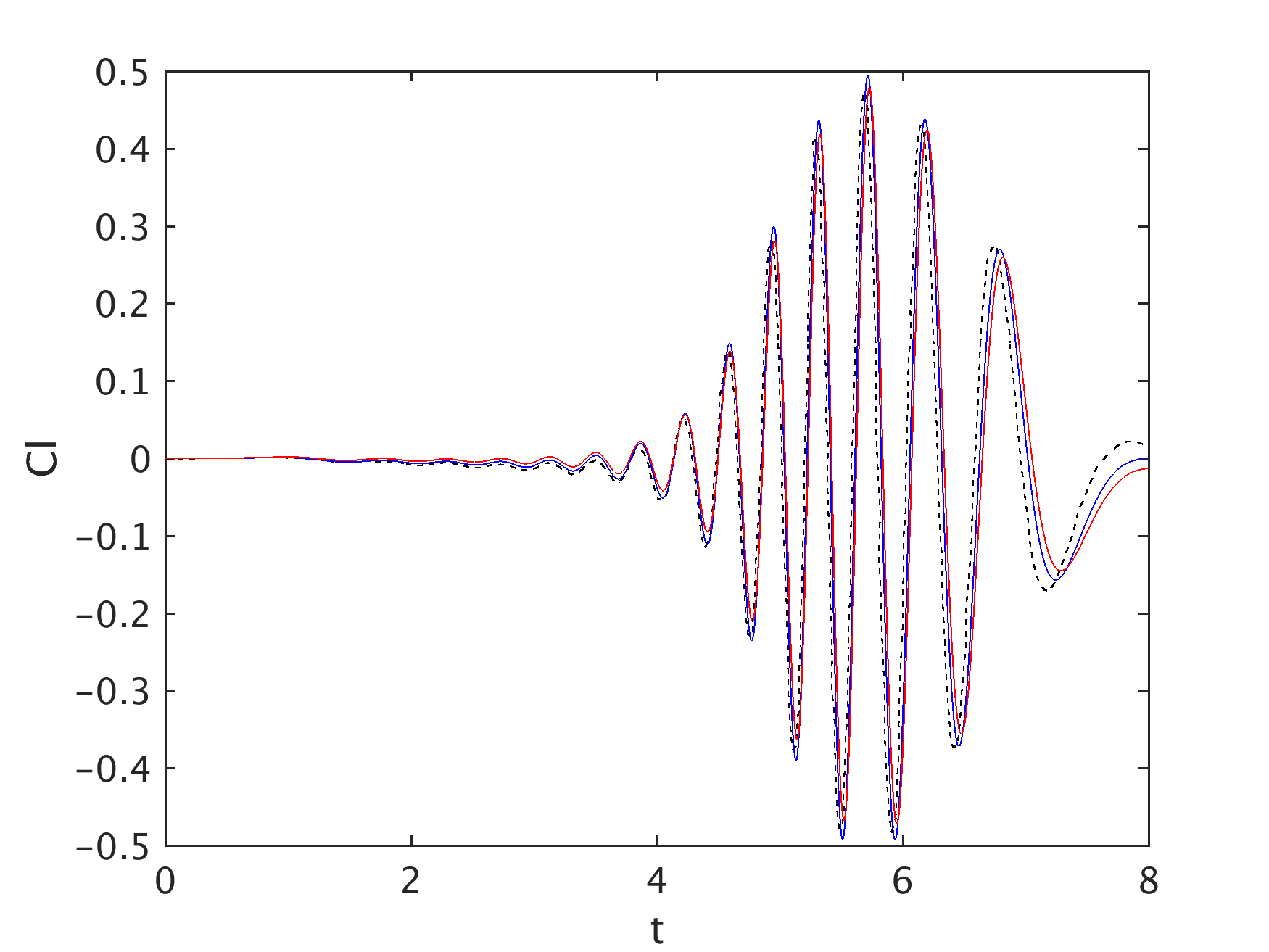}}
\subfloat[][Drag coefficient]{\includegraphics[width=0.5\textwidth]{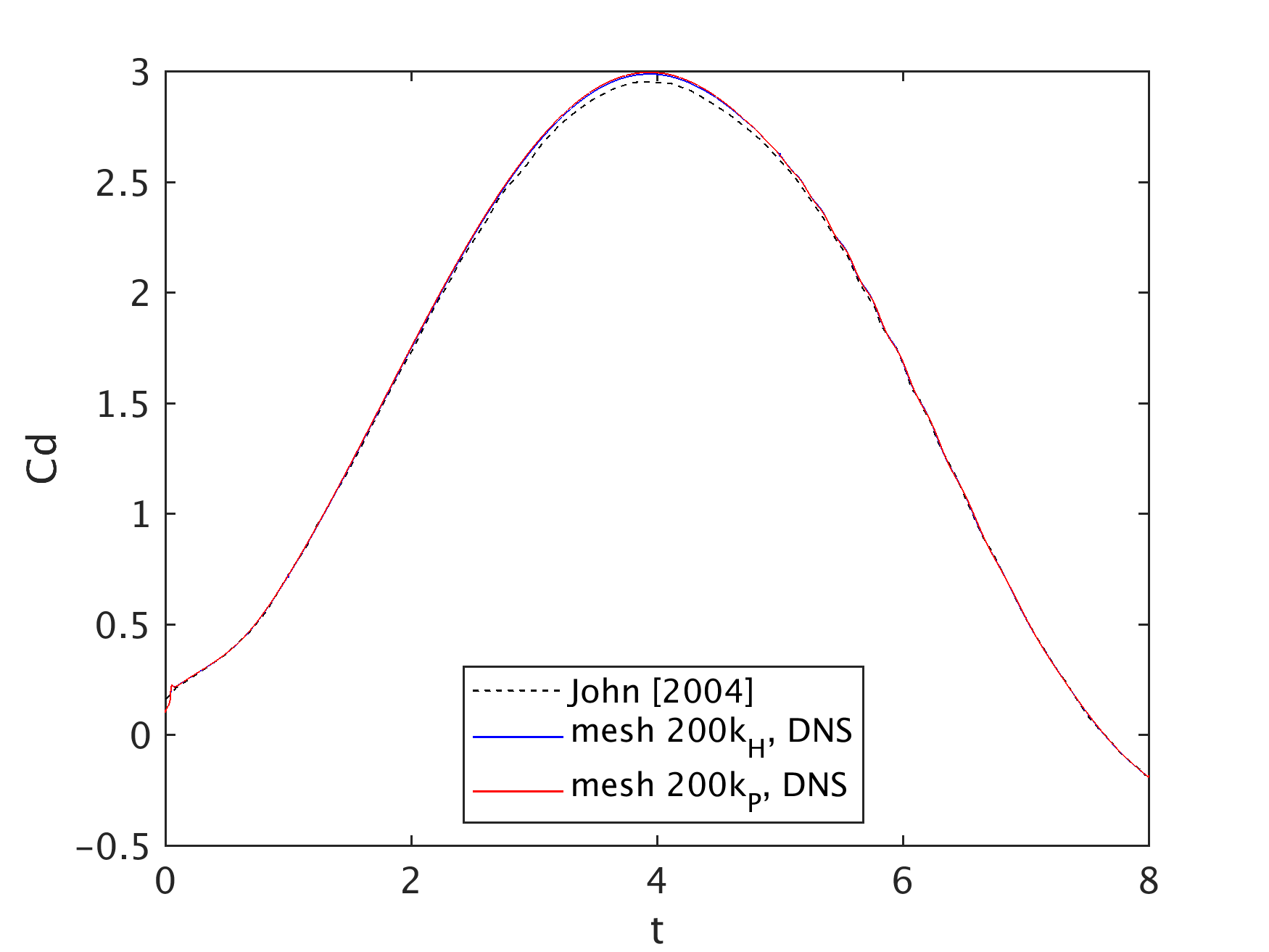}}
\caption{Evolution of lift and drag coefficients given by DNS with meshes $200k_H$ and $200k_P$
compared against the results in \cite{John2004}. The legend in (b) is common to both subfigures.}
\label{fig:DNS}
\end{figure}

Fig.~\ref{fig:DNS} compares the NSE solution obtained with meshes $200k_P$ and $200k_H$ 
with the reference curves reported in \cite{John2004}. We see that the solutions are in good agreement, so a DNS is possible both with meshes $200k_P$ and $200k_H$. Notice that the level of refinement of these meshes is in agreement with that one predicted in Remark \ref{rem:Kolm_cyl}. Hereinafter, we will refer to the solutions computed with a DNS on meshes $200k_P$ and $200k_H$ as the \emph{true} solution and we will refer to 
meshes $200k_P$ and $200k_H$ as DNS meshes.

\begin{figure}[h]
\centering
\subfloat[][Lift coefficient, $P$ meshes]{\includegraphics[width=0.5\textwidth]{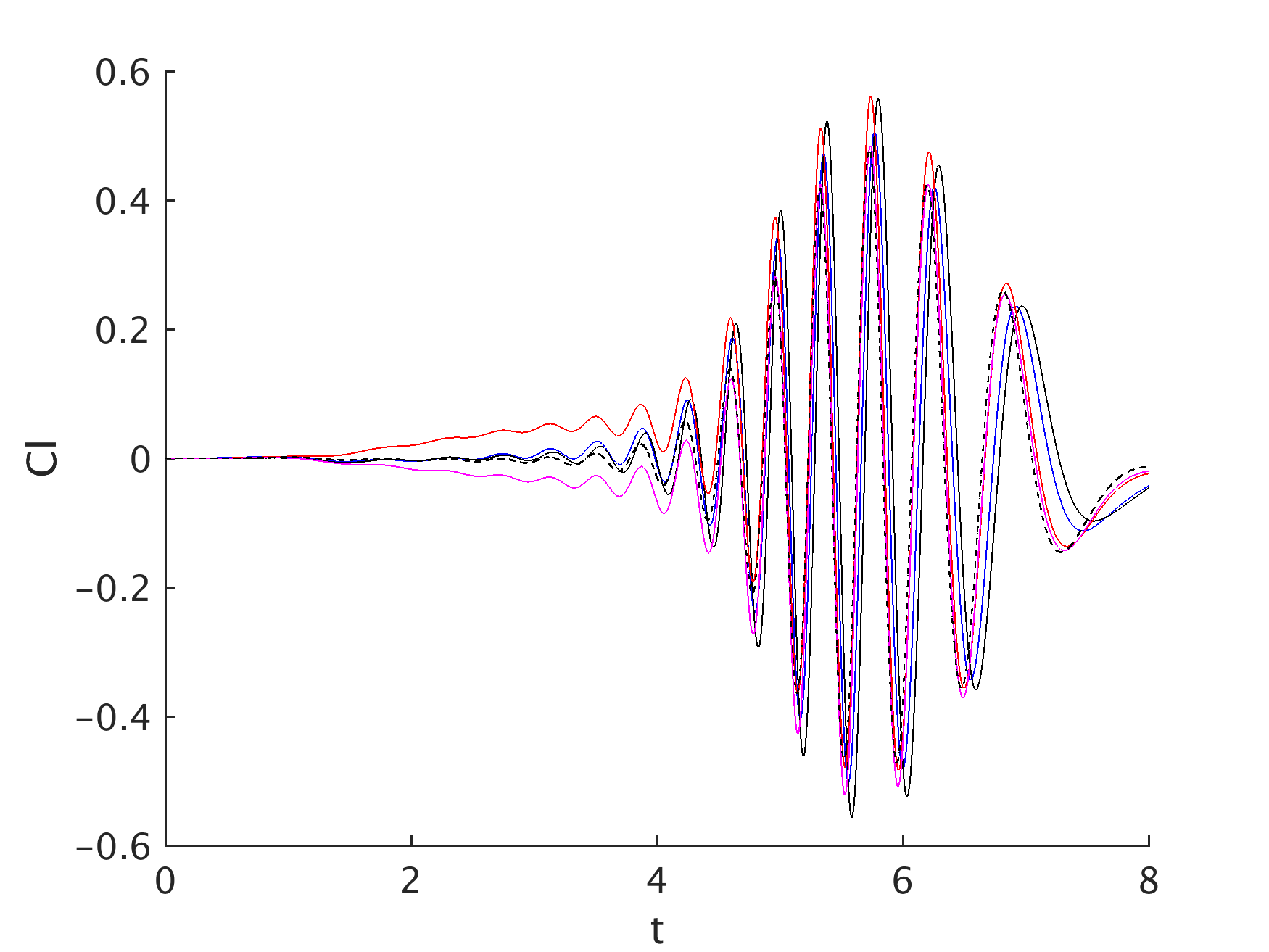}}
\subfloat[][Drag coefficient, $P$ meshes]{\includegraphics[width=0.5\textwidth]{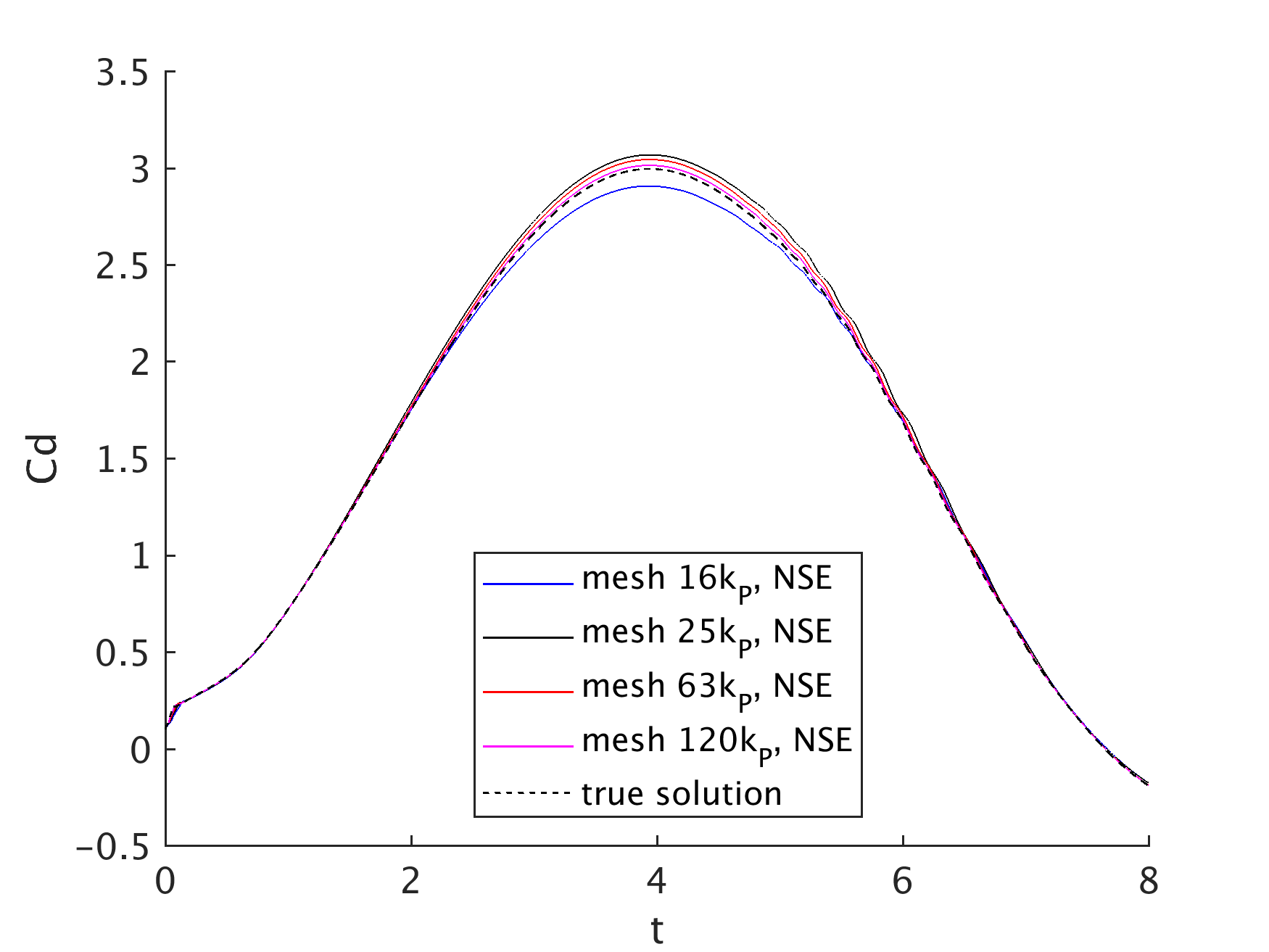}} \\
\subfloat[][Lift coefficient, $H$ meshes]{\includegraphics[width=0.5\textwidth]{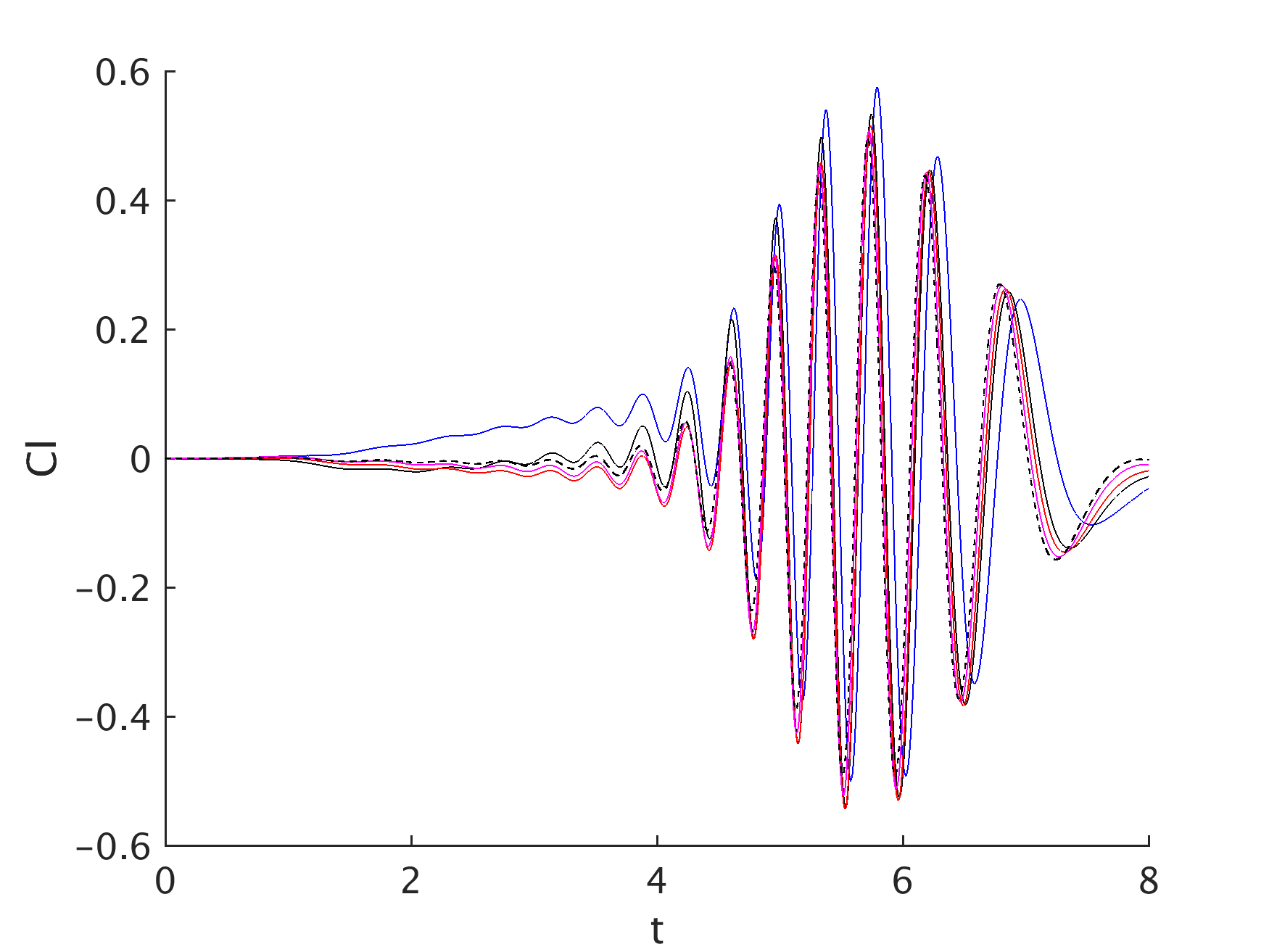}}
\subfloat[][Drag coefficient, $H$ meshes]{\includegraphics[width=0.5\textwidth]{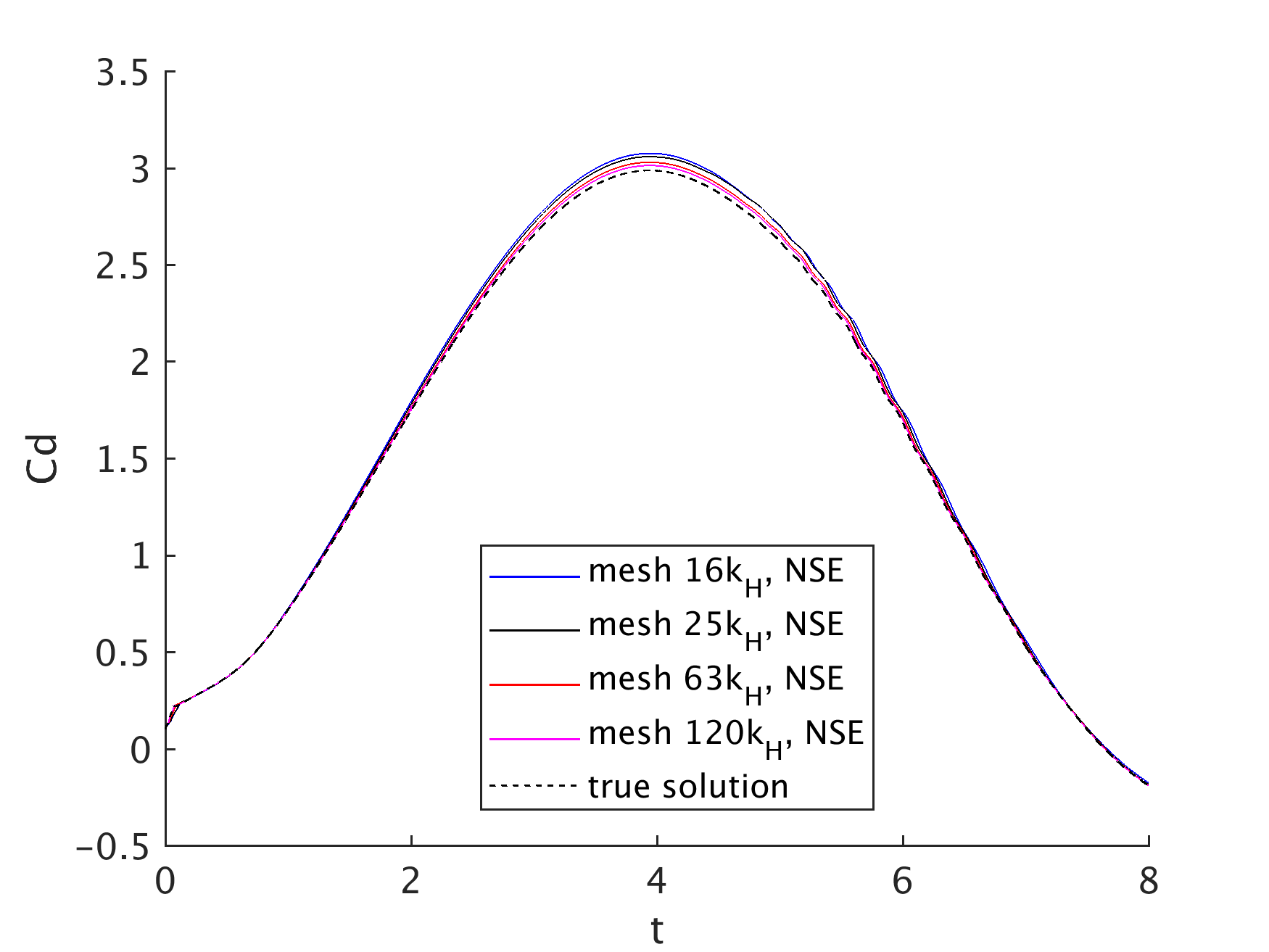}}
\caption{Evolution of lift and drag coefficients computed from the NSE solutions on the prismatic meshes
(top) and all hexahedral meshes (bottom)
reported in Table \ref{tab:1}. The legend in (b) and (d) holds also for (a) and (c).}
\label{fig:NSE_PH}
\end{figure}

Fig.~\ref{fig:NSE_PH} displays the lift and drag coefficients computed from the NSE solution on all the meshes
in Table \ref{tab:1} that are coarser than the DNS meshes. We observe that 
the lift coefficients computed from the NSE solutions on the 
coarser meshes show a large difference with respect to the \emph{true} $c_l$, 
in terms of both amplitude and phase.
To quantitate this difference, Table \ref{tab:2} reports maximum lift and drag coefficients 
obtained from the NSE solutions
and times at which the maxima occur on all the meshes, together with the corresponding
values from \cite{John2004}. We observe a monotonic convergence for the maximum values of 
$c_l$ and $c_d$ computed on hexahedral meshes, while for prismatic meshes the convergence is not monotonic. 
In Table \ref{tab:2}, our results computed with hexahedral (resp., prismatic) 
meshes are compared with the results in \cite{John2004} computed with $Q_2$-$P_1^{d}$ (resp., $P_2$-$P_1$) 
finite elements and the Crank-Nicolson
scheme for time discretization. Notice that the two sets of value from 
\cite{John2004} reported in Table \ref{tab:2} coincide up to the 
reported number of digits.

\begin{table}[h]
\centering
\begin{tabular}{ccccc|ccccc}
\multicolumn{2}{c}{} \\
\cline{1-10}
Mesh & $t(c_{l,max})$ & $c_{l,max}$ & $t(c_{d,max})$ & $c_{d,max}$ &Mesh & $t(c_{l,max})$ & $c_{l,max}$ & $t(c_{d,max})$ & $c_{d,max}$ \\
\hline
$16k_H$ & 5.79   & 0.574  & 3.947  & 3.074 &$16k_P$ & 5.769 & 0.505  & 3.928 & 2.906  \\
\hline
$25k_H$ &  5.743   & 0.533  &  3.929 & 3.057  &$25k_P$ &  5.80   & 0.558  & 3.946  & 3.066   \\
\hline
$63k_H$ &  5.739 & 0.515 &  3.945 & 3.028  & $63k_P$ &  5.737 & 0.561    &  3.943 &  3.043   \\
\hline
$120k_H$ &  5.724 & 0.507 &  3.936 & 3.012 & $120k_P$ &  5.737 &  0.484  &  3.935 & 3.012   \\
\hline
$200k_H$ &   5.716  & 0.495  &  3.934 & 2.986 & $200k_P$ &  5.726   & 0.478    & 3.939 & 2.994  \\
\hline
$Q_2$-$P_1^{d}$  & 5.694 & 0.478 & 3.936 & 2.951 & $P_2$-$P_1$ & 5.694 & 0.478 & 3.936 & 2.951  \\
\hline
\end{tabular}
\caption{Maximum lift and drag coefficients given by NSE
and times at which the maxima occur on all the meshes
under consideration. 
The bottom row reports the results from \cite{John2004}.
}
\label{tab:2}
\end{table}

\begin{figure}[h]
\centering
\subfloat[][Maximum lift coefficient errors]{\includegraphics[width=0.5\textwidth]{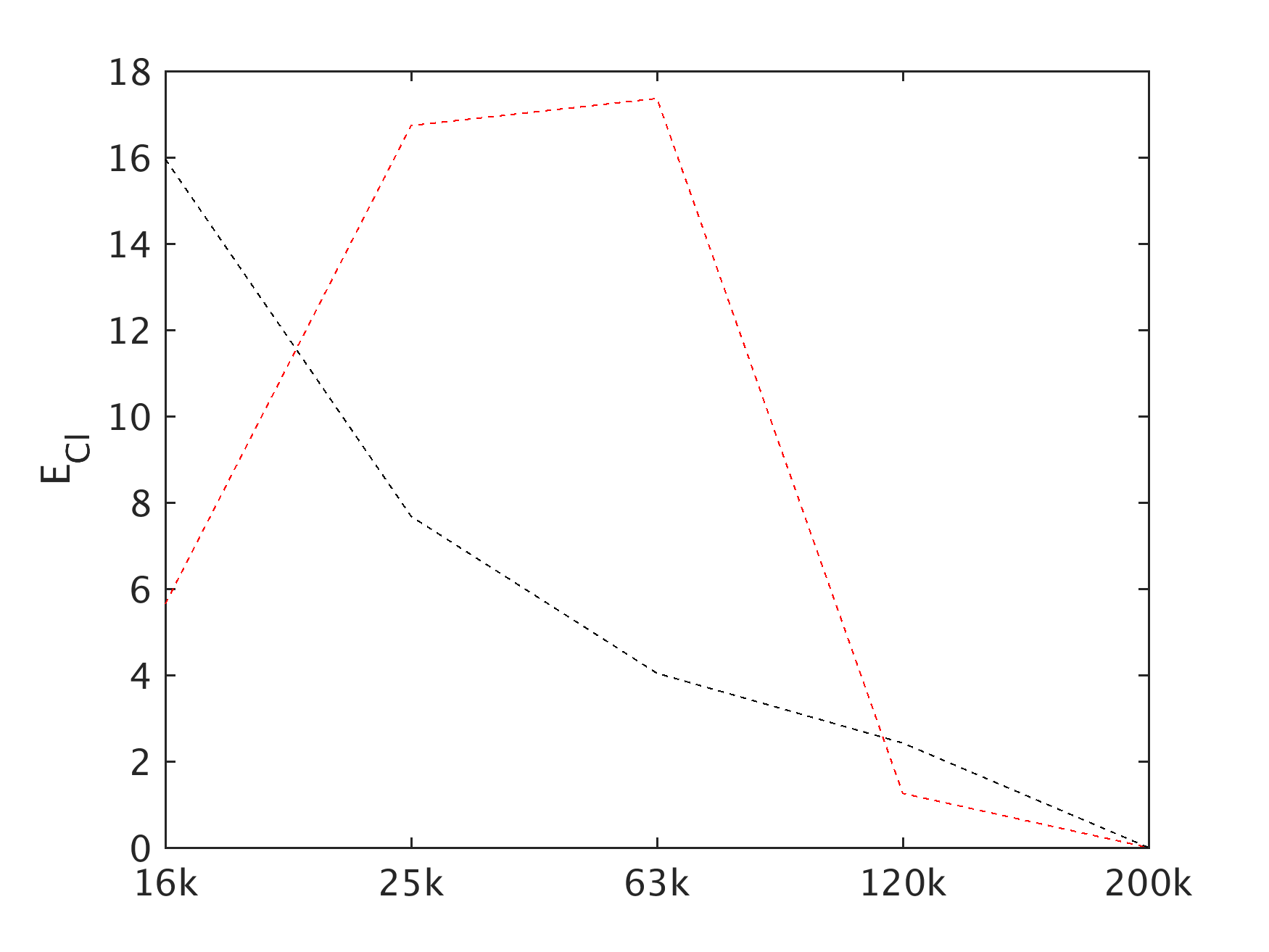}}
\subfloat[][Maximum drag coefficient errors]{\includegraphics[width=0.5\textwidth]{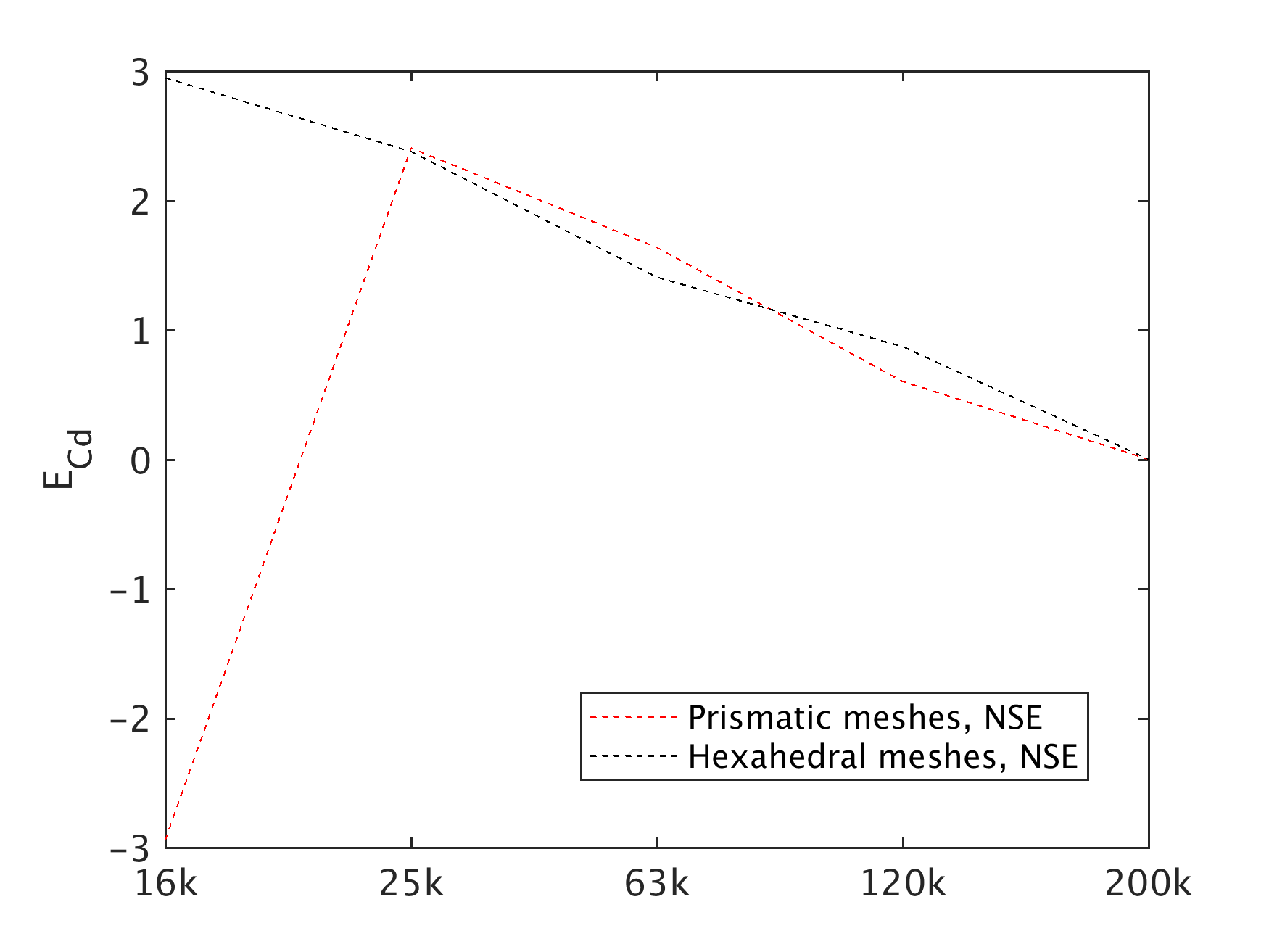}}
\caption{Maximum (a) lift and (b) drag coefficients errors as defined in \eqref{eq:E_c_1} 
for the different meshes under consideration. The legend in (b) is common to both subfigures.}
\label{fig:errs}
\end{figure}

For all simulations, we evaluated also the error for the maximum values of aerodynamic coefficients
compared to the \emph{true} values:

\begin{align}\label{eq:E_c_1}
E_{c_l} = \dfrac{c_{l,max} - c_{l,max}^{{true}}}{c_{l,max}^{{true}}} \cdot 100, \quad
E_{c_d} = \dfrac{c_{d,max} - c_{d,max}^{{true}}}{c_{d,max}^{{true}}} \cdot 100.
\end{align}
Fig.~\ref{fig:errs} confirms what already noted from Table \ref{tab:2}: for hexahedral meshes the errors 
\eqref{eq:E_c_1} decrease monotonically as the mesh
is refined, while that is not the case for prismatic meshes. In particular, we notice that the $c_{l,max}$ computed 
with mesh $16k_P$ is surprisingly good, while it worsens for meshes $25k_P$
and $63k_P$. It is only with mesh  $120k_P$ that we obtain a better $c_{l,max}$
than the one computed with mesh  $16k_P$.

Now that we have identified the right level of refinement for a DNS (i.e.,~ roughly $200k$
cells), we test the EFR algorithm on all meshes in Table \ref{tab:1} 
that have a lower level of refinement. We aim at comparing the results given by the Leray-$\alpha$ 
and Leray-$\alpha$ NL models, and understanding the role of $\chi$ and $\alpha$ for both models. 
We have chosen to set $\chi = 0$, i.e. no relaxation, and $\chi = \Delta t$, as suggested in 
\cite{layton_CMAME} to keep the numerical dissipation low.
When $\chi$ = 0, we will refer to the algorithm as Evolve-Filter (EF),
and in particular EF L (resp., EF NL) if the Leray-$\alpha$ (resp., Leray-$\alpha$ NL) model is used. 
As for $\alpha$, two characteristic values are considered: 
the length $h_{min}$ of the shortest edge in the mesh (see \cite{Bicol_Quaini2018} and Remark 4.3 in \cite{BQV}) 
and the Kolmogorov scale $\eta$ based on the maximum Reynolds number (see Remark \ref{rem:Kolm_cyl}). 
Notice that the following ``coarse'' meshes
have $h_{min}< \eta$: $63k_H$, $120k_P$, $120k_H$. Moreover, meshes 
$25k_H$ and $63k_P$ have $h_{min}$ which is comparable to $\eta$.
Thus, the choice $\alpha = \eta$ is reasonable. For all the different cases under consideration 
we will report $c_l$, which is the more critical coefficient,  
while for $c_d$ we will only report the maximum value in Table \ref{tab:3}.

\begin{figure}
\centering
 \begin{overpic}[width=0.45\textwidth]{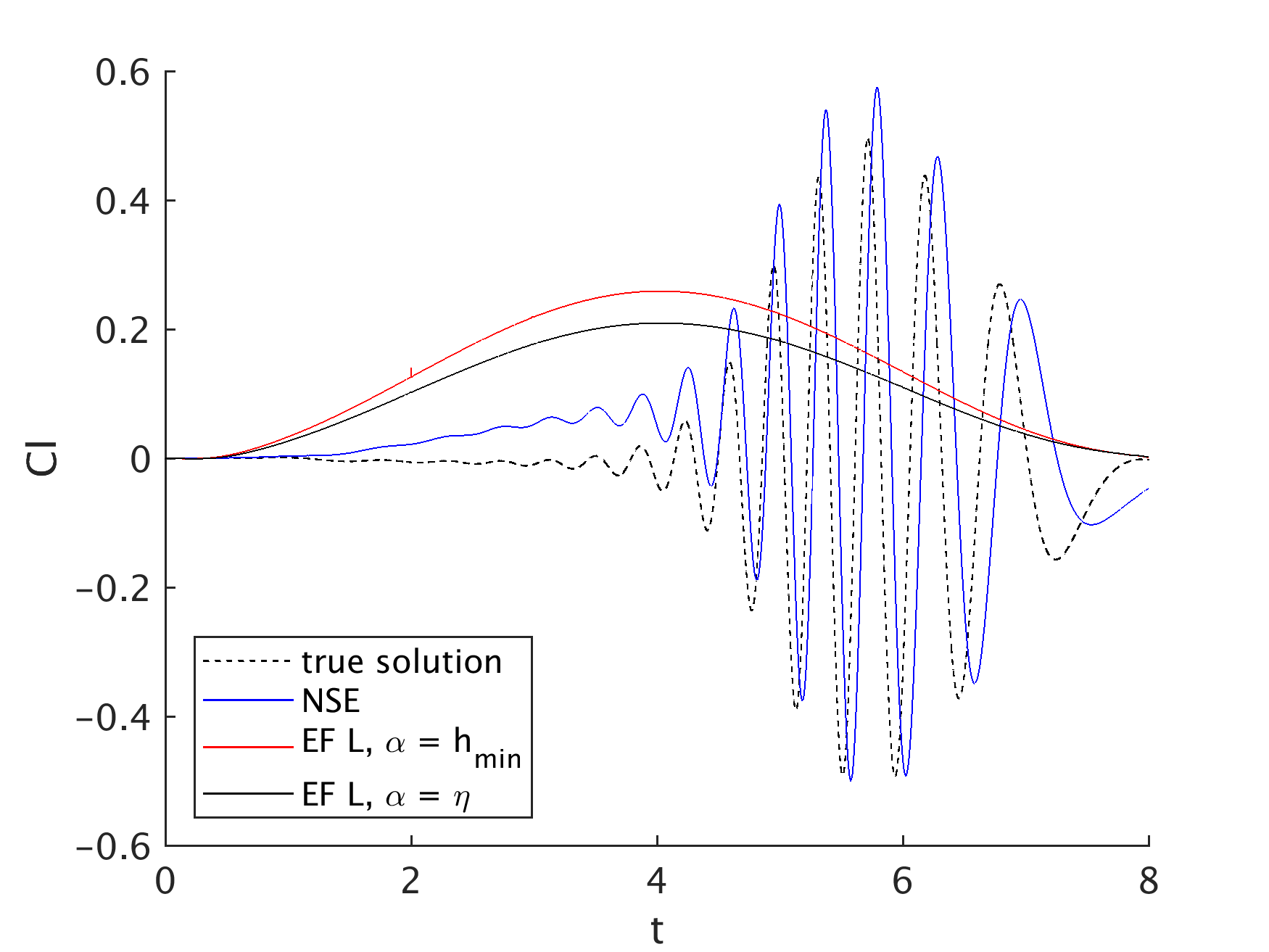}
        \put(32,70){\small{Mesh $16k_H$, EF L}}
      \end{overpic}
 \begin{overpic}[width=0.45\textwidth]{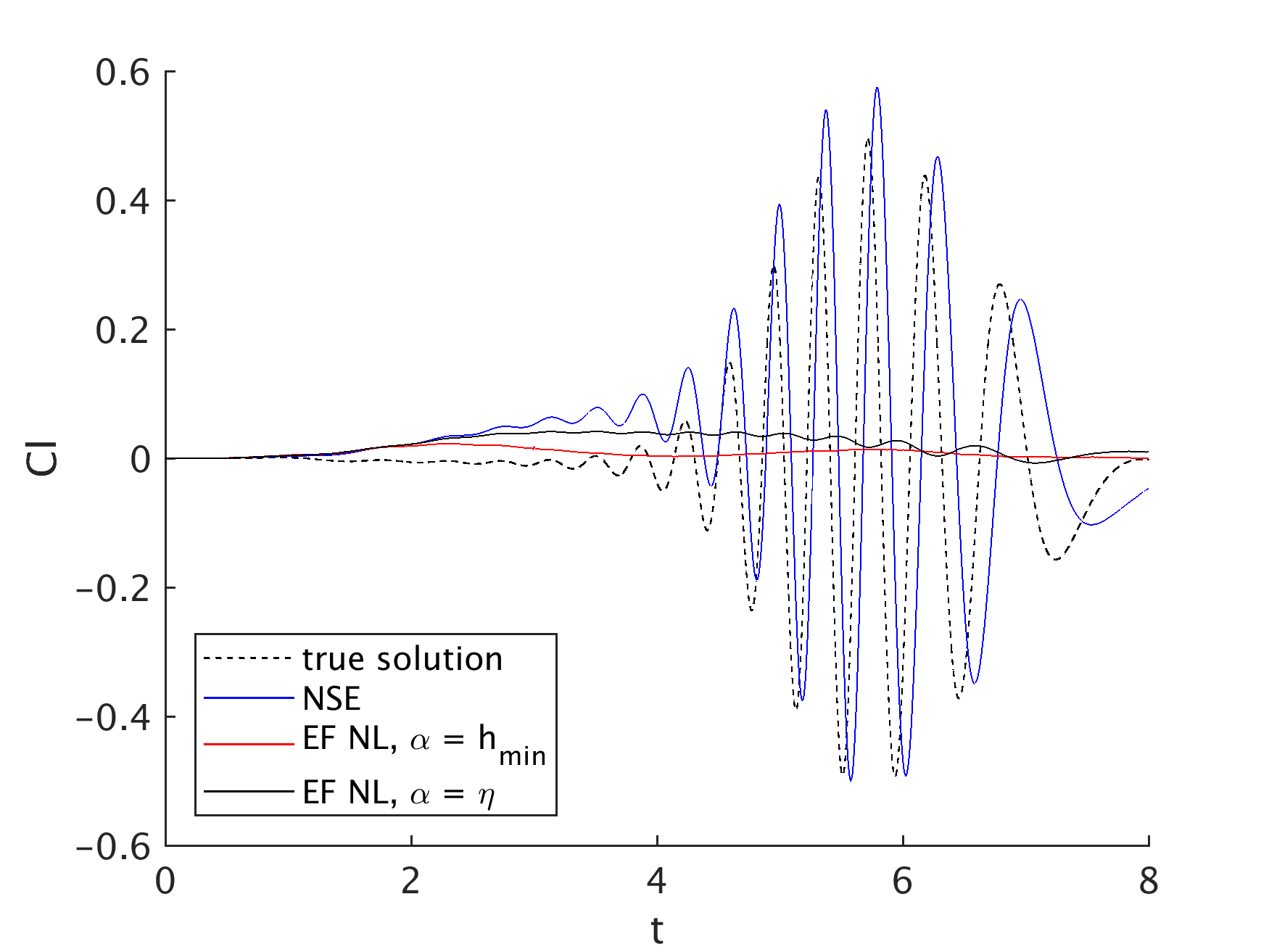}
        \put(32,70){\small{Mesh $16k_H$, EF NL}}
      \end{overpic}\\
       \begin{overpic}[width=0.45\textwidth]{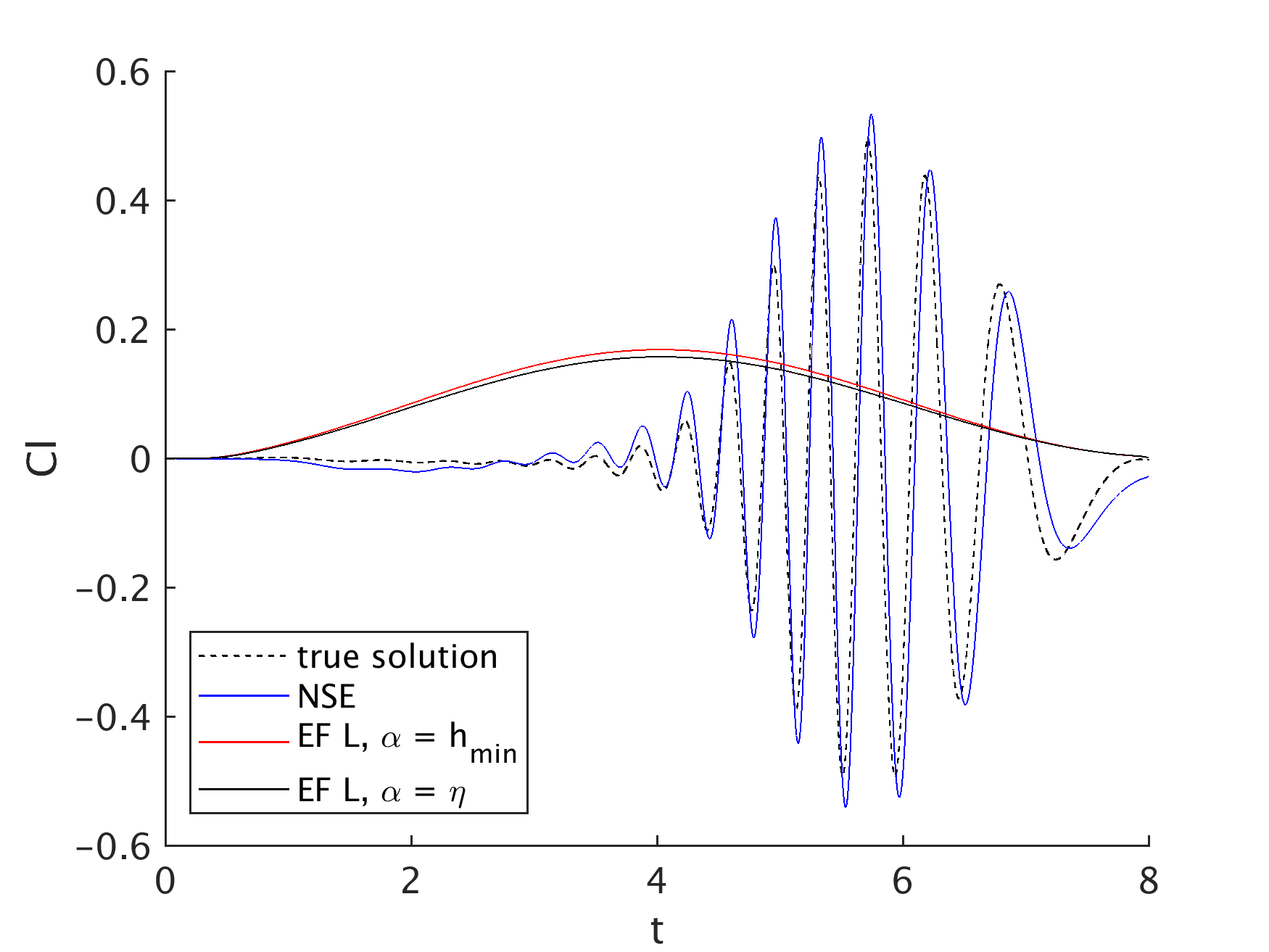}
        \put(32,70){\small{Mesh $25k_H$, EF L}}
      \end{overpic}
 \begin{overpic}[width=0.45\textwidth]{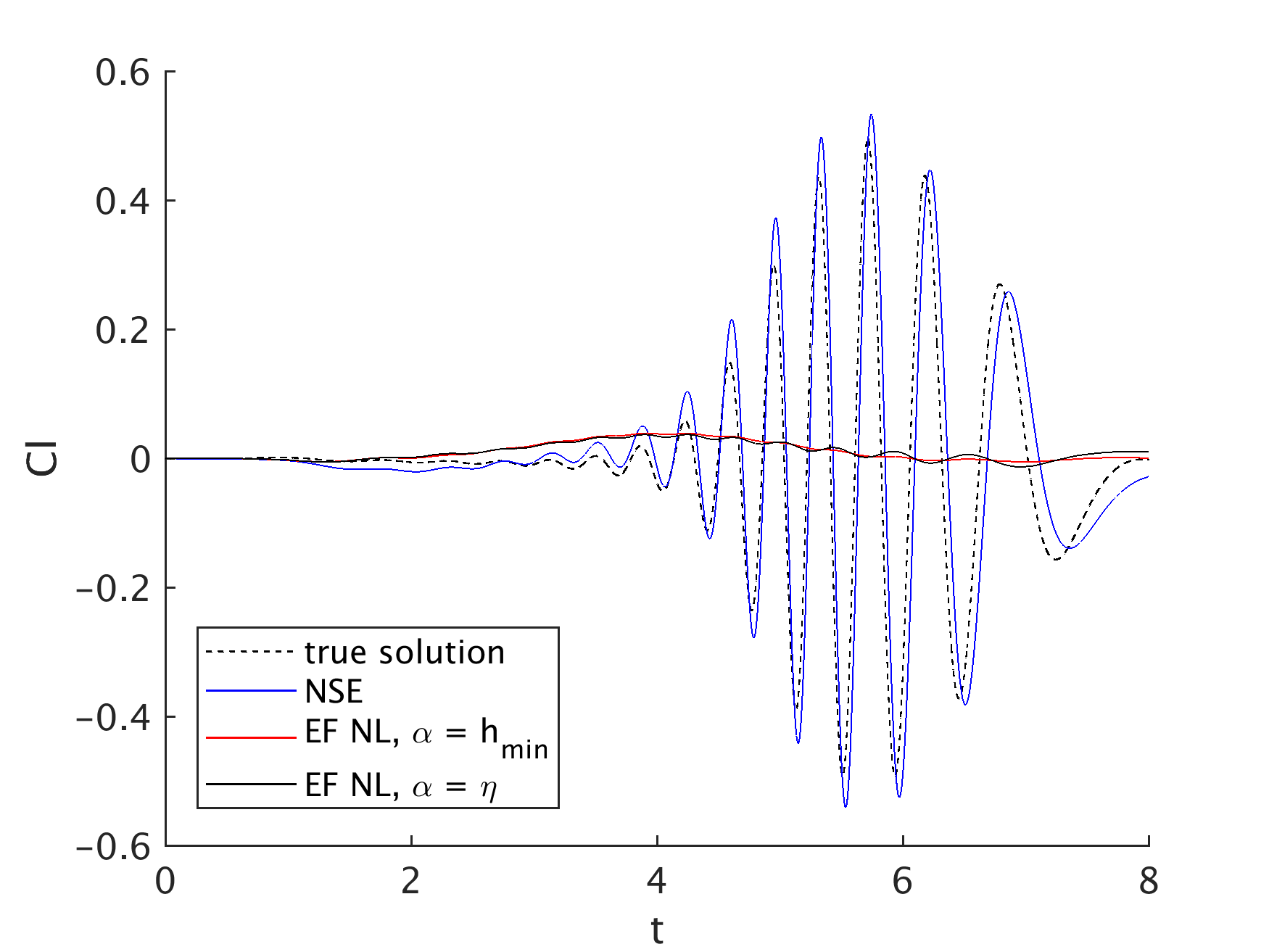}
        \put(32,70){\small{Mesh $25k_H$, EF NL}}
      \end{overpic}\\
       \begin{overpic}[width=0.45\textwidth]{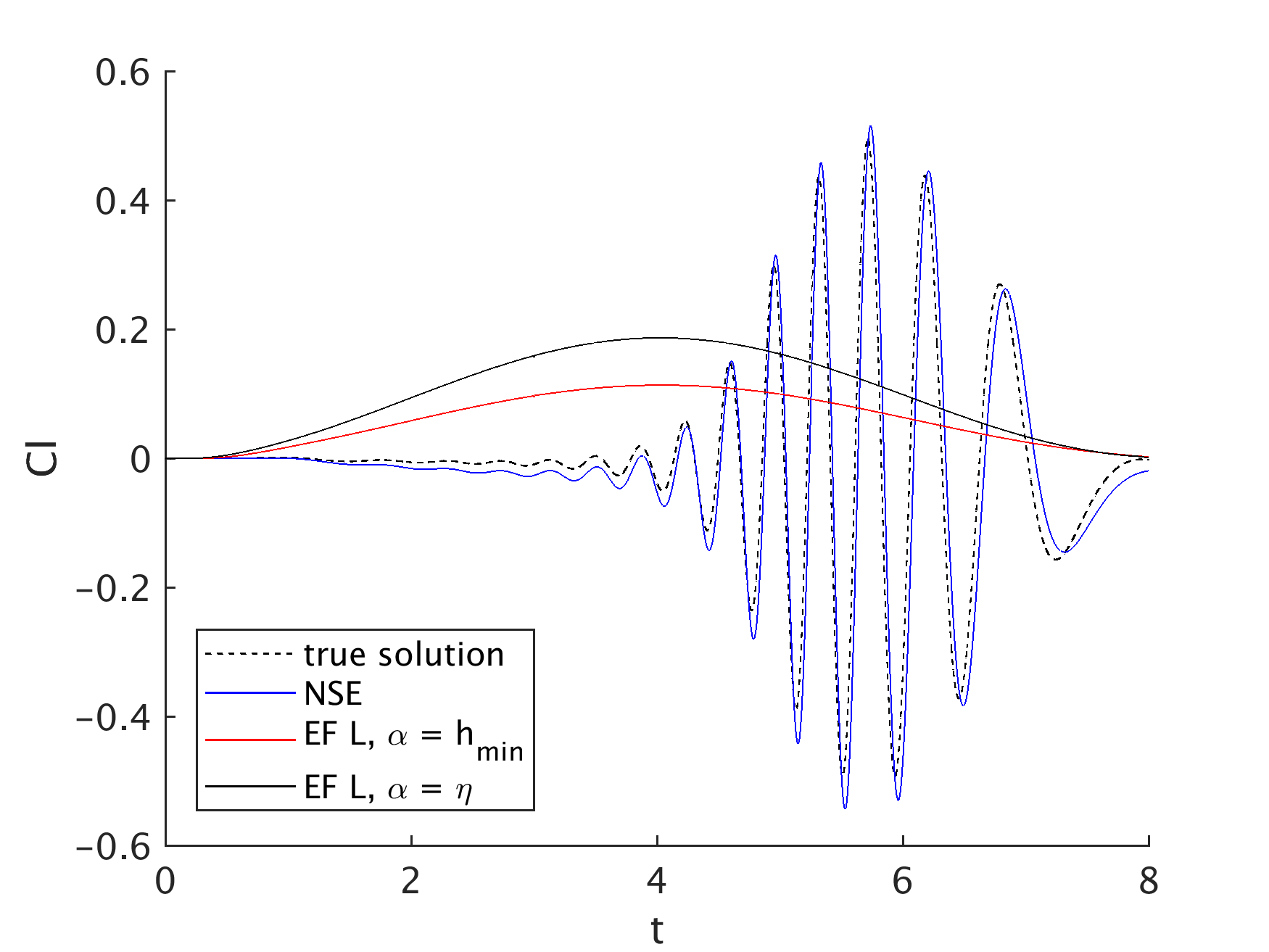}
        \put(32,70){\small{Mesh $63k_H$, EF L}}
      \end{overpic}
 \begin{overpic}[width=0.45\textwidth]{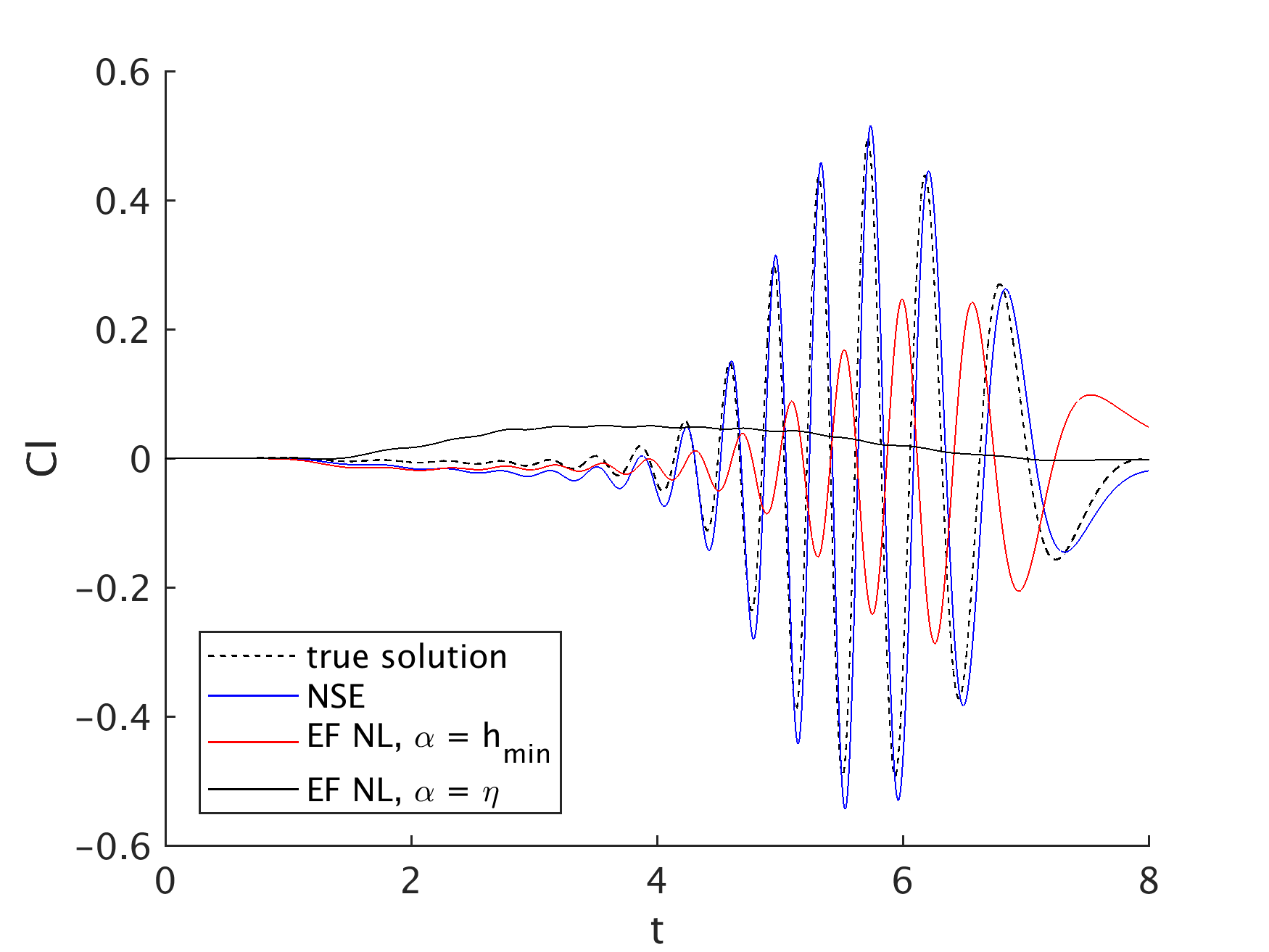}
        \put(32,70){\small{Mesh $63k_H$, EF NL}}
      \end{overpic}\\
       \begin{overpic}[width=0.45\textwidth]{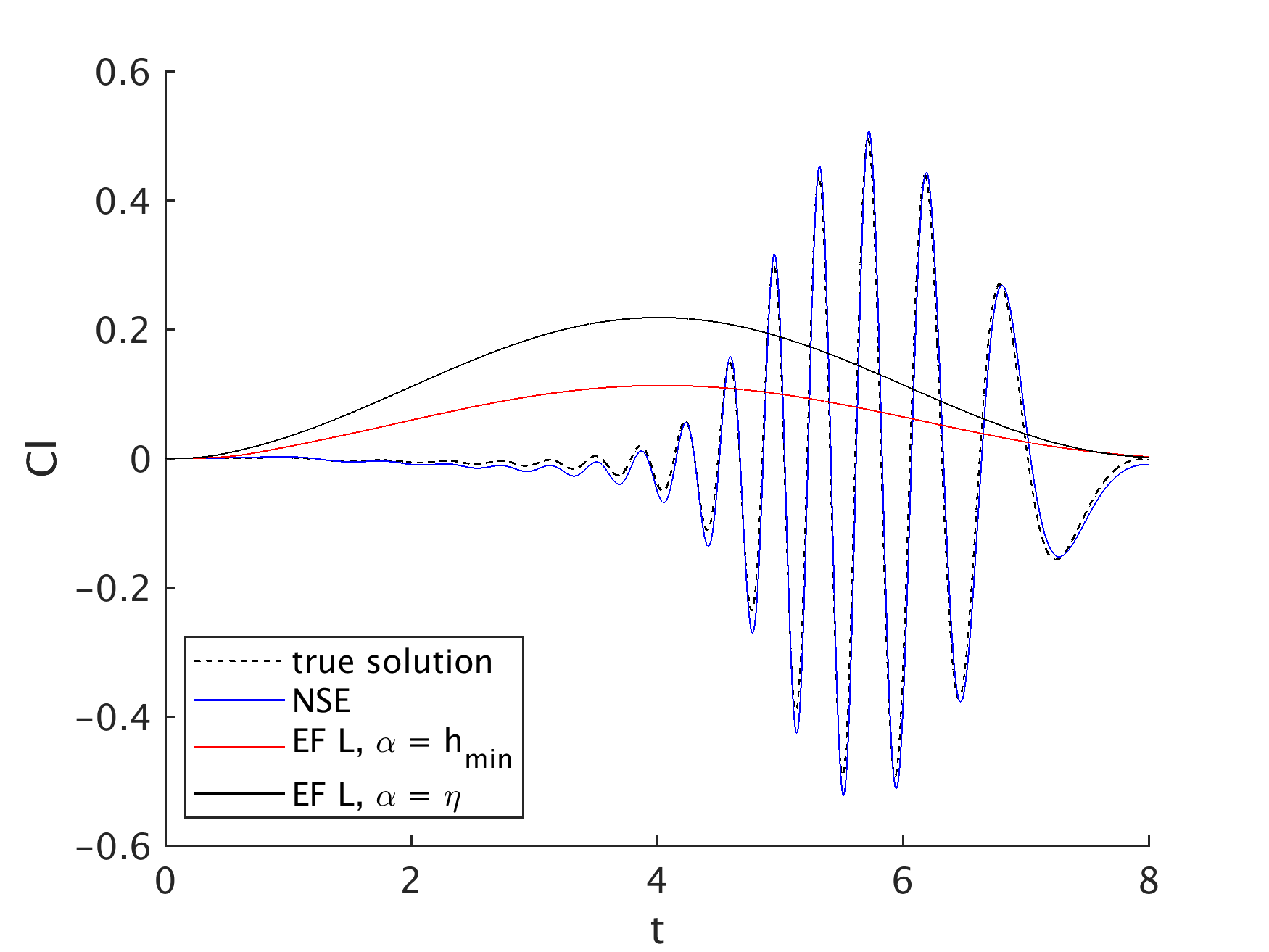}
        \put(32,70){\small{Mesh $120k_H$, EF L}}
      \end{overpic}
 \begin{overpic}[width=0.45\textwidth]{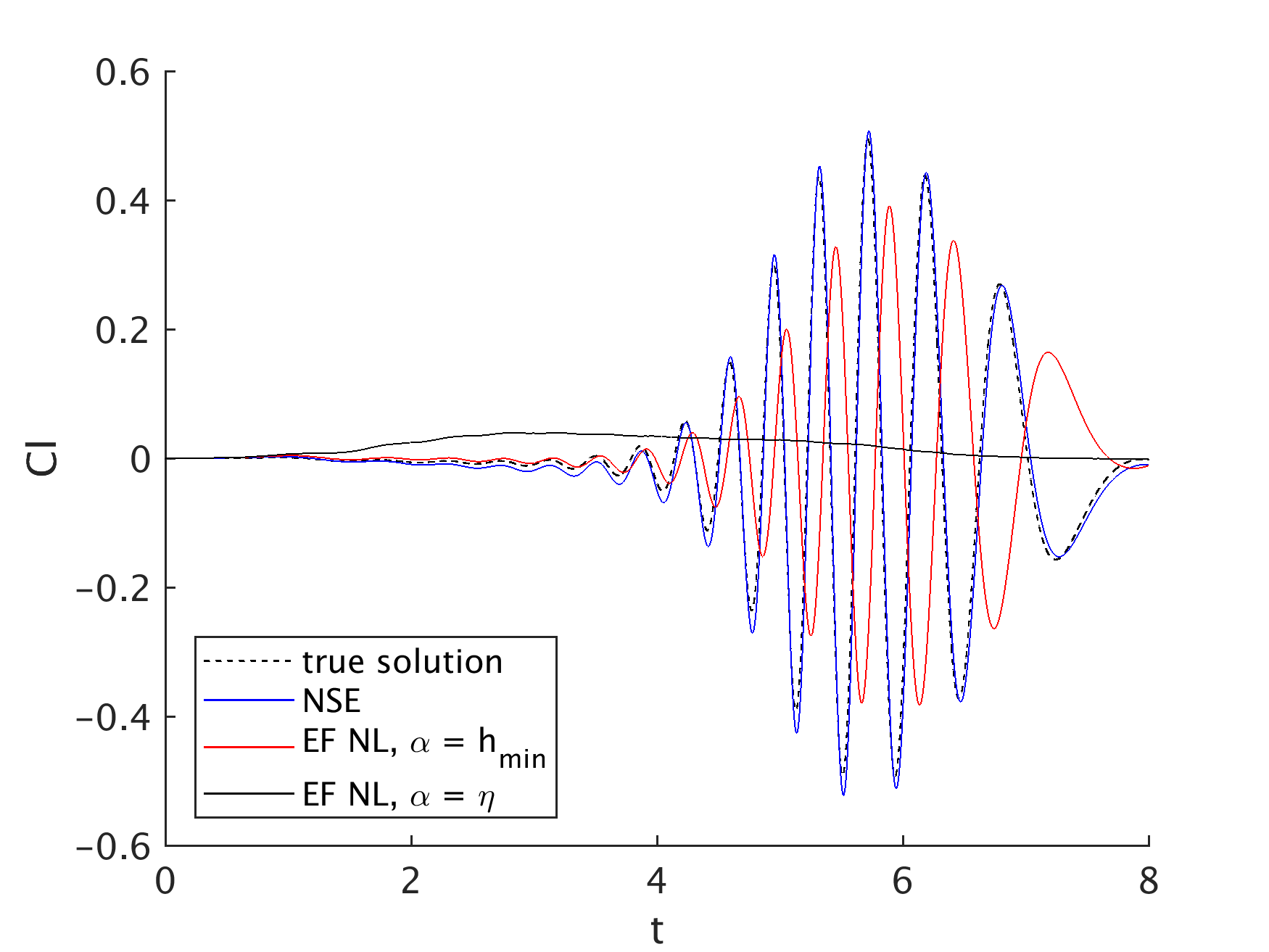}
        \put(32,70){\small{Mesh $120k_H$, EF NL}}
      \end{overpic}
\caption{Evolution over time of the lift coefficient given by the EF L (left) and EF NL (right) algorithms 
for $\alpha = h_{min}, \eta$ on all the hexahedral meshes in Table \ref{tab:1}. All the figures report also the 
NSE results.}
\label{fig:lift_H}
\end{figure}

\begin{figure}
\centering
 \begin{overpic}[width=0.45\textwidth]{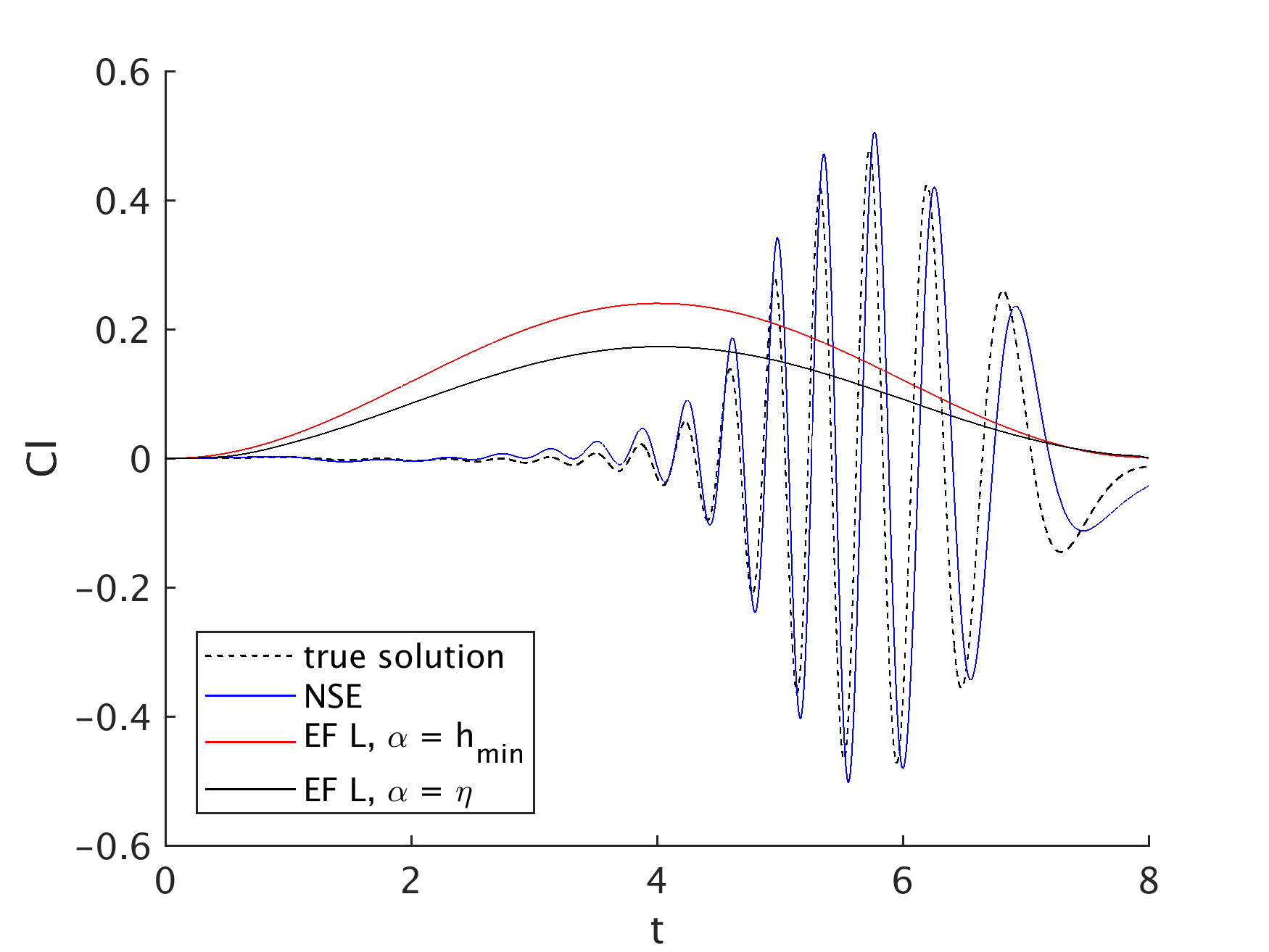}
        \put(32,70){\small{Mesh $16k_P$, EF L}}
      \end{overpic}
 \begin{overpic}[width=0.45\textwidth]{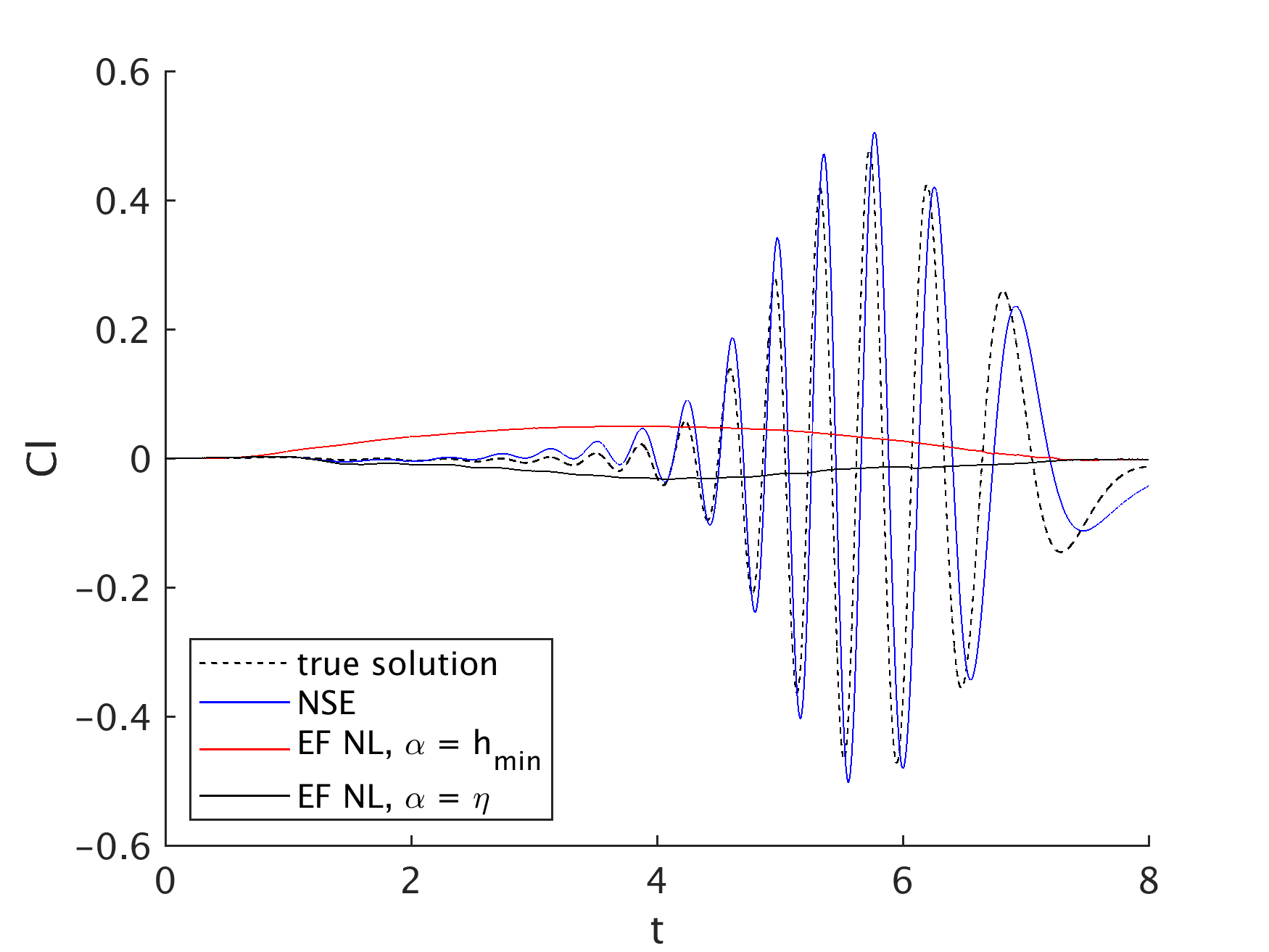}
        \put(32,70){\small{Mesh $16k_P$, EF NL}}
      \end{overpic}\\
       \begin{overpic}[width=0.45\textwidth]{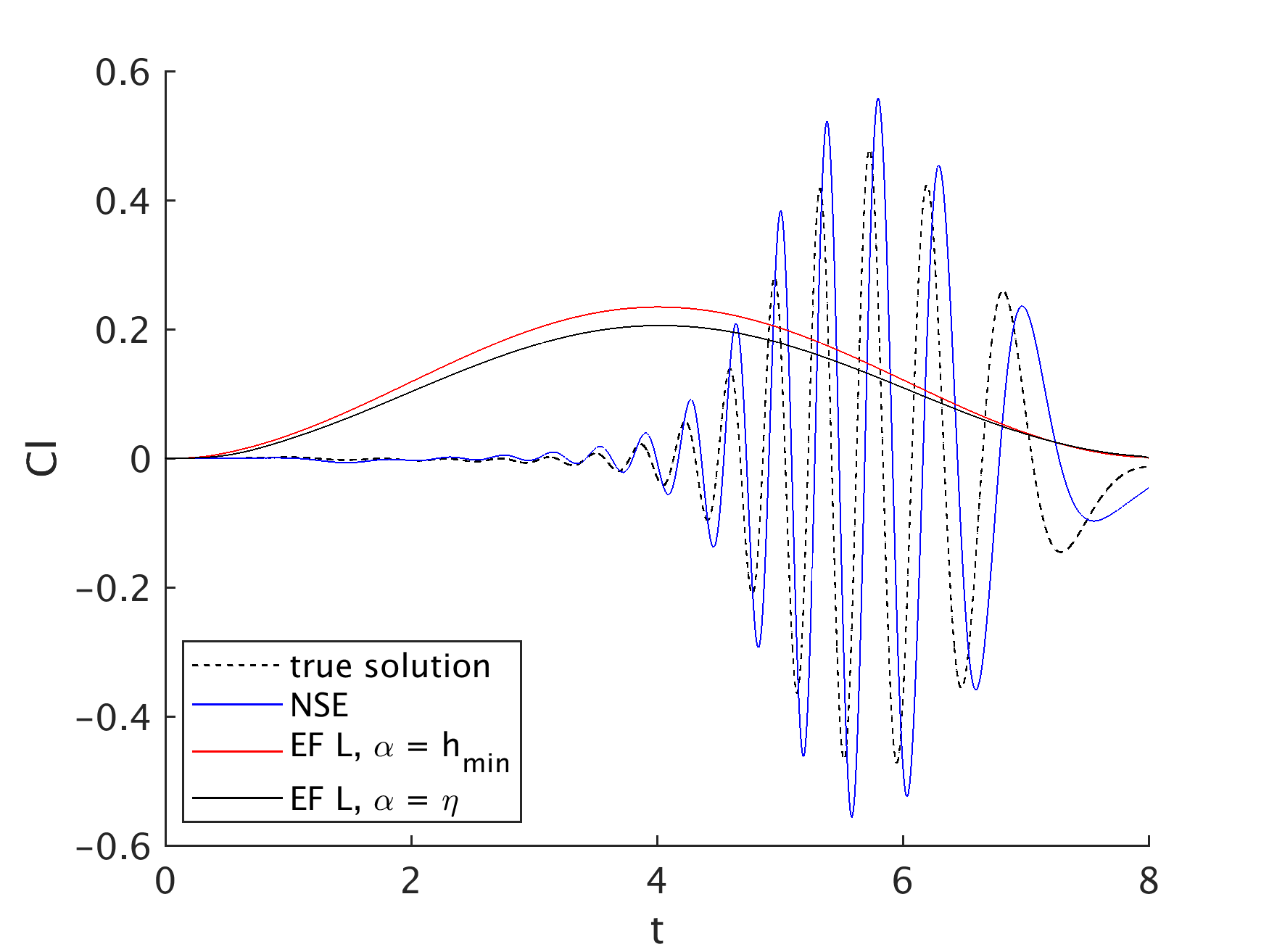}
        \put(32,70){\small{Mesh $25k_P$, EF L}}
      \end{overpic}
 \begin{overpic}[width=0.45\textwidth]{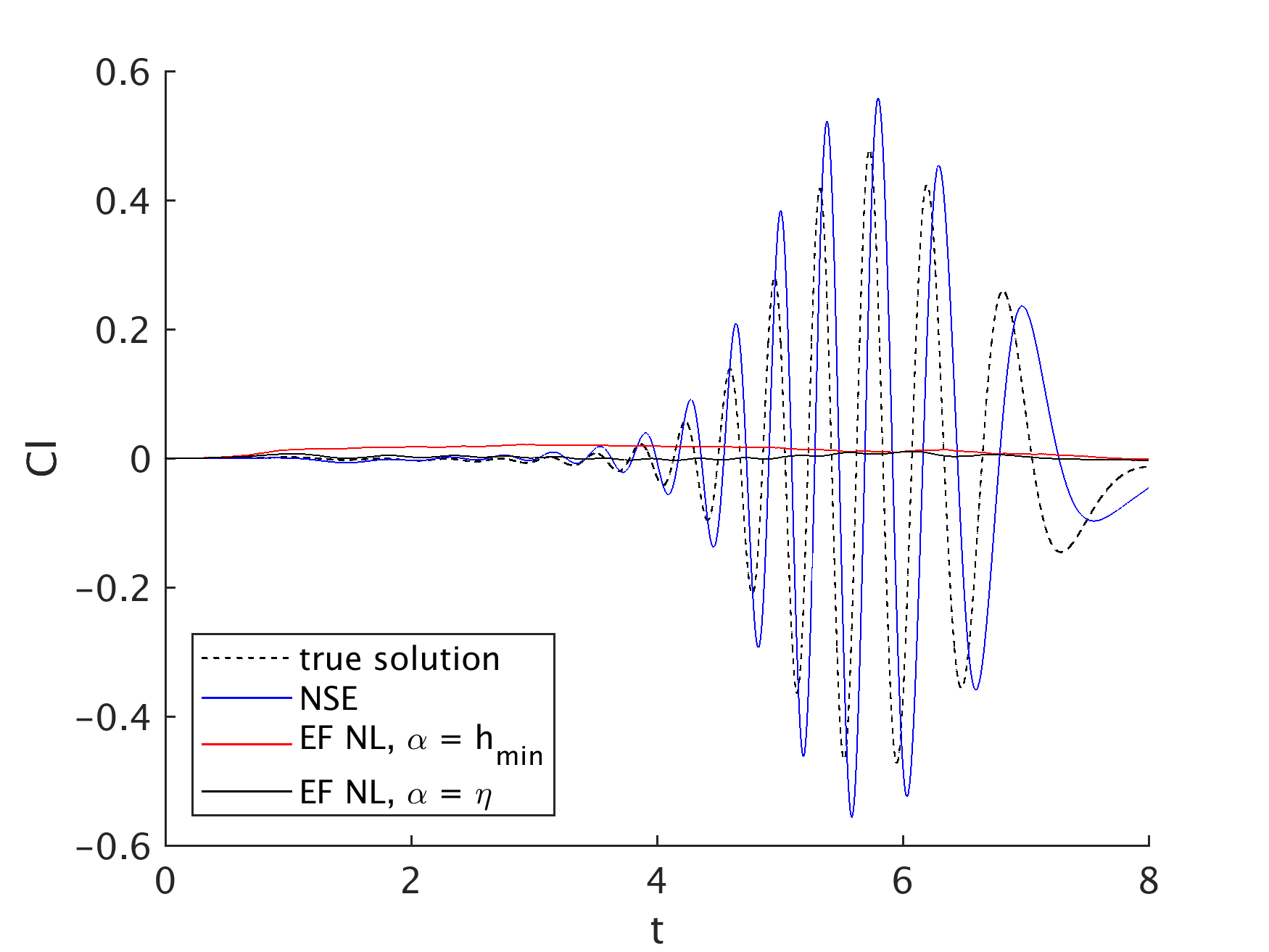}
        \put(32,70){\small{Mesh $25k_P$, EF NL}}
      \end{overpic}\\
       \begin{overpic}[width=0.45\textwidth]{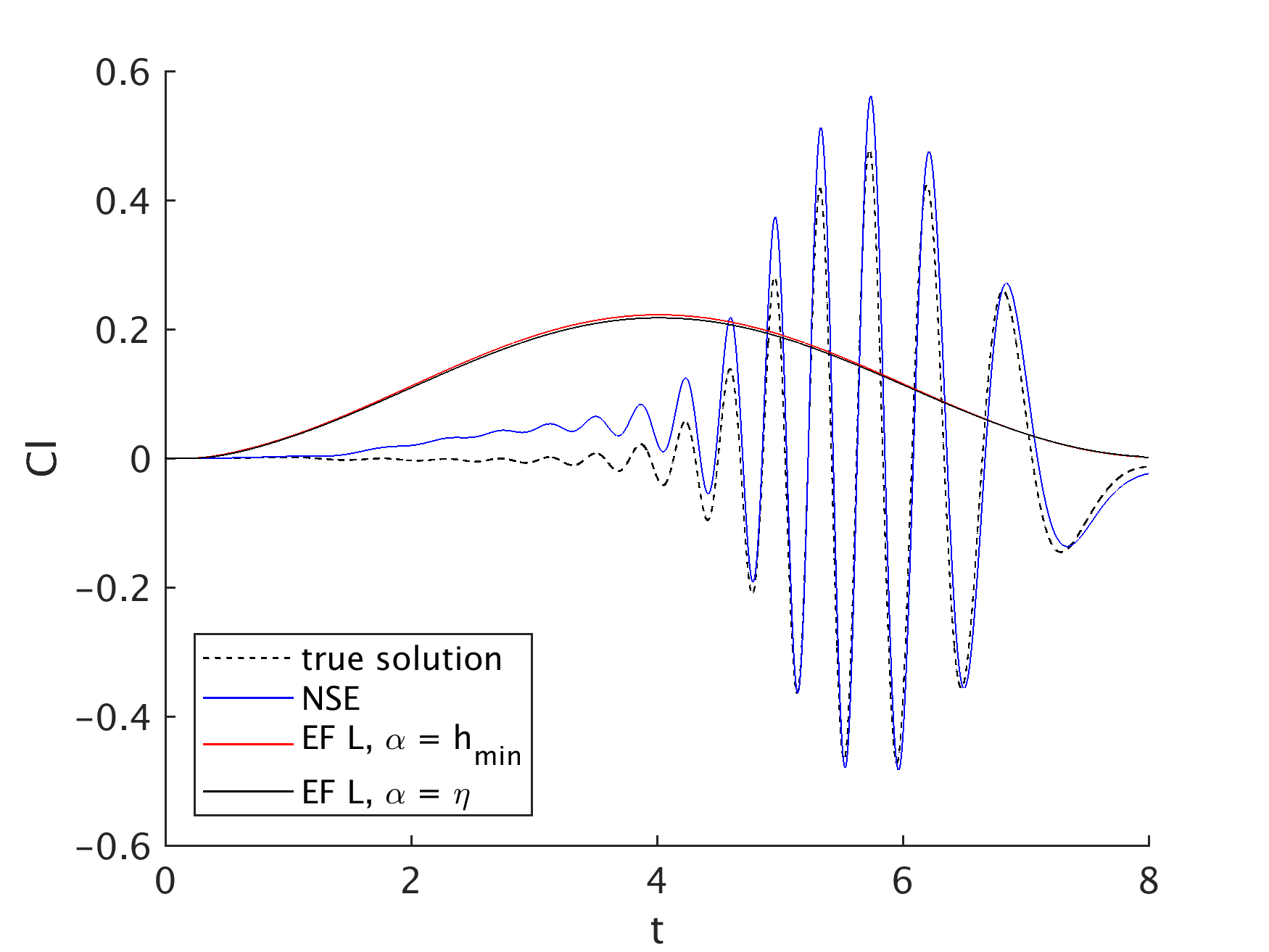}
        \put(32,70){\small{Mesh $63k_P$, EF L}}
      \end{overpic}
 \begin{overpic}[width=0.45\textwidth]{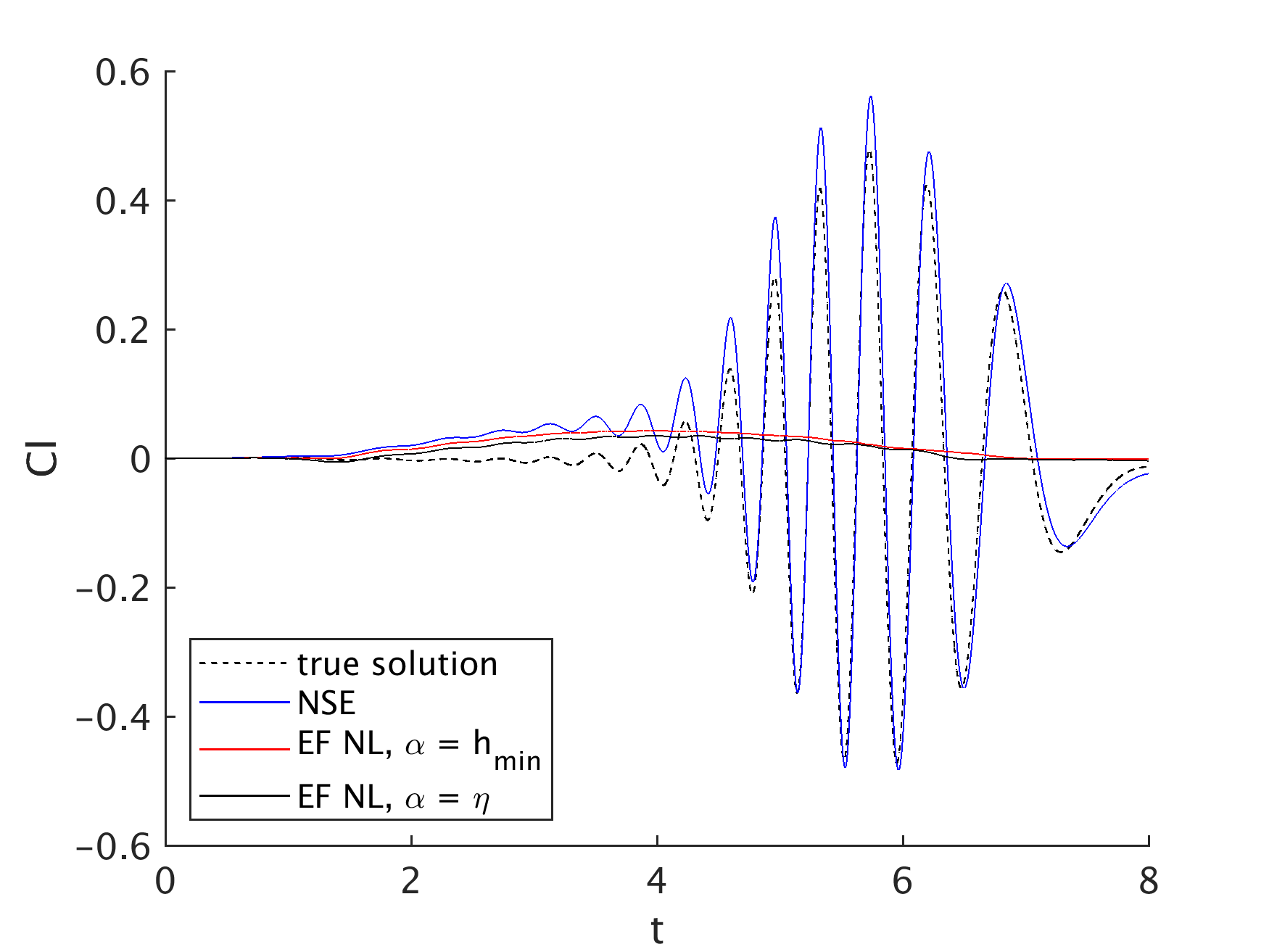}
        \put(32,70){\small{Mesh $63k_P$, EF NL}}
      \end{overpic}\\
       \begin{overpic}[width=0.45\textwidth]{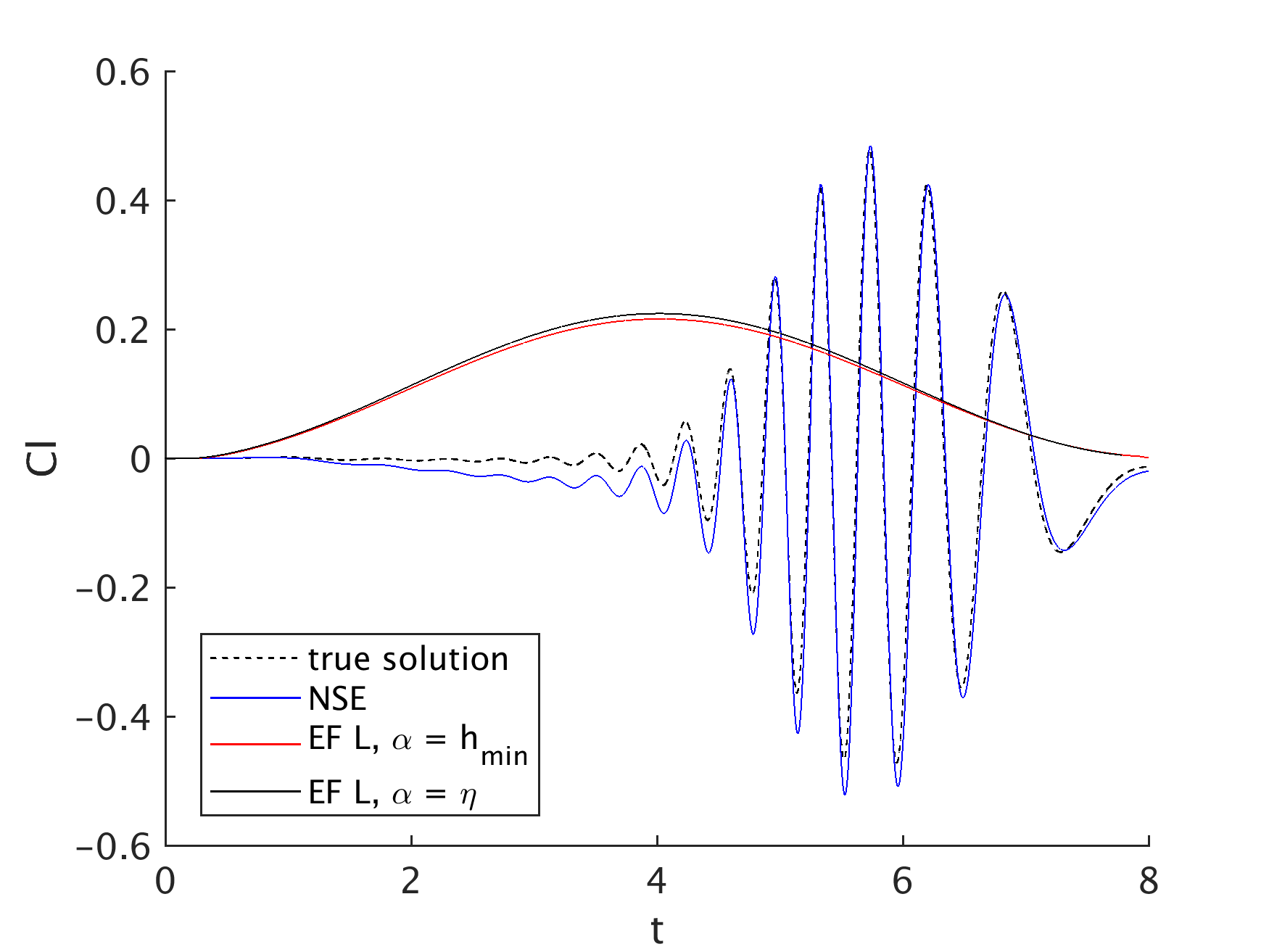}
        \put(32,70){\small{Mesh $120k_P$, EF L}}
      \end{overpic}
 \begin{overpic}[width=0.45\textwidth]{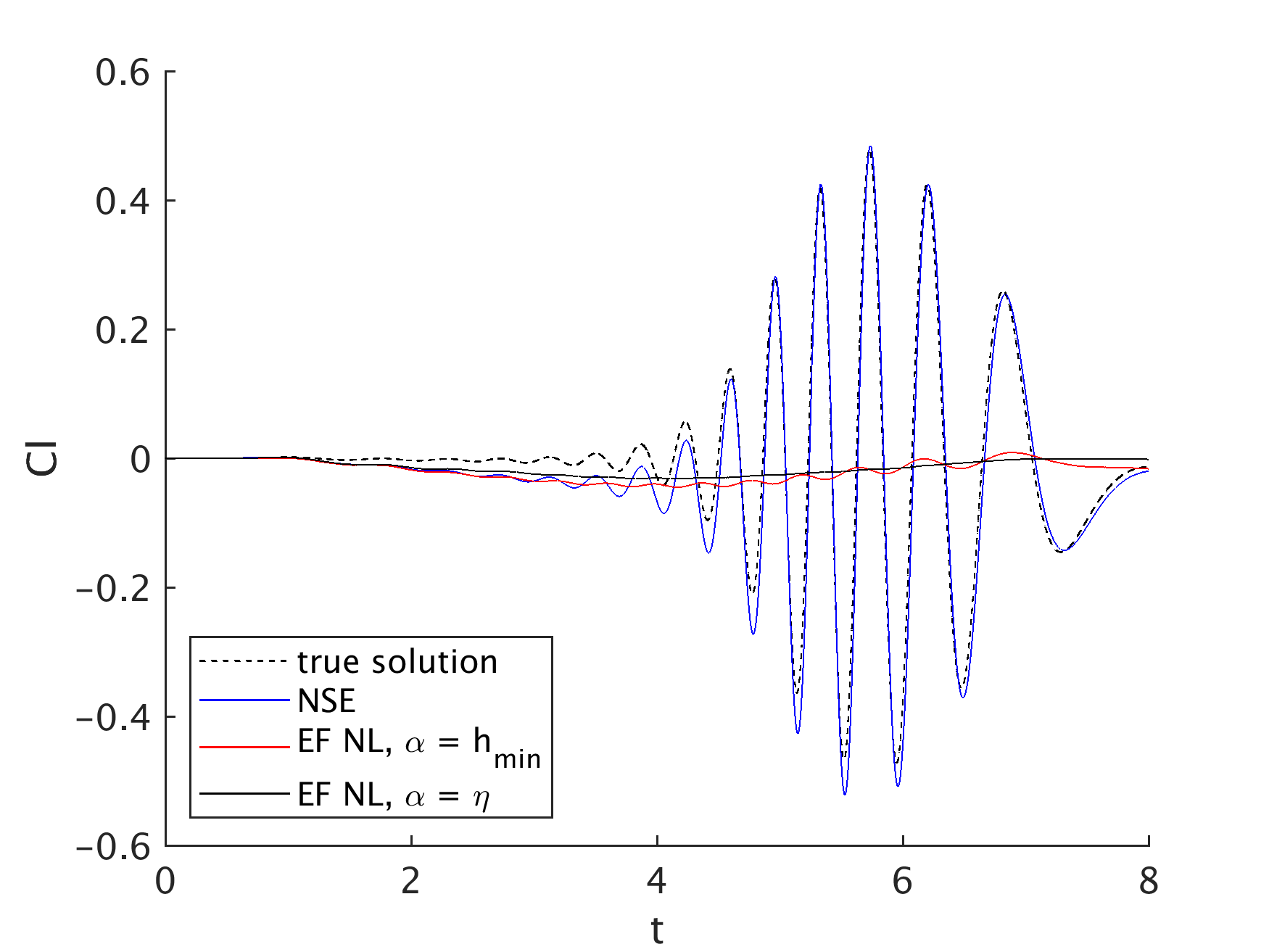}
        \put(32,70){\small{Mesh $120k_P$, EF NL}}
      \end{overpic}
\caption{Evolution of the lift coefficient given by the EF L (left) and EF NL (right)  algorithms
for $\alpha = h_{min}, \eta$ on all the prismatic meshes in Table \ref{tab:1}. All the figures report also the 
NSE results.
}
\label{fig:lift_P}
\end{figure}

Fig.~\ref{fig:lift_H} and \ref{fig:lift_P} show the evolution of $c_l$ over time
computed by the EF L and EF NL algorithms (i.e, $\chi = 0$) on all the meshes coarser than the DNS meshes
with hexahedral and prismatic elements, respectively.
Fig.~\ref{fig:lift_H} and \ref{fig:lift_P} report also the NSE results, i.e.~the same 
results shown in Fig.~\ref{fig:NSE_PH}, which were obtained with 
no filtering step. From Fig.~\ref{fig:lift_H} and \ref{fig:lift_P} 
we see that the $c_l$ computed 
on a given mesh with both EF L and EF NL algorithms is 
very far from the true $c_l$, and worse than then $c_l$ given by NSE regardless of the value of $\alpha$. 
In fact, as expected, the EF algorithm reduces the peaks due to excessive artificial diffusion. 
Moreover, the linear filter completely dampens the high frequency modes 
for either choice of $\alpha$ and for both kind of meshes.
Fig.~\ref{fig:lift_H} shows that also when using the EF NL algorithm on the hexahedral meshes
the computed results obtained get closer and closer to the true solution
as the mesh gets refined and for $\alpha = \text{min}\{h_{min},\eta\}$. 
This improvement of the computed $c_l$ is lost on the prismatic meshes. 
Focusing on the hexahedral meshes, we see
that the choice of $\alpha$ makes a visible difference in the computed $c_l$ on meshes finer that $25k$. 
In conclusion, we learned that the NSE algorithm is to be preferred to the EF algorithms and hexahedral meshes are superior to prismatic meshes.

\begin{figure}
\centering
 \begin{overpic}[width=0.45\textwidth]{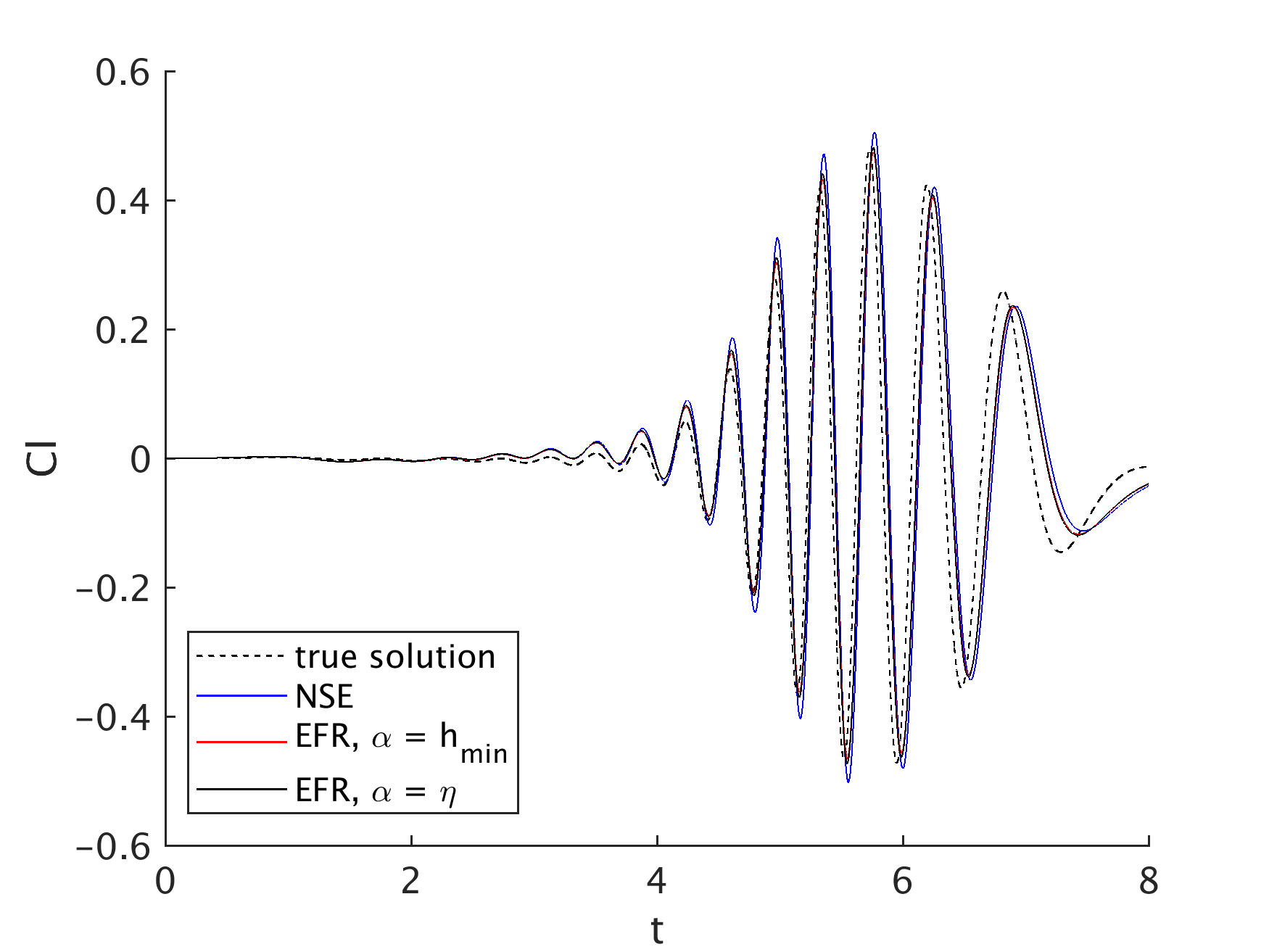}
        \put(32,70){\small{Mesh $16k_P$, EFR}}
      \end{overpic}
 \begin{overpic}[width=0.45\textwidth]{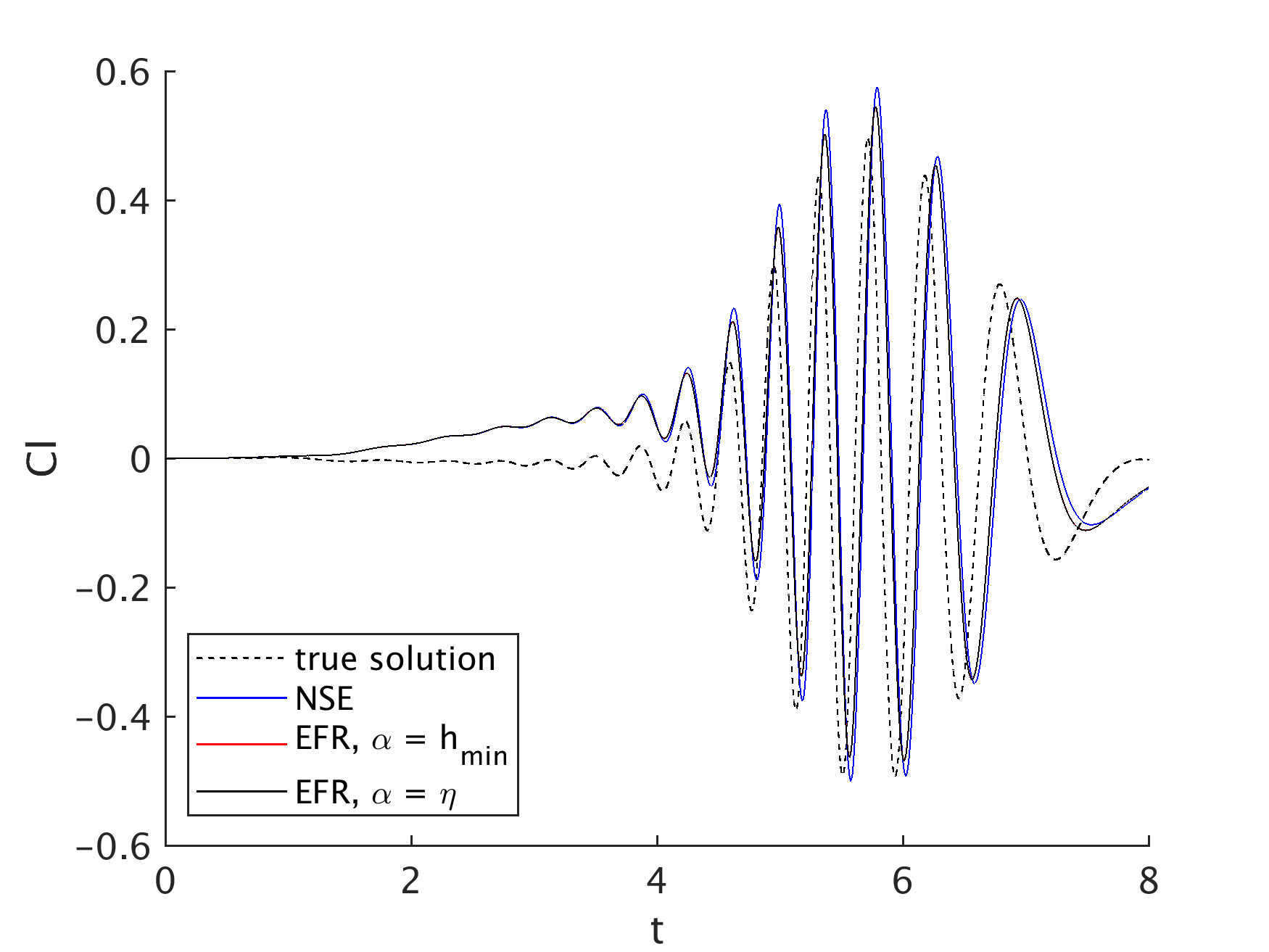}
        \put(32,70){\small{Mesh $16k_H$, EFR}}
      \end{overpic}\\
       \begin{overpic}[width=0.45\textwidth]{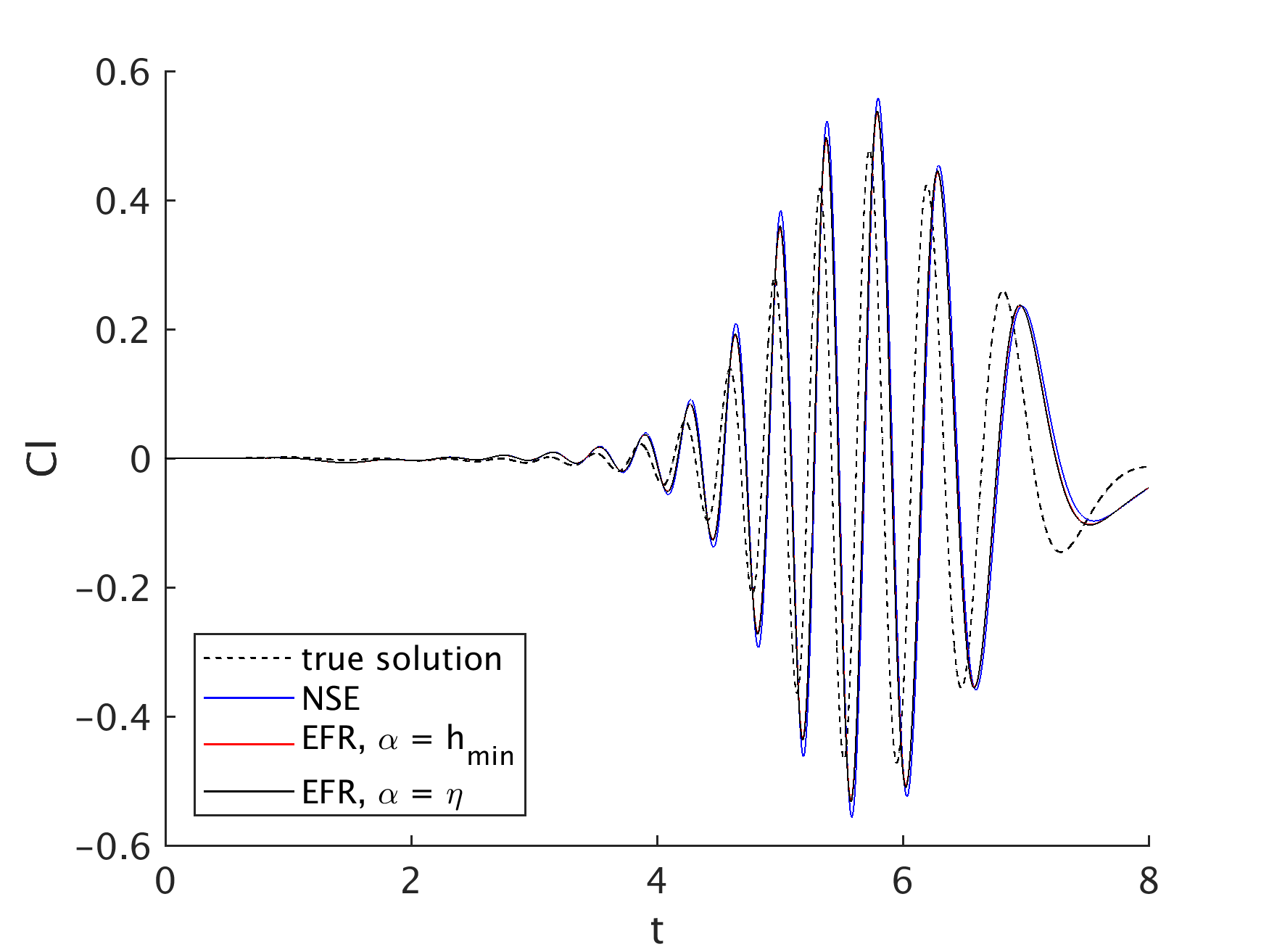}
        \put(32,70){\small{Mesh $25k_P$, EFR}}
      \end{overpic}
 \begin{overpic}[width=0.45\textwidth]{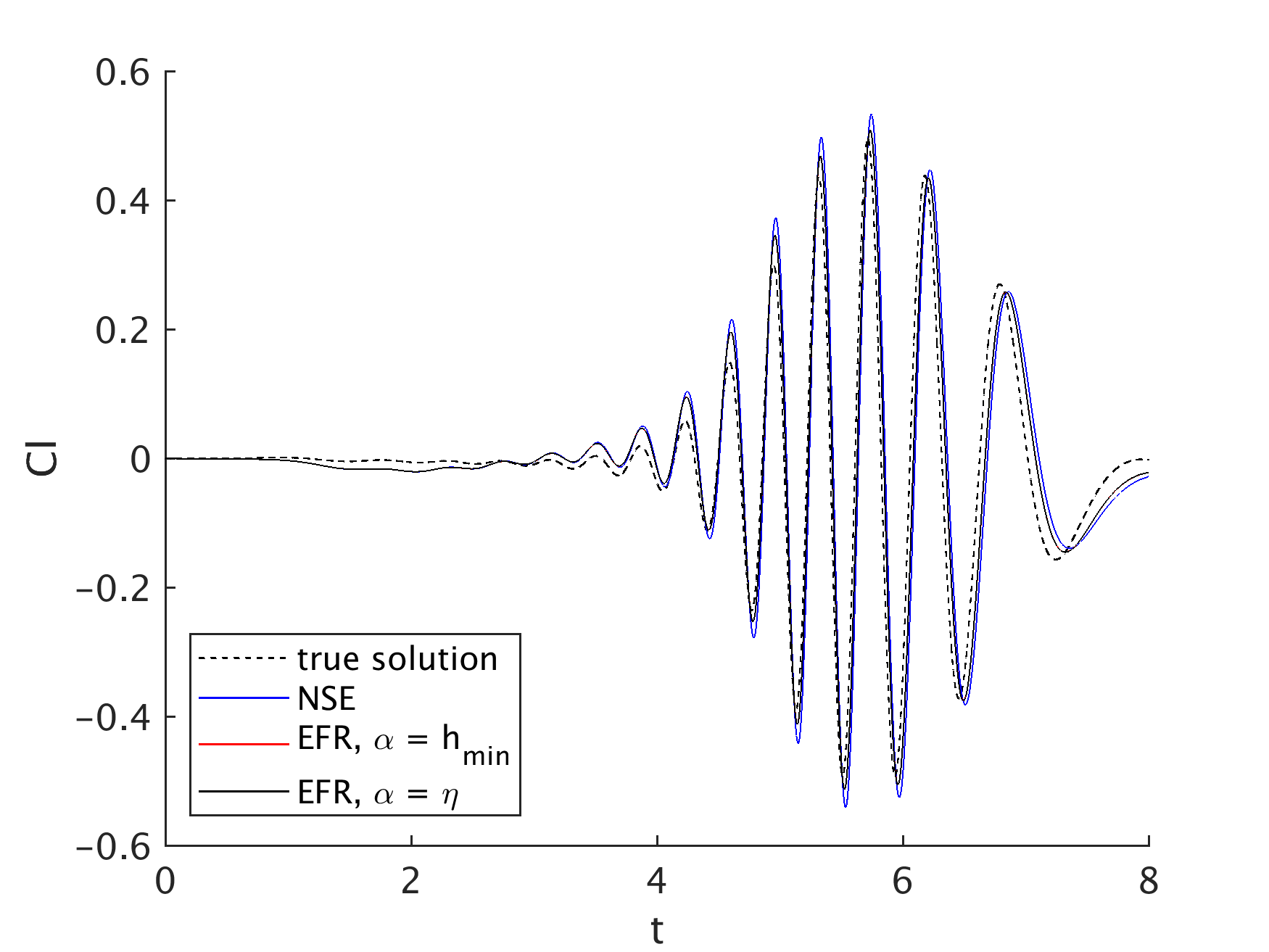}
        \put(32,70){\small{Mesh $25k_H$, EFR}}
      \end{overpic}\\
       \begin{overpic}[width=0.45\textwidth]{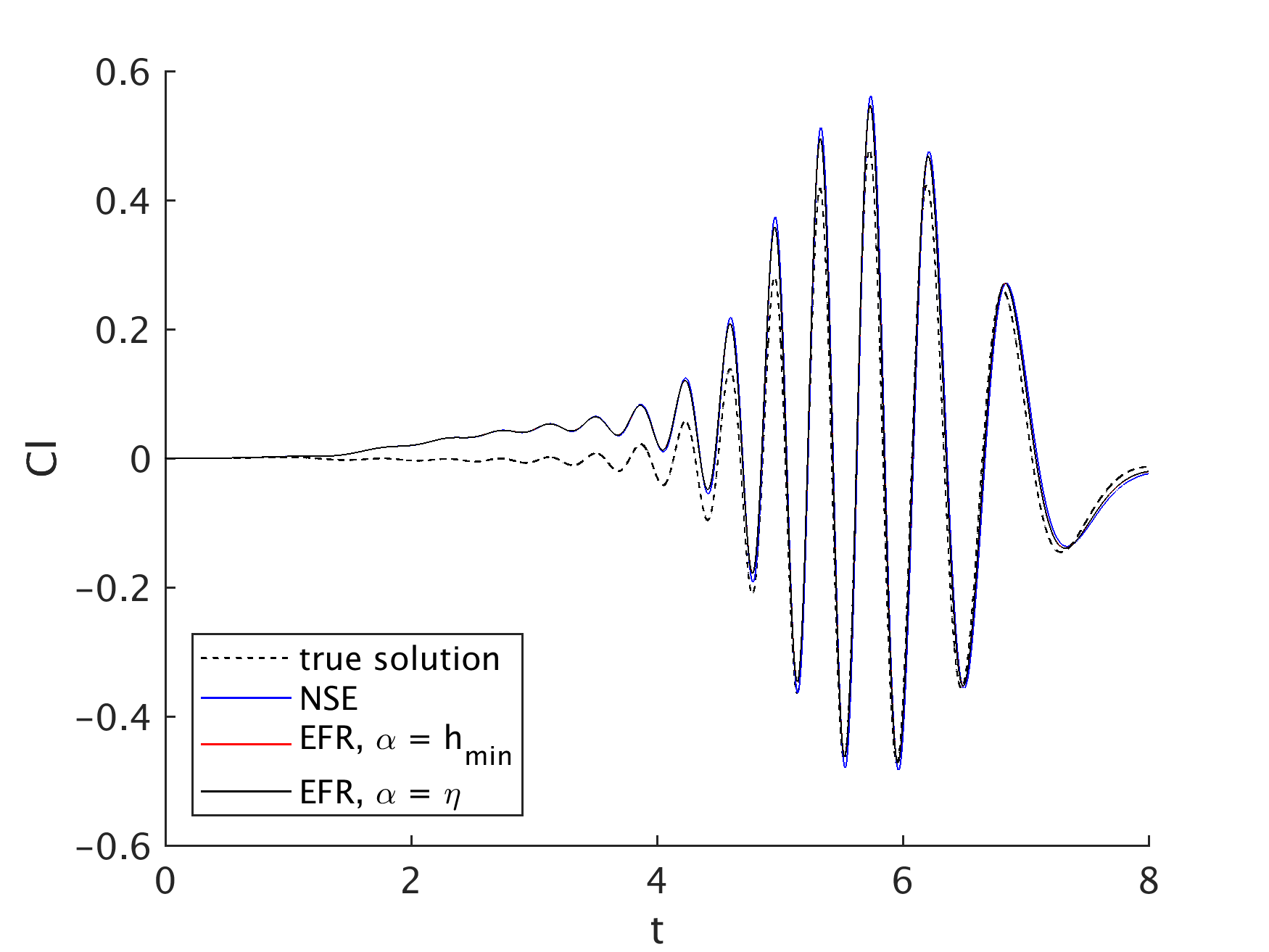}
        \put(32,70){\small{Mesh $63k_P$, EFR}}
      \end{overpic}
 \begin{overpic}[width=0.45\textwidth]{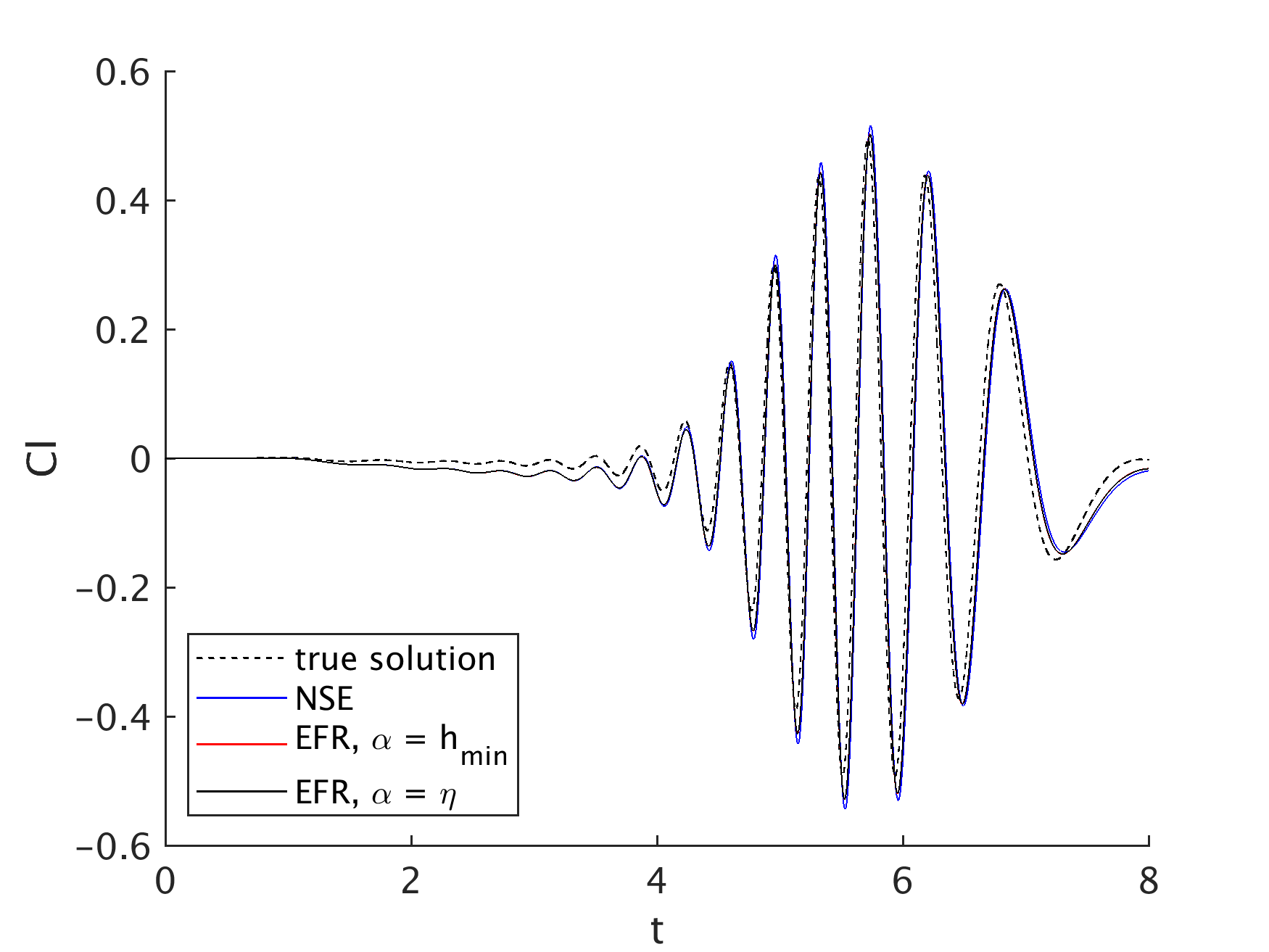}
        \put(32,70){\small{Mesh $63k_H$, EFR}}
      \end{overpic}\\
       \begin{overpic}[width=0.45\textwidth]{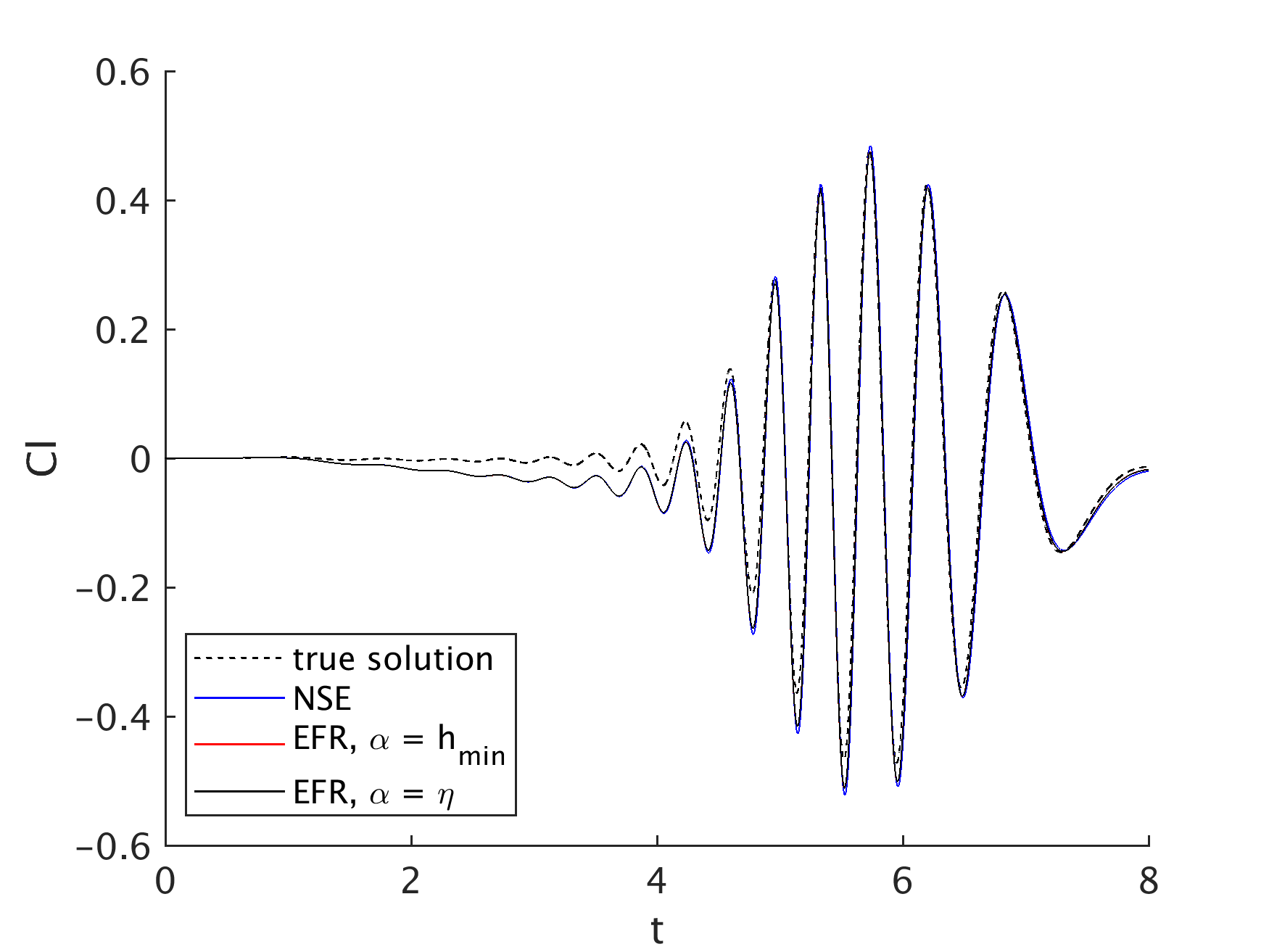}
        \put(32,70){\small{Mesh $120k_P$, EFR}}
      \end{overpic}
 \begin{overpic}[width=0.45\textwidth]{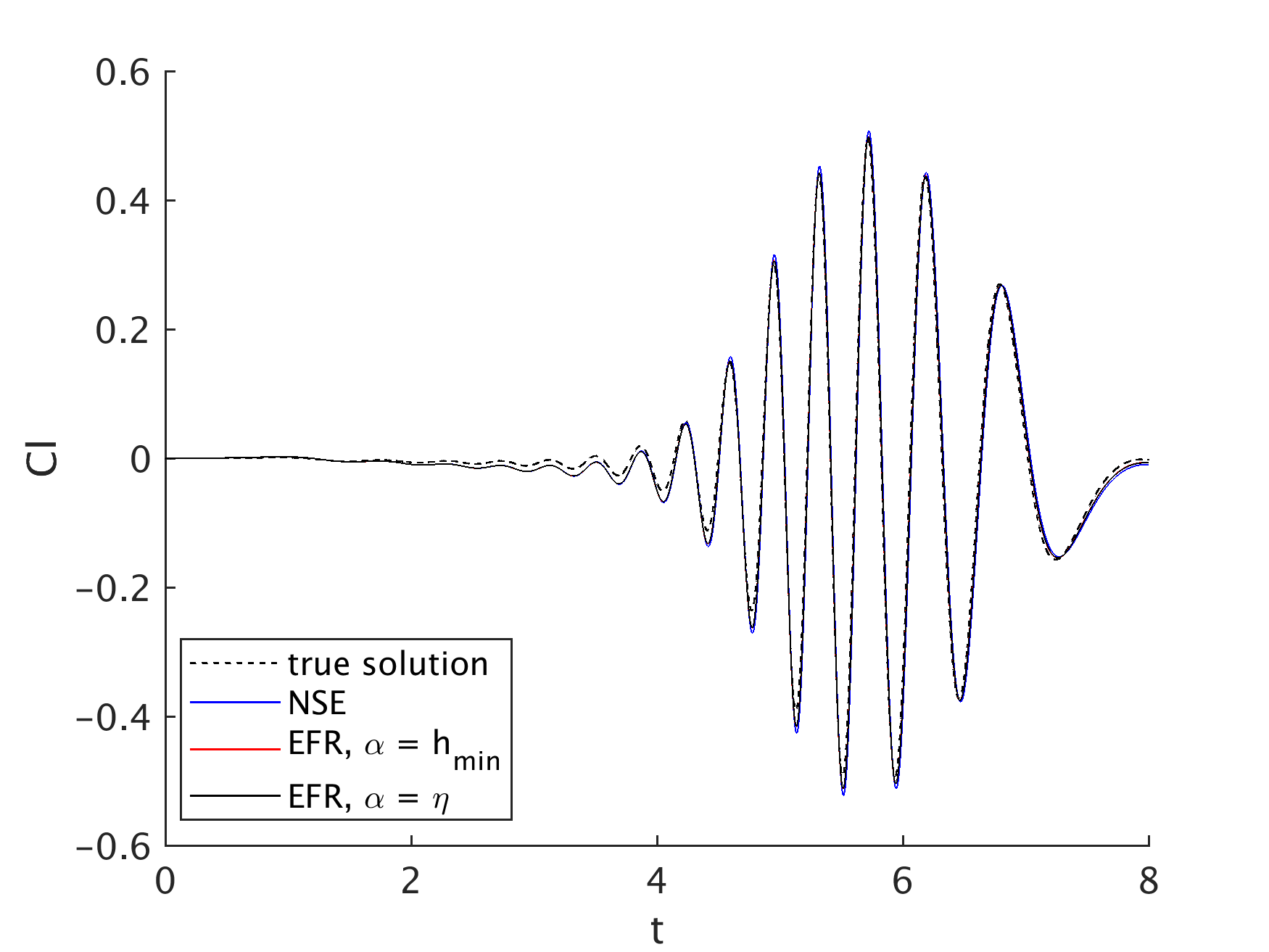}
        \put(32,70){\small{Mesh $120k_H$, EFR}}
      \end{overpic}
\caption{Evolution of the lift coefficient given by the EFR for $\alpha = h_{min}, \eta$ on 
all the prismatic (left) and hexahedral (right) meshes in Table \ref{tab:1}. All the figures report also the NSE results. }
\label{fig:lift_EFR}
\end{figure}

Since the EF L algorithm systematically introduces excessive numerical diffusion, we test the 
only the nonlinear version of the EFR algorithm, denoted hereafter
simply by EFR. We set $\chi  = \Delta t$. Fig.~\ref{fig:lift_EFR} 
shows the evolution of $c_l$ 
over time 
computed by the EFR algorithm on all the meshes coarser than the DNS meshes.
Fig.~\ref{fig:lift_EFR} 
reports also the NSE results. We notice that the EFR algorithm gives a lift coefficient 
very close to the true solution on all the meshes and for all
values of $\alpha$. We also observe that the filter helps reduce the phase advancement and the choice of $\alpha$ makes little-to-no visible difference. 
Moreover, it is confirmed that the results obtained with hexahedral meshes converge monotonically to the true solution with increase mesh refinement.

Finally, in Table \ref{tab:3} we report a quantitative comparison of all the simulations
in this section
in terms of maximum lift and drag coefficients and times at which the maxima occur.
Notice that for all the meshes the EFR algorithm with either $\alpha = h_{min}$
or $\alpha = \eta$ is the choice that minimizes $|E_{c_l}| + |E_{c_d}|$.

\begin{table}
\centering
\tiny
\begin{tabular}{cccccccc}
\cline{1-8}
Mesh name & Algorithm   & $t(c_{l,max})$ & $c_{l,max}$ & $t(c_{d,max})$ & $c_{d,max}$ & $E_{c_l}$ & $E_{c_d}$   \\
\hline
$16k_H$ & EF L, $\alpha$ = $h_{min}$ &  4.024 & 0.258  & 3.735 & 0.984  & -47.87\% &  -67.04\% \\
$16k_H$ & EF NL, $\alpha$ = $h_{min}$ & 2.303 & 0.022 & 3.943 & 3.358 & -95.55\% & 12.45\% \\
$16k_H$ & EF L, $\alpha$ = $\eta$  & 4 & 0.209 & 3.784 & 1.144 & -57.77\% & -61.68\% \\
$16k_H$ & EF NL, $\alpha$ = $\eta$  & 3.511 & 0.041 & 3.942 & 3.342 & -91.71\% & 11.92\%\\
$16k_H$ & EFR, $\alpha$ = $h_{min}$ & 5.777 & 0.544  & 3.94 & 3.074   & 9.89\% & 2.94\%\\
$16k_H$ & EFR, $\alpha$ = $\eta$  & 5.777 & 0.545 & 3.943 & 3.074 & 10.1\% & 2.94\%\\
\hline
$16k_P$ & EF L, $\alpha$ = $h_{min}$ & 3.989  & 0.24 & 3.581 & 0.734 & -49.79\% &  -75.48\% \\
$16k_P$ & EF NL, $\alpha$ = $h_{min}$ & 3.895 & 0.049 & 3.864 & 2.7 & -89.74\% & -9.82\%\\
$16k_P$ & EF L, $\alpha$ = $\eta$  & 4 & 0.173 & 3.762 & 1.026 & -63.8\% & -65.73\% \\
$16k_P$ & EF NL, $\alpha$ = $\eta$  & 0.905  & 0.0025 & 3.918 & 2.914 & -99.47\% & -2.67\% \\ 
$16k_P$ & EFR, $\alpha$ = $h_{min}$ & 5.76 & 0.474  & 3.931 & 2.905   & -0.83\% & -2.97\%\\
$16k_P$ & EFR, $\alpha$ = $\eta$  & 5.761 & 0.481 & 3.926 & 2.906 & 0.62\% & -2.94\%\\
\hline
\hline
$25k_H$ & EF L, $\alpha$ = $h_{min}$ & 4.008 & 0.168 & 3.738 & 0.993 & -66.06\% & -66.74\%  \\
$25k_H$ & EF NL, $\alpha$ = $h_{min}$ &  3.918 & 0.038 & 3.942 & 3.683 & -92.32\% & 23.34\% \\
$25k_H$ & EF L, $\alpha$ = $\eta$  &  4.012 & 0.157 & 3.763 & 1.039 & -68.28\% & -65.2\% \\
$25k_H$ & EF NL, $\alpha$ = $\eta$  &  4.249 & 0.036 & 3.974 & 3.641 & -92.72\% & 21.93\%\\
$25k_H$ & EFR, $\alpha$ = $h_{min}$ & 5.734 & 0.509  & 3.936 & 3.058   & 2.82\% & 2.41\%\\
$25k_H$ & EFR, $\alpha$ = $\eta$  & 5.734 & 0.509 & 3.935 & 3.058 & 2.82\% & 2.41\%\\
\hline
$25k_P$ & EF L, $\alpha$ = $h_{min}$  & 4.003 & 0.234 & 3.608 & 0.756 & -51.04\% & -74.75\%  \\
$25k_P$ & EF NL, $\alpha$ = $h_{min}$ & 3 & 0.021 & 3.918 & 3.526 &  -95.6\% & 17.76\%\\
$25k_P$ & EF L, $\alpha$ = $\eta$  &  4.011 & 0.205 & 3.732 & 0.954 &  -57.11\% & -68.13\%\\
$25k_P$ & EF NL, $\alpha$ = $\eta$  &  6.078 & 0.01 & 3.936 & 3.546 & -97.9\% & 18.43\% \\
$25k_P$ & EFR, $\alpha$ = $h_{min}$ & 5.792 & 0.534  & 3.942 & 3.067   & 11.71\% & 2.43\%\\
$25k_P$ & EFR, $\alpha$ = $\eta$  & 5.792 & 0.538 & 3.931 & 3.067 & 12.55\% & 2.43\%\\
\hline
\hline
$63k_H$ & EF L, $\alpha$ = $h_{min}$ & 4.03 & 0.113 & 3.803 & 1.166 & -77.17\% & -60.95\% \\
$63k_H$ & EF NL, $\alpha$ = $h_{min}$ & 5.994 & 0.246 & 3.961 & 3.573 & -50.3\% & 19.65\% \\
$63k_H$ & EF L, $\alpha$ = $\eta$  &  4.031 & 0.186 & 3.701 & 0.89 & -62.42\% & -70.19\% \\
$63k_H$ & EF NL, $\alpha$ = $\eta$  & 3.571 & 0.05 & 3.944 & 3.888 &  -89.89\% &  30.2\%\\
$63k_H$ & EFR, $\alpha$ = $h_{min}$ & 5.733 & 0.502 & 3.942 & 3.028   & 1.41\% & 1.4\%\\
$63k_H$ & EFR, $\alpha$ = $\eta$  & 5.734 & 0.502 & 3.941 & 3.029 & 1.41\% & 1.44\%\\
\hline
$63k_P$ & EF L, $\alpha$ = $h_{min}$  & 4.006 & 0.222 & 3.656 & 0.806 & -53.55\% & -73.08\%   \\
$63k_P$ & EF NL, $\alpha$ = $h_{min}$  &  3.947 & 0.042 & 3.95 & 4.118 & -91.21\% & 37.54\%\\
$63k_P$ & EF L, $\alpha$ = $\eta$  & 4 & 0.217 & 3.688 & 0.853 &  -54.6\% & -71.5\% \\
$63k_P$ & EF NL, $\alpha$ = $\eta$  &  3.98 & 0.035 & 3.95 & 4.062 & -92.67\% &  35.67\%\\
$63k_P$ & EFR, $\alpha$ = $h_{min}$ & 5.735 & 0.546 & 3.945 & 3.044   & 14.22\% & 1.67\%\\
$63k_P$ & EFR, $\alpha$ = $\eta$  & 5.734 & 0.547 & 3.935 & 3.044 & 14.43\% & 1.67\%\\
\hline
\hline
$120k_H$ & EF L, $\alpha$ = $h_{min}$  &  4.016 & 0.112 & 3.813 & 1.2 &  -77.37\% & -59.81\% \\
$120k_H$ & EF NL, $\alpha$ = $h_{min}$  & 5.891 & 0.391 & 4 & 3.37 & -21.01\% & 12.86\% \\
$120k_H$ & EF L, $\alpha$ = $\eta$  &  3.994 & 0.218 & 3.622 & 0.79 & -55.96\% & -73.54\% \\
$120k_H$ & EF NL, $\alpha$ = $\eta$  & 2.833 & 0.039 & 3.948 & 4.105 & -92.12\% & 37.47\% \\
$120k_H$ & EFR, $\alpha$ = $h_{min}$ & 5.72 & 0.498 & 3.936 & 3.012   & 0.6\% & 0.87\%\\
$120k_H$ & EFR, $\alpha$ = $\eta$  & 5.72 & 0.498 & 3.943 & 3.013 & 0.6\% & 0.9\%\\
\hline
$120k_P$ & EF L, $\alpha$ = $h_{min}$  &  4.038 & 0.215 & 3.674 & 0.832 & -55.02\% & -72.21\% \\
$120k_P$ & EF NL, $\alpha$ = $h_{min}$  &  6.883 & 0.008 & 3.966 & 4.04 &  -98.32\% & 34.93\%\\
$120k_P$ & EF L, $\alpha$ = $\eta$  &  4.04 & 0.224 & 3.658 & 0.791 & -53.13\% &  -73.58\% \\
$120k_P$ & EF NL, $\alpha$ = $\eta$  &  0.926 & 0.0014 & 3.977 & 4.189 &  -99.7\% & 39.91\%\\
$120k_P$ & EFR, $\alpha$ = $h_{min}$ & 5.734 & 0.474 & 3.94 & 3.013   & -0.83\% & 0.63\%\\
$120k_P$ & EFR, $\alpha$ = $\eta$  & 5.734 & 0.474 & 3.944 & 3.014 & -0.83\% & 0.66\%\\
\hline
\end{tabular}
\caption{Maximum lift and drag coefficients 
and times at which the maxima occur for the EF L, EF NL and EFR
algorithms with $\alpha = h_{min}, \eta$ on all the meshes
under consideration. 
}
\label{tab:3}
\end{table}

\subsection{3D benchmark: FDA}\label{sec:FDA}
 
\begin{figure}[h]
\centering
\subfloat[][sketch of the FDA nozzle geometry]{\includegraphics[width=0.48\textwidth]{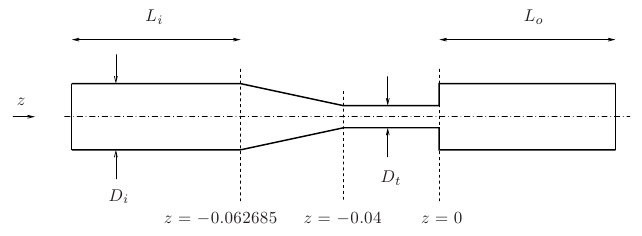}}~
\subfloat[][mesh $100k$]{\includegraphics[width=0.3\textwidth]{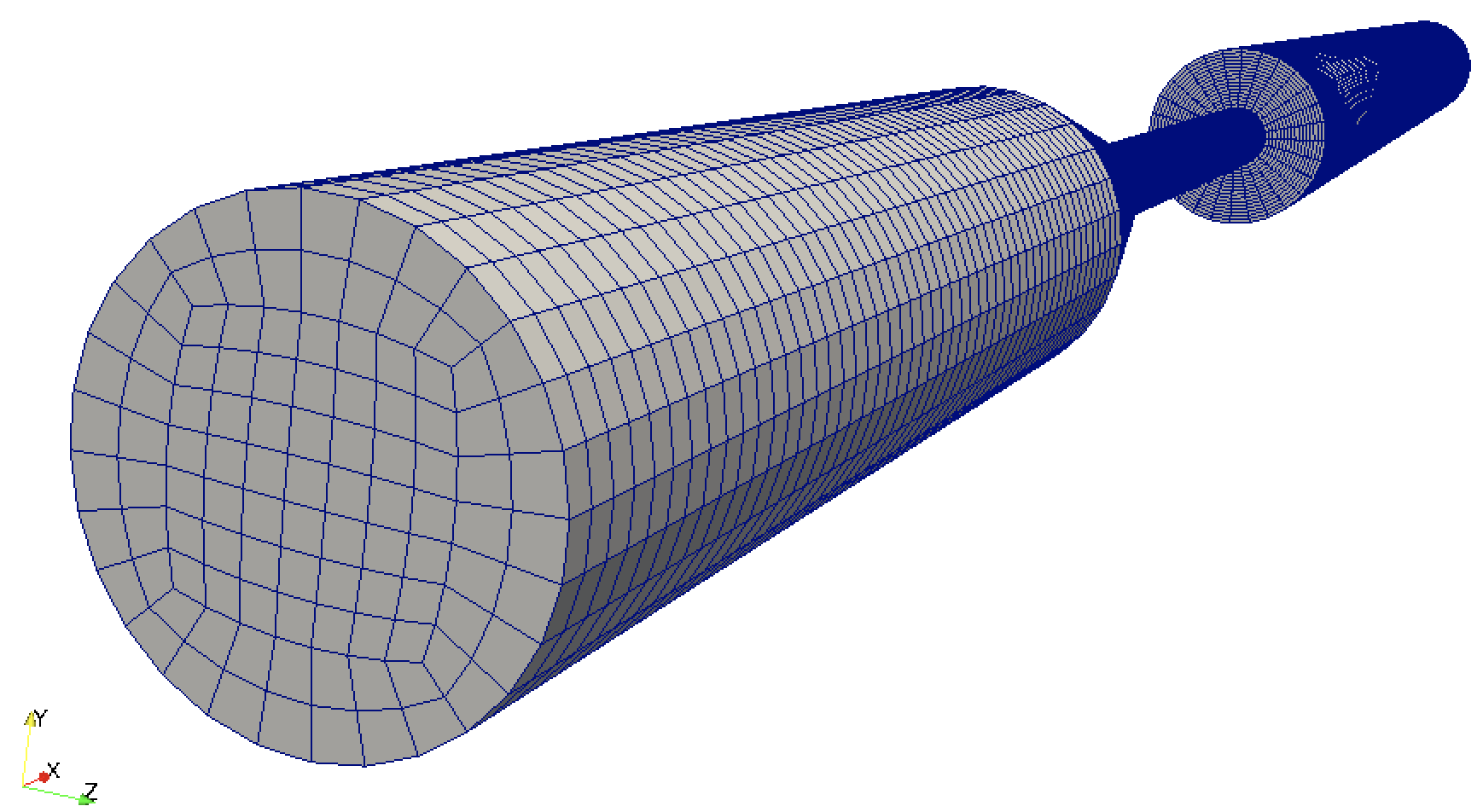}}
\caption{(a) A section of the computational domain, with $D_i = 0.012$ m, $D_t = 0.004$ m, $L_i = 4 D_i$ and $L_o = 12 D_i$ and (b) a view of mesh $100k$.}\label{fig:FDAdomain}
\end{figure}

The second test we consider is a benchmark from the U.S. Food and Drug Administration (FDA),
issued within the ``Critical Path Initiative'' program \cite{cpi}.
The objective of this benchmark is to simulate the flow of an incompressible and Newtonian fluid 
in nozzle geometry for different flow regimes, from laminar to fully turbulent. 
Despite its relative simplicity of the geometry, the nozzle contains all the features commonly encountered in medical devices,
i.e.~flow contraction and expansion, recirculation zones etc. See Fig.~\ref{fig:FDAdomain}. 
Three independent  laboratories were requested by FDA  to perform flow visualization experiments on fabricated nozzle models for different flow rates up to Reynolds number 6500 \cite{hariharang}. 
This resulted in benchmark data available online to the scientific community for the validation of CFD simulations \cite{fdacfd}. Available experimental measurements  enable us to  check 
the effectiveness of the EFR algorithm in simulating average macroscopic quantities.
The results of the published inter-laboratory experiments refer to values of the Reynolds 
numbers in the throat (defined as in \eqref{eq:re}) of ${Re_t = 500, 2000,3500, 5000, 6500}$,
$Re_t = 2000$ being the critical Reynolds number for transitional flow
in a straight pipe \cite{reynolds1883}.
To test our methodology, we focus on Reynolds numbers
$Re_t = 2000, 3500, 5000, 6500$.

In Table \ref{tab:FDA_1}, we report the throat Reynolds number $Re_t$, 
the corresponding inlet Reynolds number $Re_i$, the flow rate $Q$, the Kolmogorov 
length scale $\eta$ and the number of cells required by DNS for the considered flow regimes. 
The value of $\eta$ was found by plugging into \eqref{eq:eta-re} the value of 
$Re_t$ and the diameter of the expansion channel $D_i$ as characteristic length. 
As for the 2D benchmark, the number of cells for the DNS was calculated by assuming a mean grid 
width equal to $\eta$ as spatial resolution. We notice that the number of cells reported in Table \ref{tab:FDA_1}
assumes a uniform mesh that has the highest level of refinement possible, since $\eta$
is calculated using the highest Reynolds number in the domain. A mesh with a local refinement dictated by the local 
Reynolds number would have a more manageable number of cells.


\begin{table}[h]
\centering
\begin{tabular}{ccccc}
\multicolumn{2}{c}{} \\
\cline{1-5}
$Re_t$    & $Re_i$ & $Q$ (m$^3$/s) & $\eta$ (m) & No.~of cells for DNS \\
\hline
2000     &  667 & 2.0825e-5  & 4.0124e-5  & 3.6285e8   \\
3500      & 1167  & 3.6444e-5   & 2.6371e-5  &  1.2781e9 \\
5000       &  1667  &  5.2062e-5  &  2.0182e-5  & 2.8513e9 \\
6500       &  2167  &  6.7681e-5  &  1.6577e-5  &  5.1455e9\\
\hline
\end{tabular}
\caption{Throat Reynolds number $Re_t$, inlet Reynolds number $Re_i$, flow rate $Q$, Kolmogorov length scale $\eta$ and number of cells for a DNS for all the flow regimes under consideration.}
\label{tab:FDA_1}
\end{table}

The fluid used in the experiments has $\rho = 1056$ Kg/m$^3$ and $\mu = 0.0035$ Pa/s.
Following \cite{BQV}, the length of the inlet chamber $L_i$ was set to four times its diameter, 
while the length of the expansion channel ($L_o$ in Fig.~\ref{fig:FDAdomain}) was set to
12 times its diameter. On the lateral surface of the computational domain we prescribe a no-slip boundary
condition. 
At the inflow we prescribe a Poiseuille velocity profile:
\begin{align}\label{eq:Poiselle}
\v (r, \theta, z) = \left(0, 0, 2v_{mean}\left(1 - \dfrac{r^2}{R_{i}^2}\right)\right),
\end{align}
where $r$, $\theta$ and $z$ are the radial, polar and axial locations respectively, $v_{mean}$ is the mean inlet velocity
magnitude to obtain the desired flow rate,
and $R_i = D_i/2$. 
At the outlet section, we prescribe an advective outflow condition:
\begin{align}\label{eq:free-stress}
\dfrac{\partial \v}{\partial t} + a \dfrac{\partial \v}{\partial \n} = 0,
\end{align}
where $a$ is the advection velocity magnitude computed so that the total mass is conserved.
Boundary condition \eqref{eq:free-stress} preserves the numerical stability
in case the jet starts to break down too close to the outlet section \cite{Nicoud}.
The actual experimental set up of the FDA benchmark is 
a closed loop \cite{hariharang} and it is not reflected by our condition. However, 
this is expected to introduce a minimal error localized only in the close neighborhood of the outlet section. 
The results of the flow analysis are not affected,
provided that the expansion channel is long enough.
We recall that
the partitioned algorithms we use require a boundary condition for the pressure as well: 
we impose ${\partial q}/{\partial \n} = {\partial \qbar}/{\partial \n} = 0$ at the inlet and the wall and $p = \qbar = 0$ at the outlet.


We start all the simulations from fluid at rest and use a smooth transition of the inlet velocity profile 
to regime flow conditions. 
We let the simulations run till $T = 0.8$ s
and then we compare with experimental data the solution obtained by averaging the snapshots collected 
between $t = 0.4$ s and $t = 0.8$ s
\cite{Passerini2013}. 
The comparison between computational and experimental data is made in terms of normalized axial component of the velocity and normalized pressure difference along the centerline. The axial component of the velocity $u_z$ is normalized with respect to the average axial velocity at the inlet $\bar{u}_i$:
\begin{align}\label{eq:uz}
{u_z}^{norm} = \dfrac{u_z}{\bar{u}_i}, ~\text{with} \quad \bar{u}_i = \dfrac{Q}{\pi {D_i}^2/4}.
\end{align}
The pressure difference is normalized with respect to the dynamic pressure in the throat:
\begin{align}\label{eq:pz}
\Delta p^{norm}= \dfrac{p_z - p_{z=0}}{1/2 \rho {\bar{u}_t}^2}, ~\text{with} \quad \bar{u}_t = \dfrac{Q}{\pi {D_t}^2/4},
\end{align}
where $p_z$ denotes the wall pressure along the $z$-axis and $p_{z=0}$ is the wall pressure at $z = 0$, 
i.e.~at the sudden expansion.

In the following, we compare our computational results with the experimental data from
\cite{hariharang,fdacfd}. Moreover, we compare our findings with \cite{BQV} and \cite{Nicoud}. 
In \cite{BQV}, the same LES method is used but it is combined with a Finite Element 
method ($\mathbb{P}_2$-$\mathbb{P}_1$ elements) implemented in \emph{LifeV} \cite{lifev}. 
In \cite{Nicoud} the authors use a different LES method, called $\sigma$-model \cite{Nicoud2011}, and combine it
with a central fourth-order Finite Volume scheme implemented in \emph{YALES2BIO} \cite{YALES2BIO}.
In \cite{BQV} BDF2 is used for the time discretization, while an explicit Runge-Kutta scheme 
of the fourth order (RK4) and an advanced fourth-order accurate time scheme (TFV4A)
are used in \cite{Nicoud}.

Based on the results in Sec.~\ref{sec:cylinder}, we consider only meshes with hexahedral elements. 
For the mesh generation we used the mesh generation utility, \emph{blockMesh} within 
OpenFOAM. Table \ref{tab:mesh_FDA} reports
the mesh name, minimum and maximum diameter, number of cells, and average number of cells along 
a throat diameter for all the meshes under consideration.
Just like in Sec.~\ref{sec:cylinder}, the name of each grid refers to the number of cells (Figure ~\ref{fig:FDAdomain}).
We generated the meshes to have a similar refinement to the tetrahedral meshes 
used in \cite{BQV}, while the tetrahedral meshes considered in \cite{Nicoud} have a much larger number of cells.
See Table \ref{tab:mesh_BQV_Nicoud}.
Also for this benchmark, we made sure that all the meshes feature high quality: 
low values of maximum non-orthogonality (20$^\circ$ to 25$^\circ$),  average 
non-orthogonality (around 4$^\circ$), skewness (around 1), and maximum aspect ratio (up to 7). 
Meshes were selectively refined in the convergent, throat, and sudden expansion.

\begin{table}[h]
\centering
\begin{tabular}{ccccc}
\hline
mesh name    & $h_{min}$ & $h_{max}$ & No. of cells & $D_t/h$   \\
\hline
$100k$      & 2.6784e-4   & 1.1e-3   & 1.03e5 &  11 \\
$320k$     &  1.8186e-4   & 7.7958e-4   & 3.21e5 &  16 \\
$770k$      &  1.3611e-4  &  5.9389e-4  & 7.72e5 &  21   \\
$1600k$      & 1.0613e-4    & 4.6966e-4  &  1.59e6 & 27\\
$2700k$      &  8.8289e-5   & 3.947e-4  & 2.66e6 &  32\\
$6400k$      &   6.529e-5  &  2.957e-4 &   6.38e6 & 43 \\
\hline
\end{tabular}
\caption{Mesh name, minimum diameter $h_{min}$, maximum diameter $h_{max}$, 
number of cells, and average number of cells along a throat diameter $D_t/h$ for all the meshes 
under consideration.}
\label{tab:mesh_FDA}
\end{table}

\begin{table}[h]
\centering
\begin{tabular}{cccc|cccc}
\multicolumn{2}{c}{} \\
\hline
mesh name    & $h_{min}$ & $h_{max}$ &  $D_t/h$ & mesh name  & $h_{min}$ & $h_{max}$ & $D_t/h$ \\
\hline
$140k$      & 3.39e-4   & 3.09e-3    & 10 & $5000k$    & 3.2e-4   & 3.6e-4    & 13 \\
$330k$     &  2.23e-4   & 1.93e-3    & 11 & $15000k$     &  1.9e-4   & 2.3-4    & 21    \\
$900k$      &  1.09e-4  &  1.87e-3   & 19  &  $50000k$      &  1.2e-4  &  1.7e-4   & 34 \\
$1200k$      & 1.08e-4    & 1.63e-3  & 19 & & &  &  \\
$1900k$      & 1.06e-4    & 1.49e-3   & 20 & & & &  \\
$3000k$      &  1.17e-4   & 9.64e-4  & 21  & & & & \\
\hline
\end{tabular}
\caption{Mesh name, minimum diameter $h_{min}$, maximum diameter $h_{max}$, 
and average number of cells along a throat diameter $D_t/h$ for the tetrahedral meshes 
used in \cite{BQV} (left) and \cite{Nicoud} (right).}
\label{tab:mesh_BQV_Nicoud}
\end{table}

For the sake of completeness, we performed additional tests with the tetrahedral grids used in \cite{BQV} (see Table \ref{tab:mesh_BQV_Nicoud} left) that worked well within a finite element framework. We found 
that the high non-orthogonality of those meshes (up to 55$^\circ$) significantly affects 
the computational time and the accuracy of the computed solutions. Thus, we decided to disregard those
results and focus on low non-orthogonality hexahedral meshes that generally represent the best choice 
in the finite volume context \cite{jasakphd}.
However, we observe that tetrahedral grids with very good features in terms of skewness and aspect ratio were adopted for finite volume 
simulations in \cite{Nicoud}. 
Moreover, the authors state that additional tests on hexahedral meshes (not reported in the manuscript) 
did not show substantial differences, confirming that the quality of their tetrahedral meshes
is high.

We start to investigate the cases of $Re_t = 3500$ and above. We present the results 
for the transitional case ($Re_t = 2000$) last, since it is known being a tough test 
both from the experimental and the numerical point of view.  

\subsubsection{Case $Re_t = 3500$}\label{sec:3500}

In \cite{Passerini2013,BQV}, it was shown that tetrahedral mesh $1200k$ (see Table \ref{tab:mesh_BQV_Nicoud} (left)) 
has a level of refinement that gives numerical results obtained with a finite element method
in excellent agreement with the experimental data, despite the average mesh size being roughly 20 times larger than the Kolmogorov length scale at $Re_t = 3500$. 
Thus, we are going to consider meshes with a similar level of refinement or coarser meshes. 

\begin{figure}[h]
\centering
 \begin{overpic}[width=0.45\textwidth]{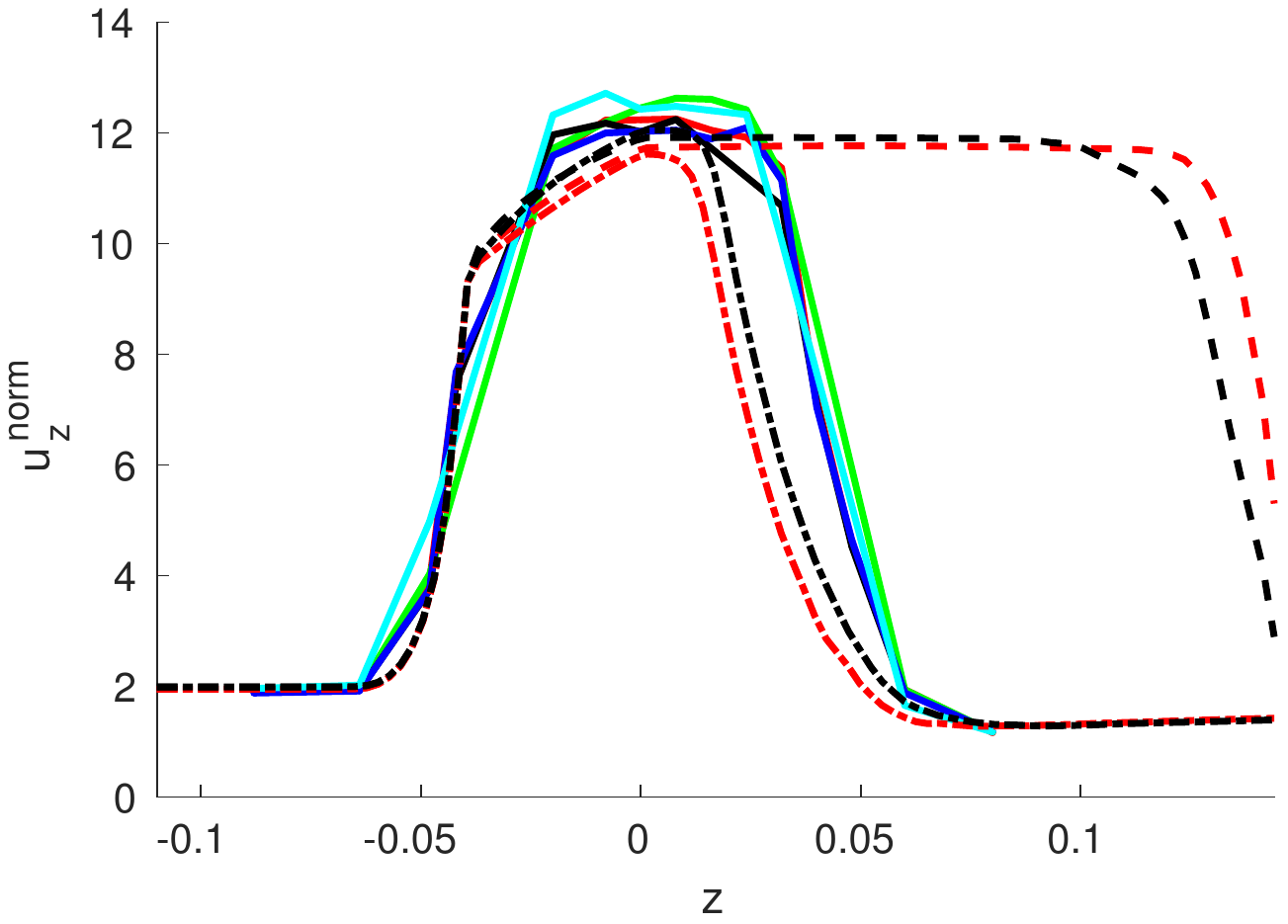}
        \put(15,70){\small{normalized axial velocity along $z$}}
      \end{overpic}
 \begin{overpic}[width=0.45\textwidth]{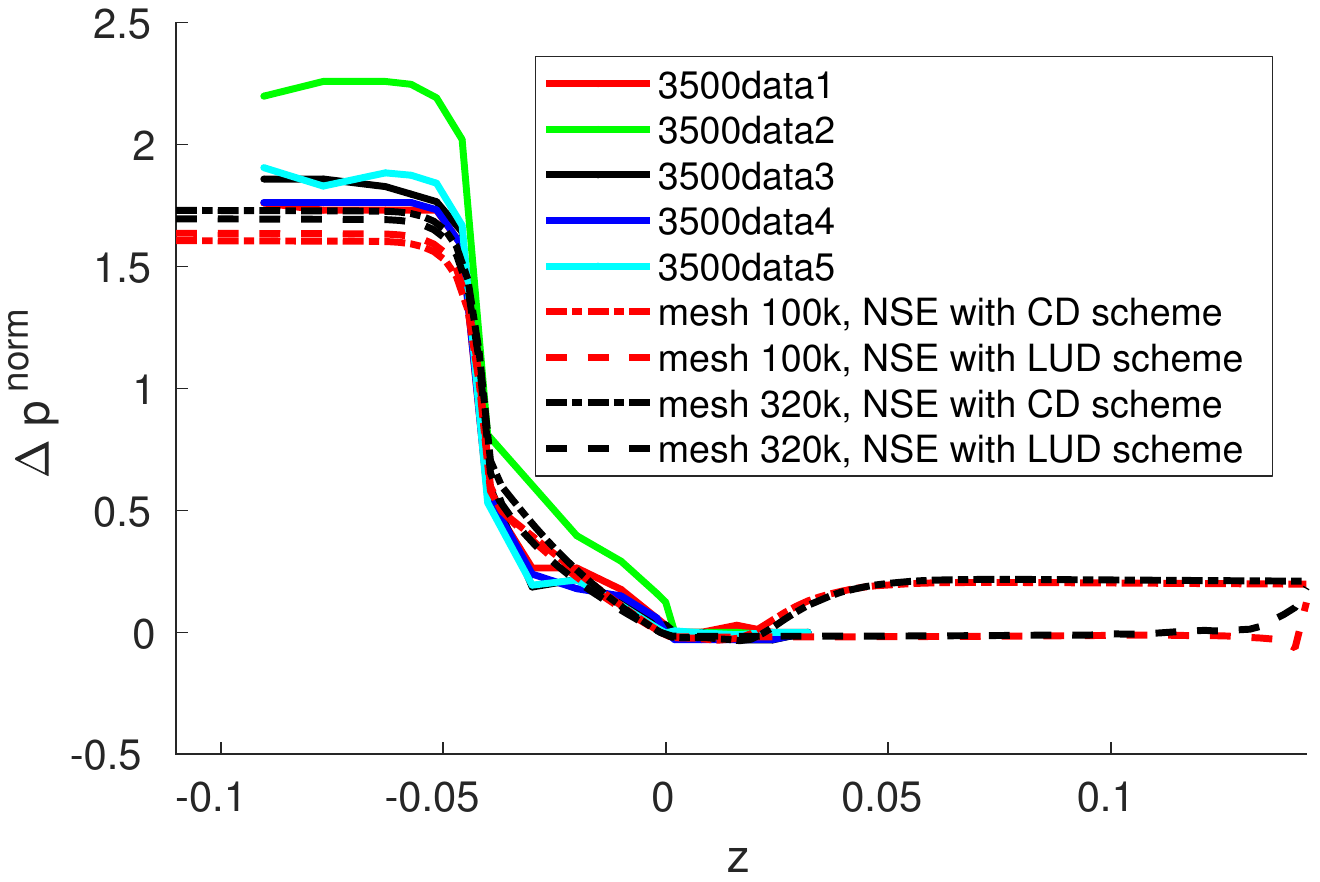}
       \put(10,70){\small{normalized pressure difference along $z$}}
      \end{overpic}
\caption{Case $Re_t = 3500$, NSE: comparison between experimental data (solid lines) 
and numerical results (dashed lines) for the normalized axial velocity \eqref{eq:uz} (left)
and the normalized pressure drop \eqref{eq:pz} (right) along the $z$. The numerical 
results were obtained with meshes $100k$ and $320k$ and
two different interpolation schemes for the convective term.
The legend on the right is common to both graphs.}
\label{fig:FDA_lUCD}
\end{figure}

It is known that the numerical results for the FDA benchmark test are particularly
sensitive to small perturbations in the mesh, geometry, discretization parameters,
and numerical schemes \cite{stewartp,Passerini2013,BQV,Nicoud,Bergensen2018,Janiga2014,Delorme2013,Bhushan2013,Chabannes2017,Fehn2018,Nicoud2018,Wu2019}. 
We start with investigating the effect of interpolation scheme
for the convective term. 
In particular, we consider two second-order accurate schemes: Linear Upwind Differencing (LUD) 
\cite{Warming1976}, which was used in Sec. \ref{sec:cylinder}, and 
Central Differencing (CD) \cite{Lax1960}.
Fig.~\ref{fig:FDA_lUCD} reports the normalized axial velocity \eqref{eq:uz}
and the normalized pressure drop \eqref{eq:pz} 
along the $z$ axis from the NSE solution on meshes $100k$ and $320k$ (see Table \ref{tab:mesh_FDA})
obtained by using the LUD and the CD convective schemes. We set $CFL = 0.6$ as in \cite{Nicoud}. 
The jet given by the NSE algorithm with the LUD convective scheme is much longer than 
the experimental jet for both meshes, while the NSE algorithm with the CD convective scheme 
provides a smaller jet length. 
It is known that upwind schemes introduce large diffusive errors and have a significant damping effect on the energy spectrum of the flow \cite{simonsphd}.
Therefore, central difference schemes are to be preferred for LES \cite{simonsphd}.  
From now on, we will use the CD scheme for both NSE and EFR algorithms. 
In addition, we will use the original version of PISO since the flux correction term
has been shown to introduce artificial dissipation \cite{ddtPhiCorr}.
From Fig.~\ref{fig:FDA_lUCD} (left) we observe that the pressure computed with either scheme and mesh
is in excellent agreement with the experimental data for $z < 0.3$. For $z > 0.3$ the difference in the
pressures computed with the two schemes for the convective term becomes large
and experimental data are not available.

\begin{figure}[h]
\centering
 \begin{overpic}[width=0.45\textwidth]{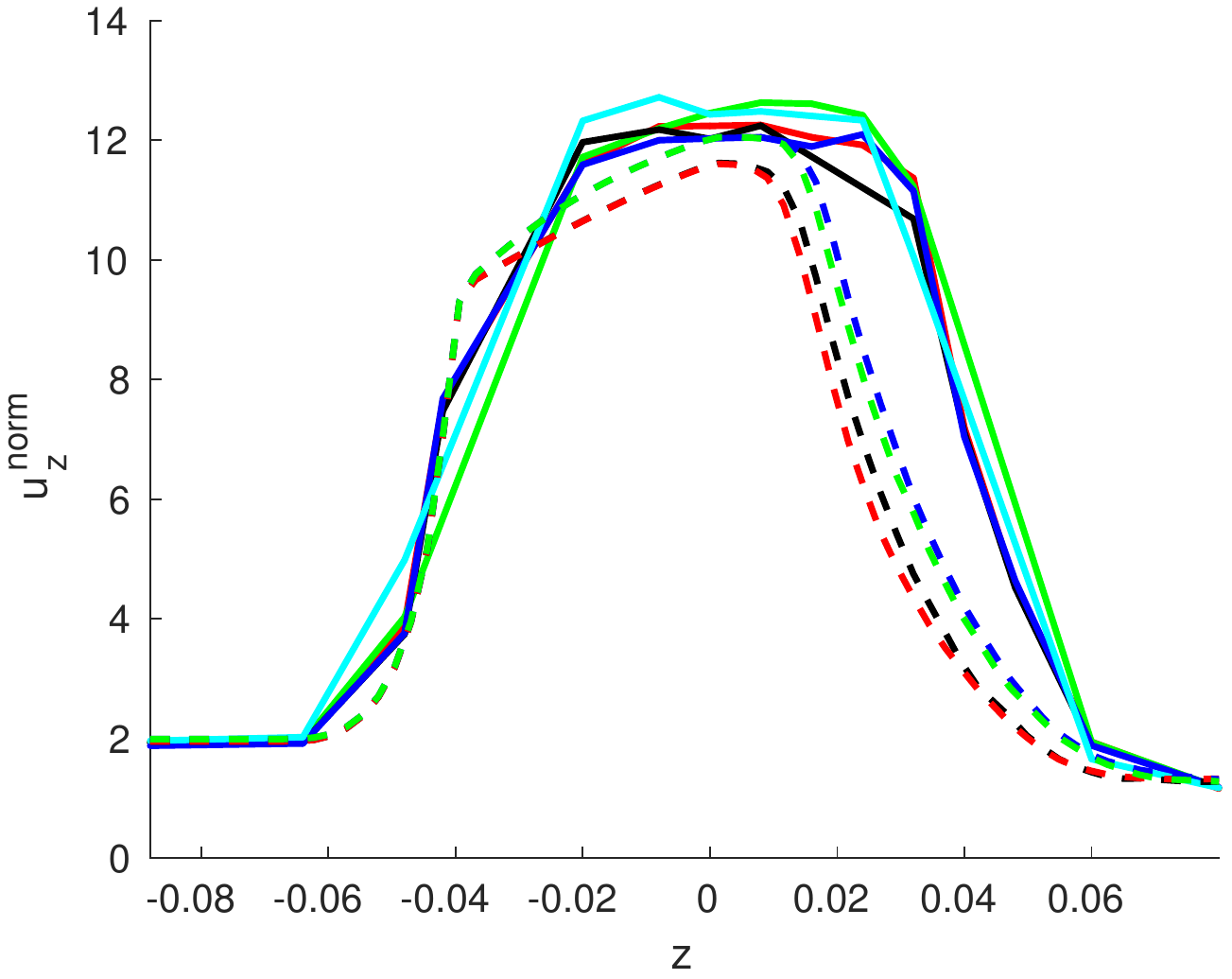}
      \put(20,78){\small{normalized axial velocity along $z$}}
      \end{overpic}
 \begin{overpic}[width=0.45\textwidth]{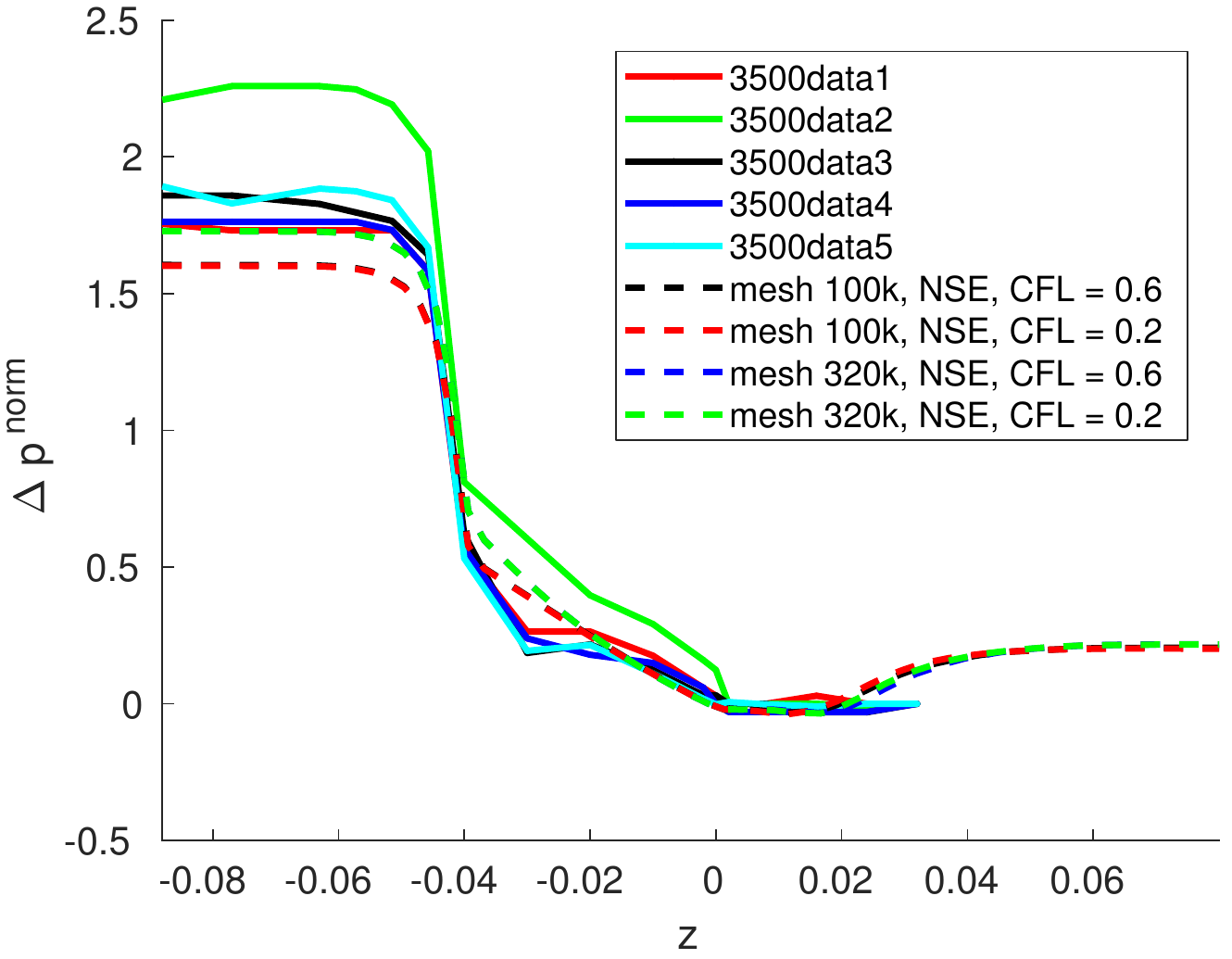}
       \put(20,78){\small{normalized pressure difference along $z$}}
      \end{overpic}
\caption{Case $Re_t = 3500$, NSE: comparison between experimental data (solid lines) 
and numerical results (dashed lines) for the normalized axial velocity \eqref{eq:uz} (left)
and the normalized pressure drop \eqref{eq:pz} (right) along the $z$. The numerical 
results were obtained with meshes $100k$ and $320k$ and
two different $CFL$ numbers.
The legend on the right is common to both graphs.}
\label{fig:FDA_NSE_3500_CFL}
\end{figure}

Next, we vary the $CFL$ number. Fig.~\ref{fig:FDA_NSE_3500_CFL} shows 
the normalized axial velocity \eqref{eq:uz} and the normalized pressure drop \eqref{eq:pz} 
along the $z$ axis from the NSE solution on meshes $100k$ and $320k$ (see Table \ref{tab:mesh_FDA})
obtained with $CFL$ number 0.2 and 0.6. We see that a lower $CFL$ number produces a slightly shorter jet.
In \cite{Nicoud} the same low sensitivity to the $CFL$ number was observed when the 
RK4 scheme is used for time differentiation, while the numerical results obtained with the
TFV4A scheme are extremely sensitive to the choice of $CFL$ number.
From now on, we will set the $CFL$ number to 0.6.


\begin{figure}[h]
\centering
 \begin{overpic}[width=0.45\textwidth]{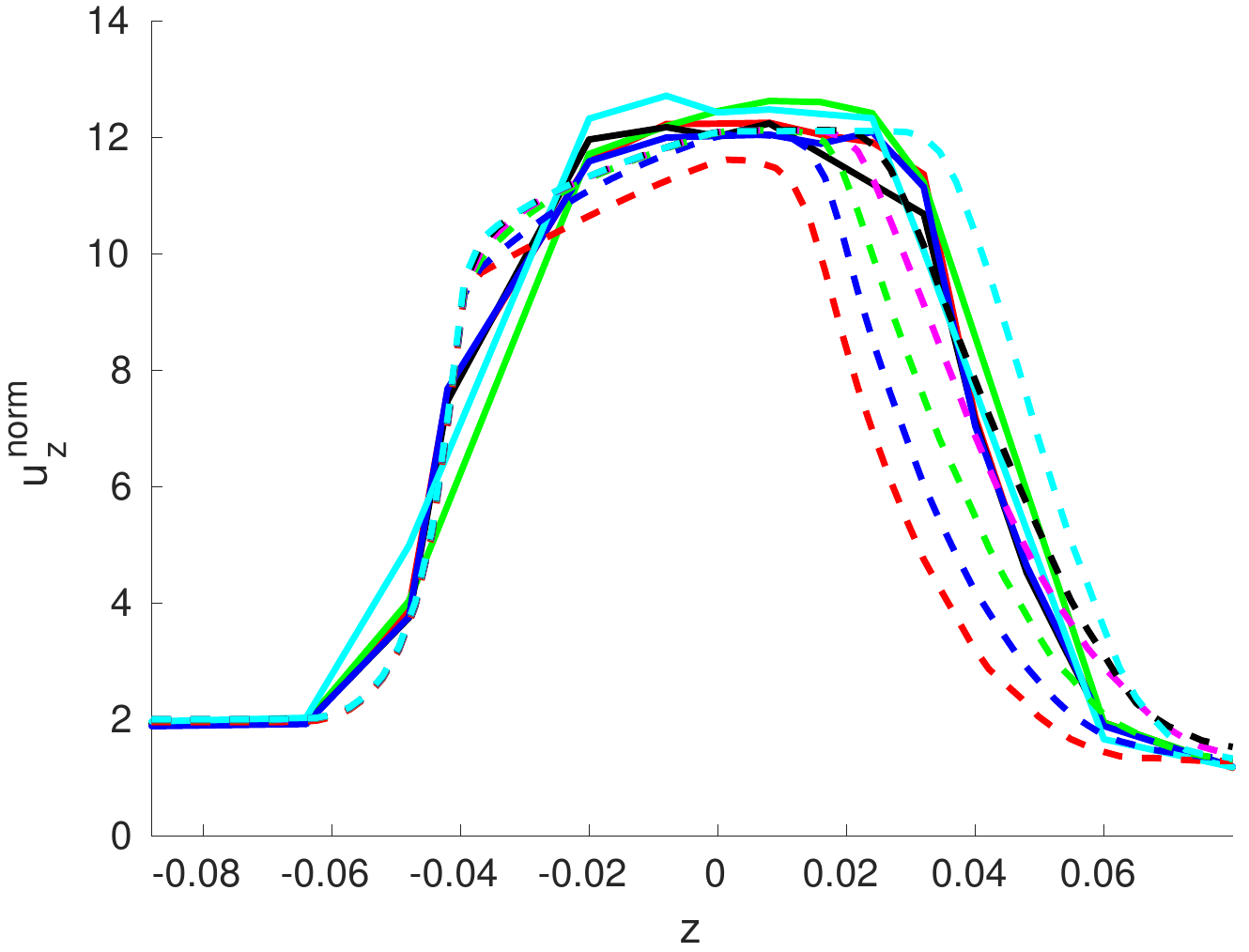}
      \put(17,78){\small{normalized axial velocity along $z$}}
      \end{overpic}
 \begin{overpic}[width=0.45\textwidth]{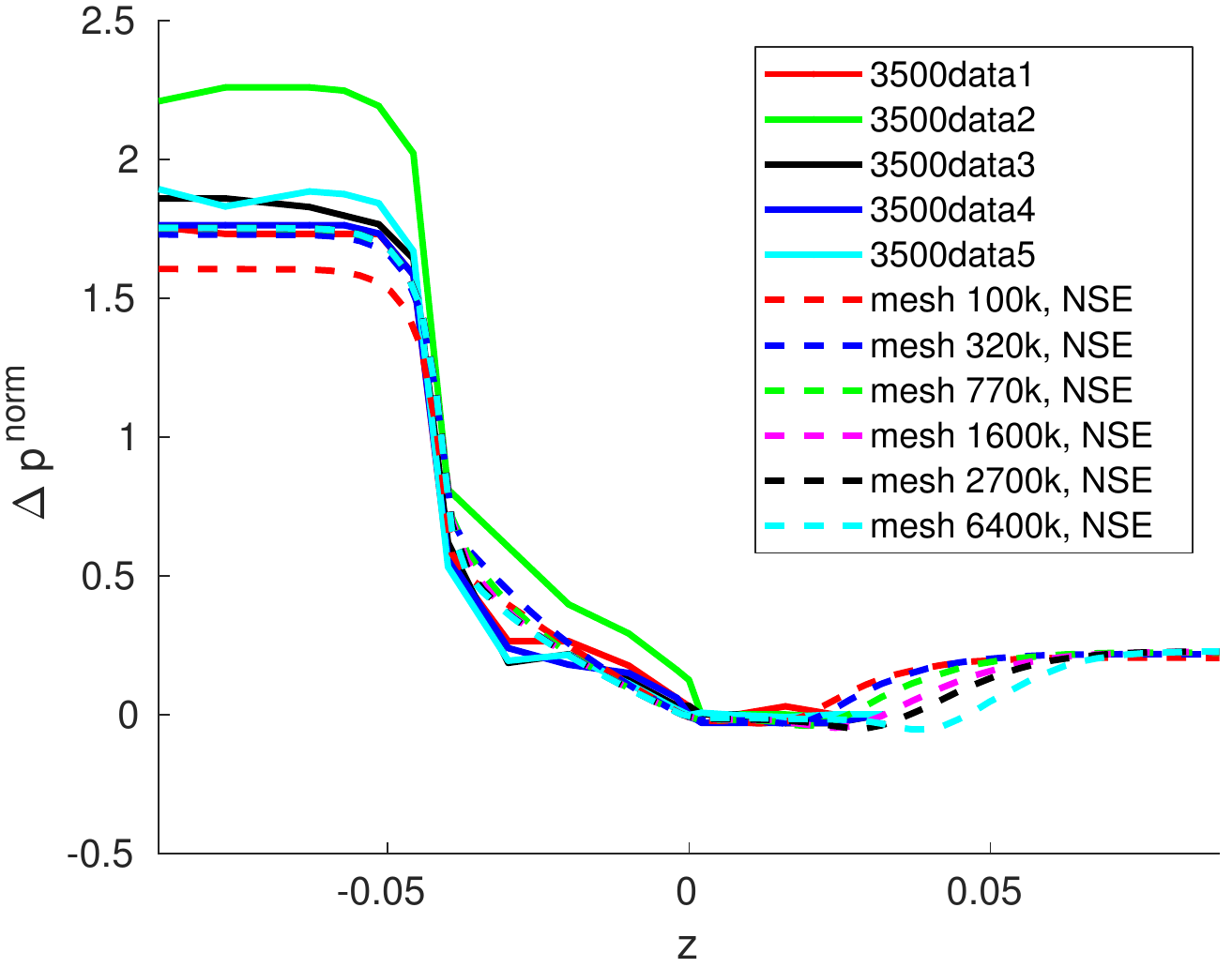}
       \put(13,78){\small{normalized pressure difference along $z$}}
      \end{overpic}
\caption{Case $Re_t = 3500$, NSE: comparison between experimental data (solid lines) 
and numerical results (dashed lines) for the normalized axial velocity \eqref{eq:uz} (left)
and the normalized pressure drop \eqref{eq:pz} (right) along the $z$. The numerical 
results were obtained with all the meshes in Table \ref{tab:mesh_FDA}.
The legend on the right is common to both graphs.}
\label{fig:FDA_NSE_3500}
\end{figure}

We now consider all the meshes in Table \ref{tab:mesh_FDA}. The associated results from the NSE solution
are displayed in Fig.~\ref{fig:FDA_NSE_3500}.
We see that the axial velocities computed with meshes $1600k$ and $2700k$ match quite well the measurements 
all along the portion of the $z$-axis under consideration ($-0.088<z<0.08$). 
These results are in agreement with those in \cite{BQV}, where it is shown that a 
DNS is possible with the tetrahedral meshes $1200k$ and $1900k$ (see Table \ref{tab:mesh_BQV_Nicoud} (left)).
Such meshes have number of cells, average number of cells along a throat diameter,
and values of $h_{min}$ comparable to meshes $1600k$ and $2700k$ used here. 
However, while in \cite{BQV} it is observed that instabilities in the NSE algorithm lead to 
a simulation breakdown for the meshes coarser than mesh $1200k$, no
instability arises with coarser meshes when a finite volume method is adopted.
From Fig.~\ref{fig:FDA_NSE_3500}, we see that the results obtained with meshes
$1600k$ and $2700k$ are in good agreement with the experimental data. However, 
we do not observe a clear convergence for the average velocity, i.e. the curves
do not get closer to each other as the mesh is refined. 
This suggests that at $Re_t$ = 3500 the flow is still very sensitive to perturbations, 
since it is close to the transitional regime. See also Sec.~\ref{sec:2000}.
Thus, much finer meshes are needed for a DNS, as reported in \cite{Nicoud}.
As for the pressure, we observe excellent agreement with the experimental
data for all the meshes but the coarsest, which underestimates the pressure drop
in the entrance ragion. In particular, Fig.~\ref{fig:FDA_NSE_3500} (right) shows
great agreement in the conical convergent, which could never be achieved with
the finite element method used in \cite{BQV}.

\begin{rem}
Additional tests (not reported for brevity) showed that small perturbations in the mesh 
strongly affect the results. For example, meshes very similar to $1600k$ and $2700k$
in terms of maximum non-orthogonality, average number of cells along the throat diameter, and mesh size
but with slightly different refinement near the sudden expansion produced very long jets. 
Thus, we stress the importance of choosing the correct features in mesh generation, as already pointed
out in \cite{stewartp,Passerini2013,BQV,Nicoud,Bergensen2018,Janiga2014,Delorme2013,Bhushan2013,Chabannes2017,Fehn2018,Nicoud2018,Wu2019}.
\end{rem}

\begin{figure}[h]
\centering
 \begin{overpic}[width=0.45\textwidth]{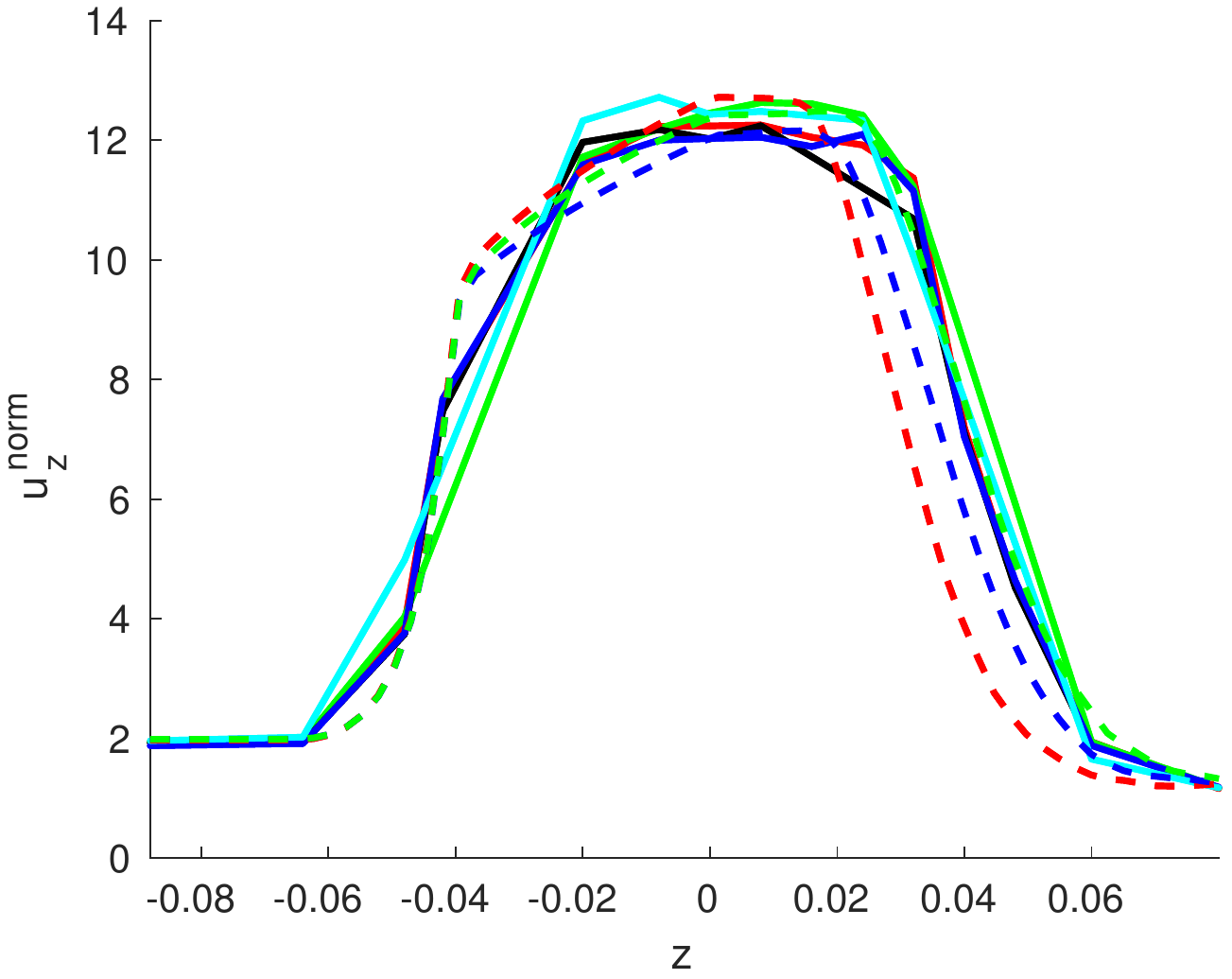}
      \put(17,78){\small{normalized axial velocity along $z$}}
      \end{overpic}
 \begin{overpic}[width=0.45\textwidth]{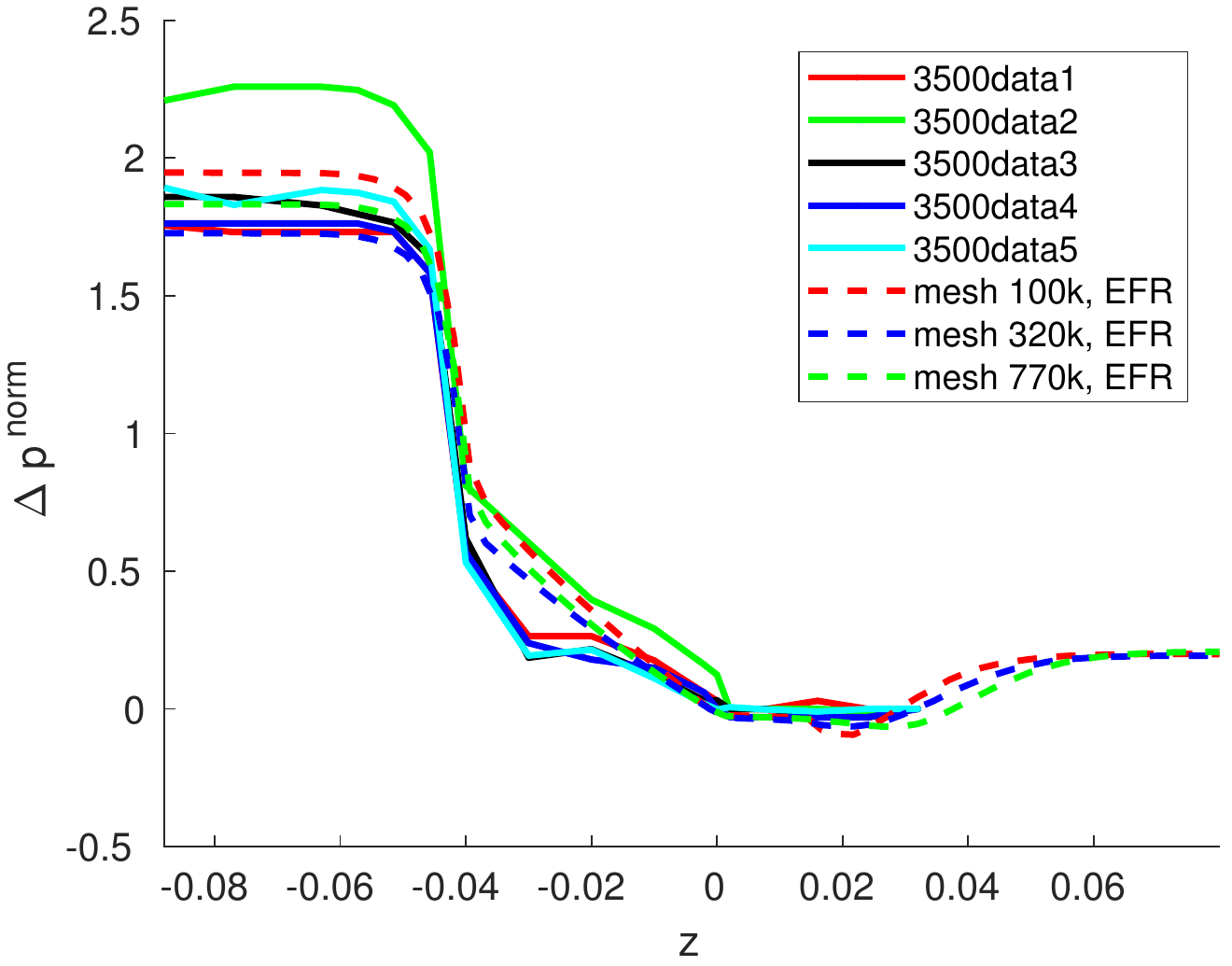}
        \put(13,78){\small{normalized pressure difference along $z$}}
      \end{overpic}
\caption{
Case $Re_t = 3500$, EFR: comparison between experimental data (solid lines) 
and numerical results (dashed lines) for the normalized axial velocity \eqref{eq:uz} (left)
and the normalized pressure drop \eqref{eq:pz} (right) along the $z$. The numerical 
results were obtained for all the meshes in Table \ref{tab:mesh_FDA}
coarser than mesh $1600k$.
The legend on the right is common to both graphs.}
\label{fig:FDA_EFR_3500}
\end{figure}

We have noticed that the choice $\chi = \chi_2$ introduces slightly more
diffusion that $\chi = \chi_1$ for all the coarser meshes under consideration, thereby providing
results in better agreement with the experimental data. Thus, we set $\chi = \chi_2$ whenever
we use the EFR algorithm.
In Fig.~\ref{fig:FDA_EFR_3500}, we report the comparison between EFR solutions and experimental data in terms of 
normalized axial velocity \eqref{eq:uz} and pressure difference \eqref{eq:pz}
for all the meshes in Table \ref{tab:mesh_FDA} coarser than mesh $1600k$. 
We see that the results computed with mesh $770k$ are 
in very good agreement with the measurements all along the $z$-axis. 
With meshes $100k$ and $320k$, the EFR algorithm provides a jet that starts to break 
down closer to the sudden expansion.  
Nevertheless, for all three meshes we notice that the EFR algorithm significantly improves the agreement with the 
experimental data with respect to NSE: compare Fig.~\ref{fig:FDA_NSE_3500} with \ref{fig:FDA_EFR_3500}. 

Despite the fact that a Finite Volume approximation is used also in \cite{Nicoud}, 
a direct comparison between our results and those in \cite{Nicoud} is complicated 
because we use difference LES approaches. 
In \cite{Nicoud}, the authors chose to use meshes tetrahedral elements and
avoid local refinement.
As a result, their meshes feature a higher level of refinement than ours. 
Moreover, the results in \cite{Nicoud} are very sensitive to the time discretization scheme, 
$CFL$ condition and grid resolution. 
However, it is shown that introducing fluctuations in the upstream boundary condition 
drastically reduces this sensitivity and makes the prediction of the jet break-down point very robust.


\subsubsection{Case $Re_t = 5000$}\label{sec:5000}

At $Re_t = 5000$ all the experiments in \cite{hariharang} observed 
turbulence downstream of the sudden expansion with a reproducible jet breakdown point.
Just like for the $Re_t = 3500$ case, we will first check the results given by the NSE algorithm 
and then consider the EFR algorithm.

\begin{figure}[h]
\centering
 \begin{overpic}[width=0.45\textwidth]{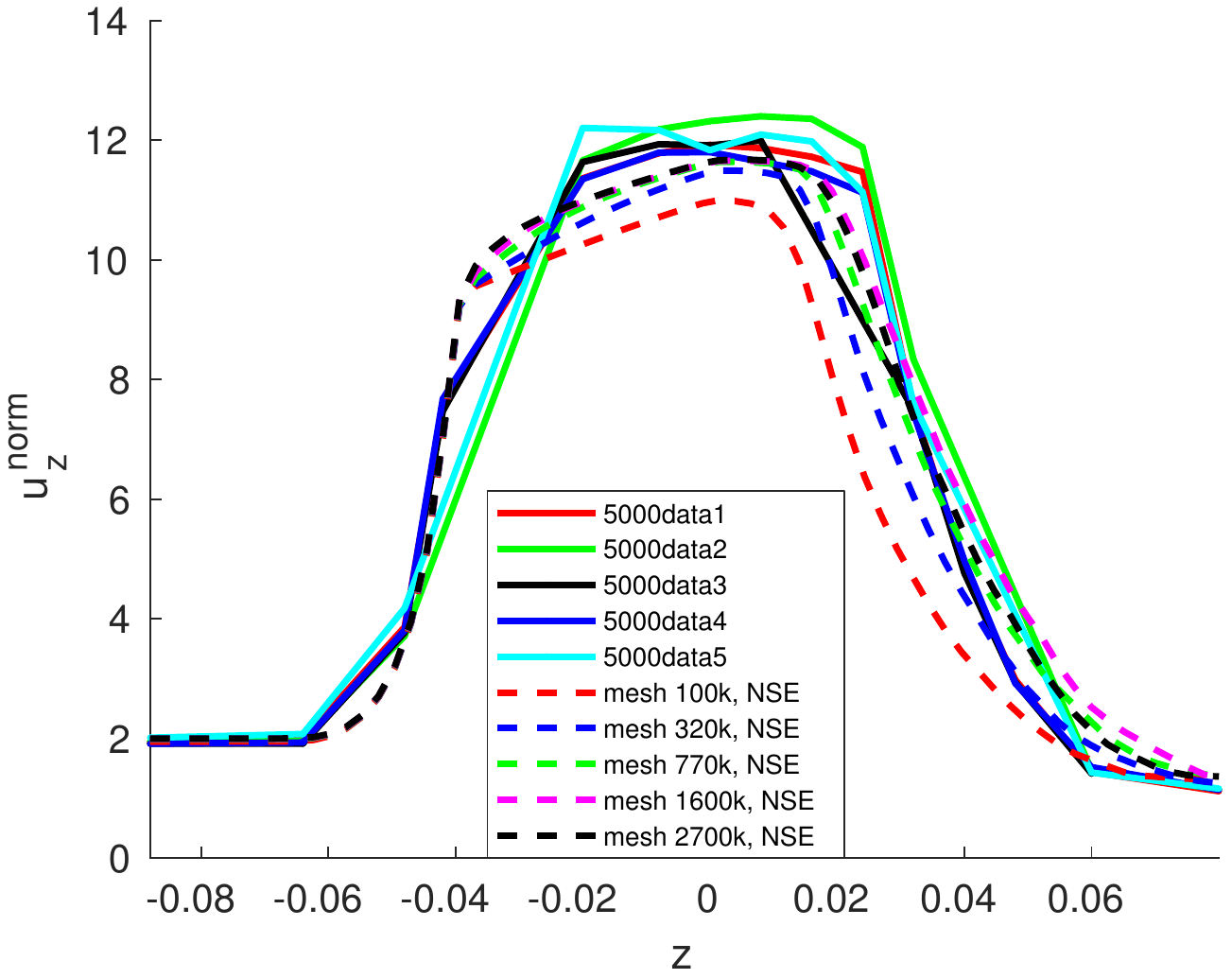}
      \put(17,78){\small{normalized axial velocity along $z$}}
      \end{overpic}
 \begin{overpic}[width=0.45\textwidth]{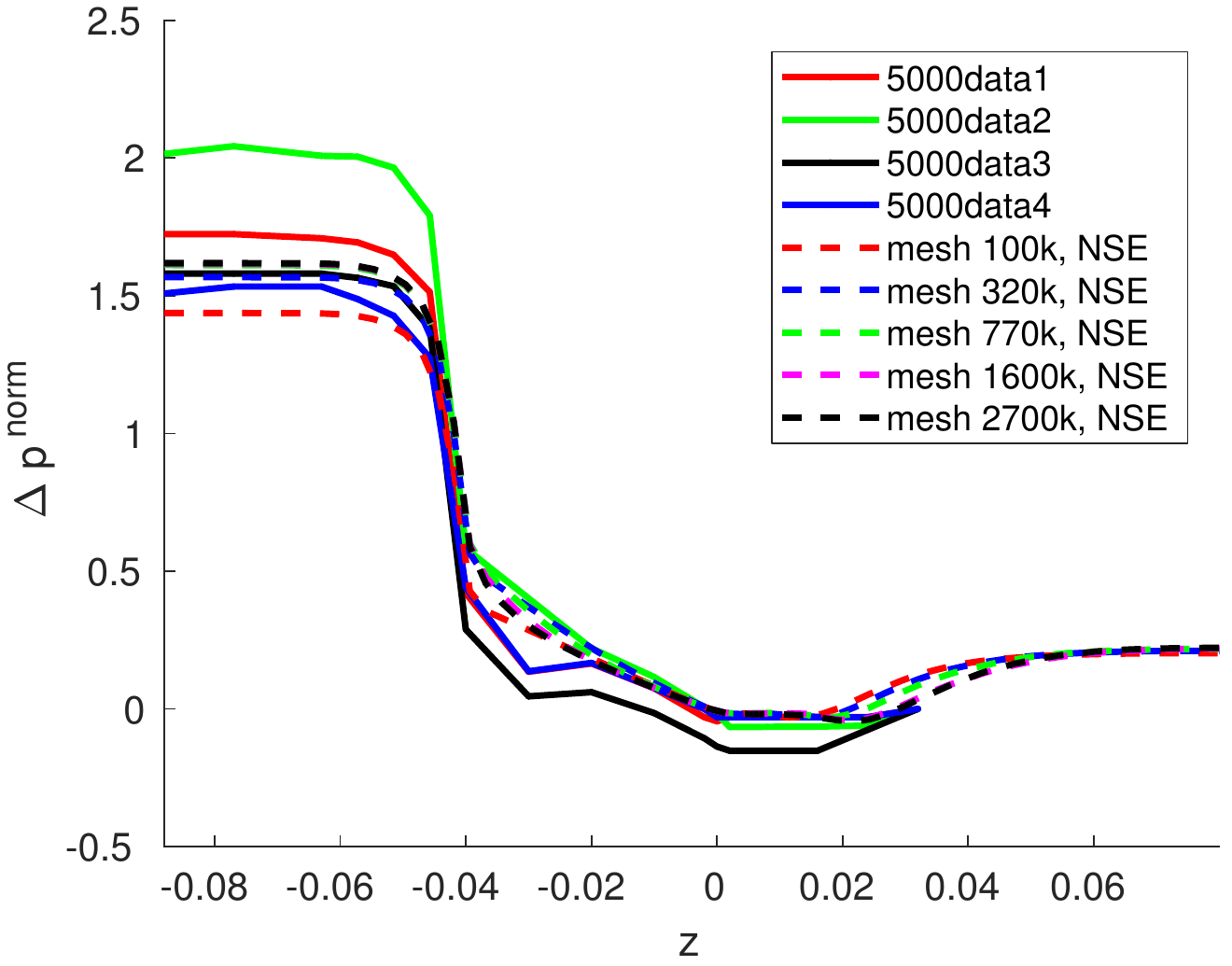}
       \put(13,78){\small{normalized pressure difference along $z$}}
      \end{overpic}
\caption{Case $Re_t = 5000$, NSE: comparison between experimental data (solid lines) 
and numerical results (dashed lines) for the normalized axial velocity \eqref{eq:uz} (left)
and the normalized pressure drop \eqref{eq:pz} (right) along the $z$. The numerical 
results were obtained for mesh $2700k$ and coarser (see Table \ref{tab:mesh_FDA}).}
\label{fig:FDA_NSE_5000}
\end{figure}

Fig.~\ref{fig:FDA_NSE_5000} shows the normalized axial velocity \eqref{eq:uz} 
and the normalized pressure drop \eqref{eq:pz} along the $z$ axis computed from 
the NSE solution on all the meshes in Table \ref{tab:mesh_FDA} coarser than 
mesh $2700k$. First, we observe that that the axial velocities computed 
with meshes $1600k$ and $2700k$ match quite well the measurements all along 
the portion of the $z$-axis under consideration ($-0.088<z<0.08$), just like 
the $Re_t$ = 3500 case. 
Unlike the $Re_t$ = 3500 case, 
we do observe convergence as the mesh is refined. 
We remark that finite element method in \cite{BQV} 
did not provide results in good agreement with the experimental data with a
tetrahedral mesh $3000k$ (see Table \ref{tab:mesh_BQV_Nicoud}). We believe
this is due to the fact that mesh $3000k$ in \cite{BQV} has a lower average number of cells throat diameter 
and a greater $h_{max}$ than meshes $1600k$ and $2700k$.
See Tables \ref{tab:mesh_FDA} and \ref{tab:mesh_BQV_Nicoud}.
As for the pressure, from Fig.~\ref{fig:FDA_NSE_5000} (right)
we see excellent agreement with the experimental data for all the meshes but the coarsest, which 
slightly underestimates the pressure drop in the entrance ragion.

\begin{figure}[h]
\centering
 \begin{overpic}[width=0.45\textwidth]{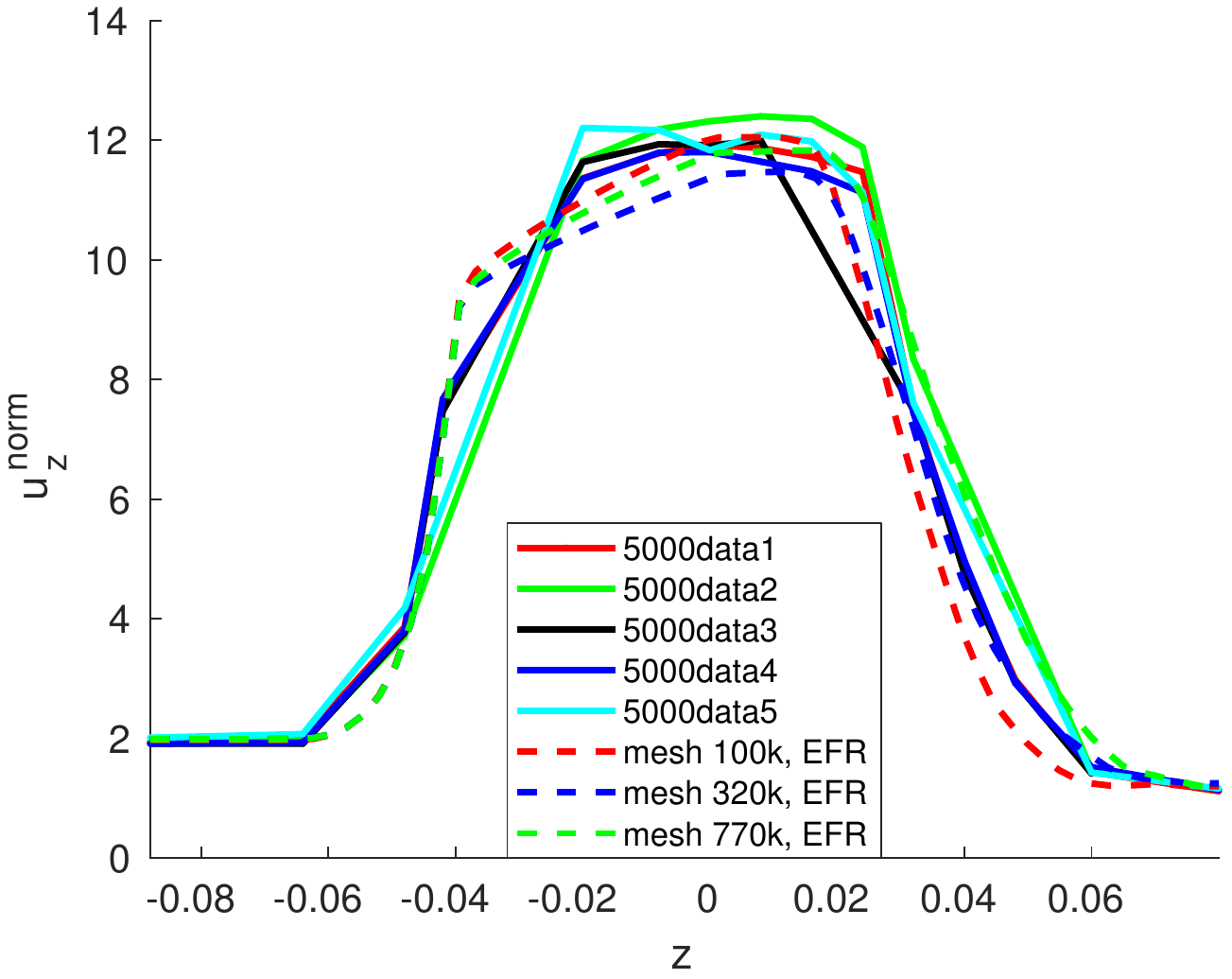}
      \put(17,78){\small{normalized axial velocity along $z$}}
      \end{overpic}
 \begin{overpic}[width=0.45\textwidth]{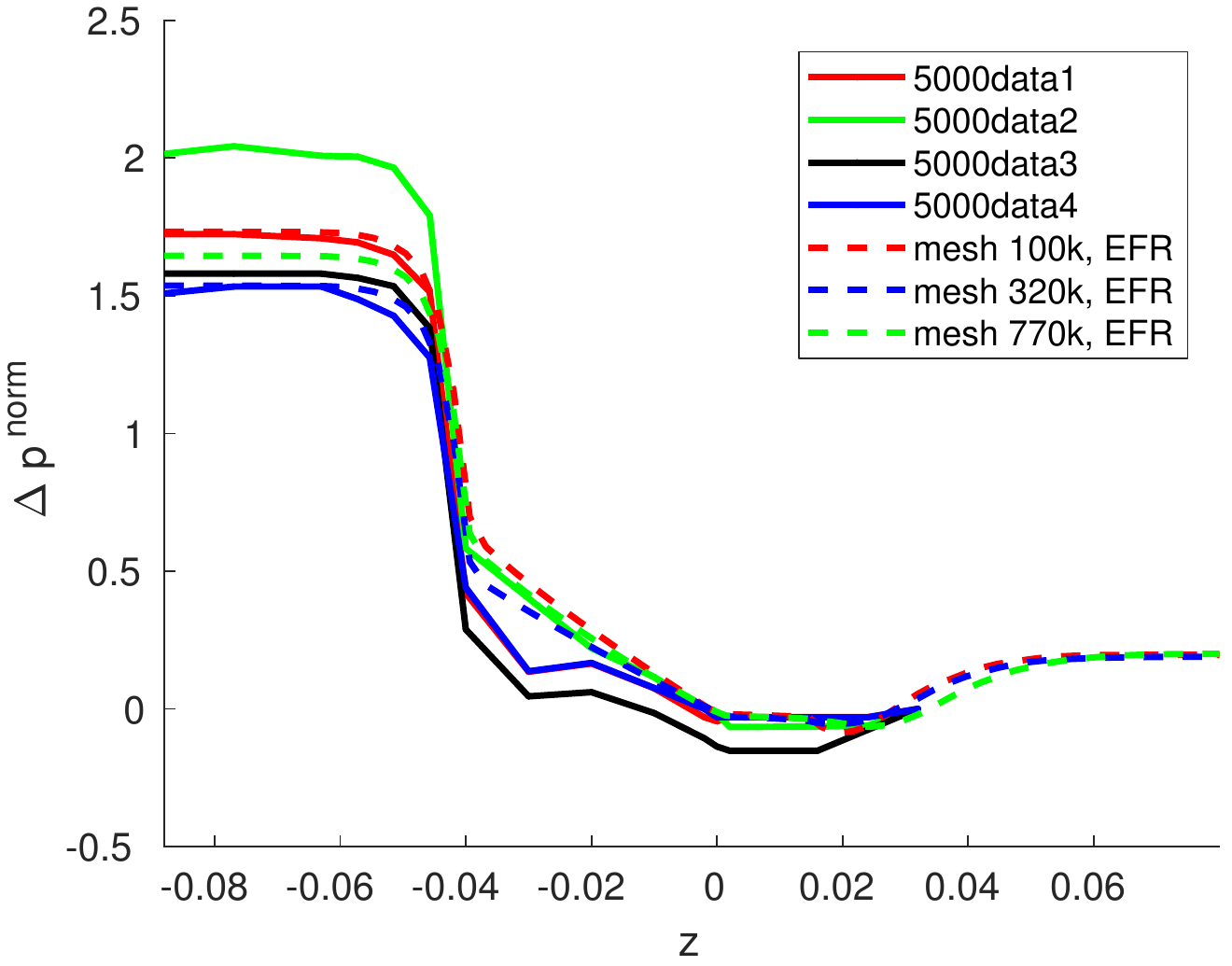}
       \put(13,78){\small{normalized pressure difference along $z$}}
      \end{overpic}
\caption{
Case $Re_t = 5000$, EFR: comparison between experimental data (solid lines) 
and numerical results (dashed lines) for the normalized axial velocity \eqref{eq:uz} (left)
and the normalized pressure drop \eqref{eq:pz} (right) along the $z$. The numerical 
results were obtained for all the meshes in Table \ref{tab:mesh_FDA}
coarser than mesh $1600k$.}
\label{fig:FDA_EFR_5000}
\end{figure}

Next, we consider the EFR algorithm. We report the comparison between computed 
and measured normalized axial velocity \eqref{eq:uz} and pressure drop \eqref{eq:pz} 
in Fig.~\ref{fig:FDA_EFR_5000}. 
We see that the results computed with meshes $320k$ and $770k$ are in very good 
agreement with the measurements all along the $z$-axis. 
It is remarkable that the results given by the EFR algorithm on a very coarse mesh
like $320k$ match so well the experimental data. As we will see in the next subsection, 
this is the case also for $Re_t = 6500$. A possible reason why the results given by the EFR algorithm
with mesh $320k$ were not as good for $Re_t = 3500$ might be that $Re_t = 3500$ is still close to 
the transitional regime. 
From Fig.~\ref{fig:FDA_EFR_5000} (left) we see that the EFR algorithm
with mesh $100k$ provides a jet that starts to break down closer to the sudden expansion,
just like the $Re_t = 3500$ case.

From the tests performed for the $Re_t = 3500$ and $Re_t = 5000$ cases, we learned that 
the computed jet breakdown location is very sensitive to 
grid features (kind of element, non-orthogonality, local refinement, aspect ratio) and 
interpolation scheme for the convective term.  
Based on our experience, a robust FV-based CFD framework requires: (i) 
non-orthogonality and low skewness hexahedral meshes having cells with a contained element aspect ratio, and 
(ii) the CD scheme for the interpolation of the convective term. With this framework, 
we obtained favorable comparisons with the experimental data 
without introducing inlet perturbations as in \cite{Nicoud}.

\subsubsection{Case $Re_t = 6500$}

The third flow regime we focus on features the highest Reynolds number 
considered by FDA: $Re_t = 6500$. Also in this case, turbulence downstream of the sudden expansion 
was observed in all the experiments reported in \cite{hariharang}
with a reproducible jet breakdown point. 

The $Re_t = 6500$ case is thoroughly investigated in \cite{Janiga2014}. 
The author uses the Smagorinsky subgrid-scale model \cite{Smagorinsky1963} 
and combines it with a FV method implemented in the commercial CFD code ANSYS 
Fluent. Great agreement with the experimental data is obtained 
with a very fine hexahedral block-structured mesh that has roughly 9 millions
elements ($h_{min}$ =2.5e-5) and excellent features in terms of non-orthogonality and skewness.
We notice that the meshes used in \cite{Janiga2014} are very similar to those ones 
we use in this work. 

\begin{figure}[h]
\centering
 \begin{overpic}[width=0.45\textwidth]{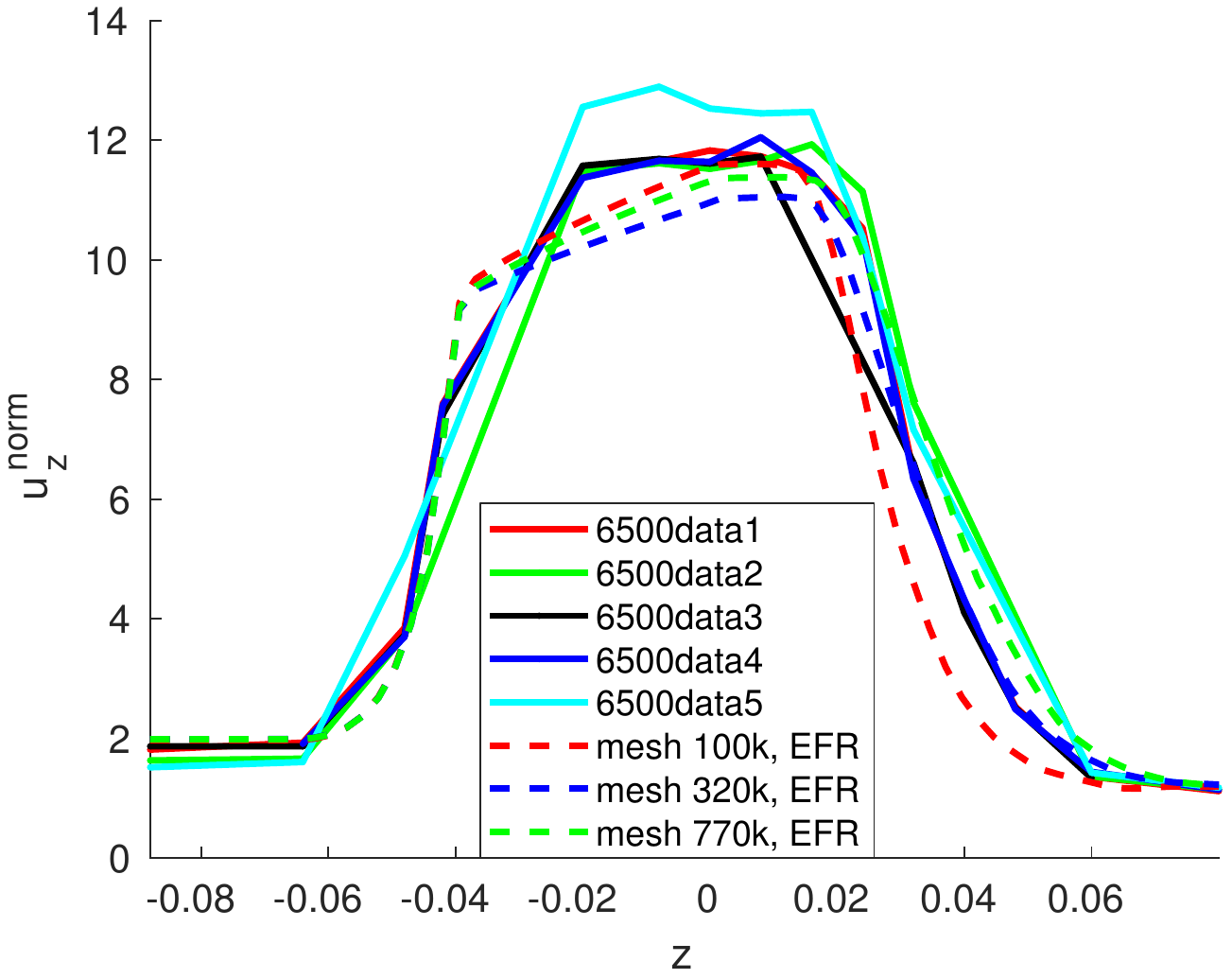}
      \put(17,78){\small{normalized axial velocity along $z$}}
      \end{overpic}
 \begin{overpic}[width=0.45\textwidth]{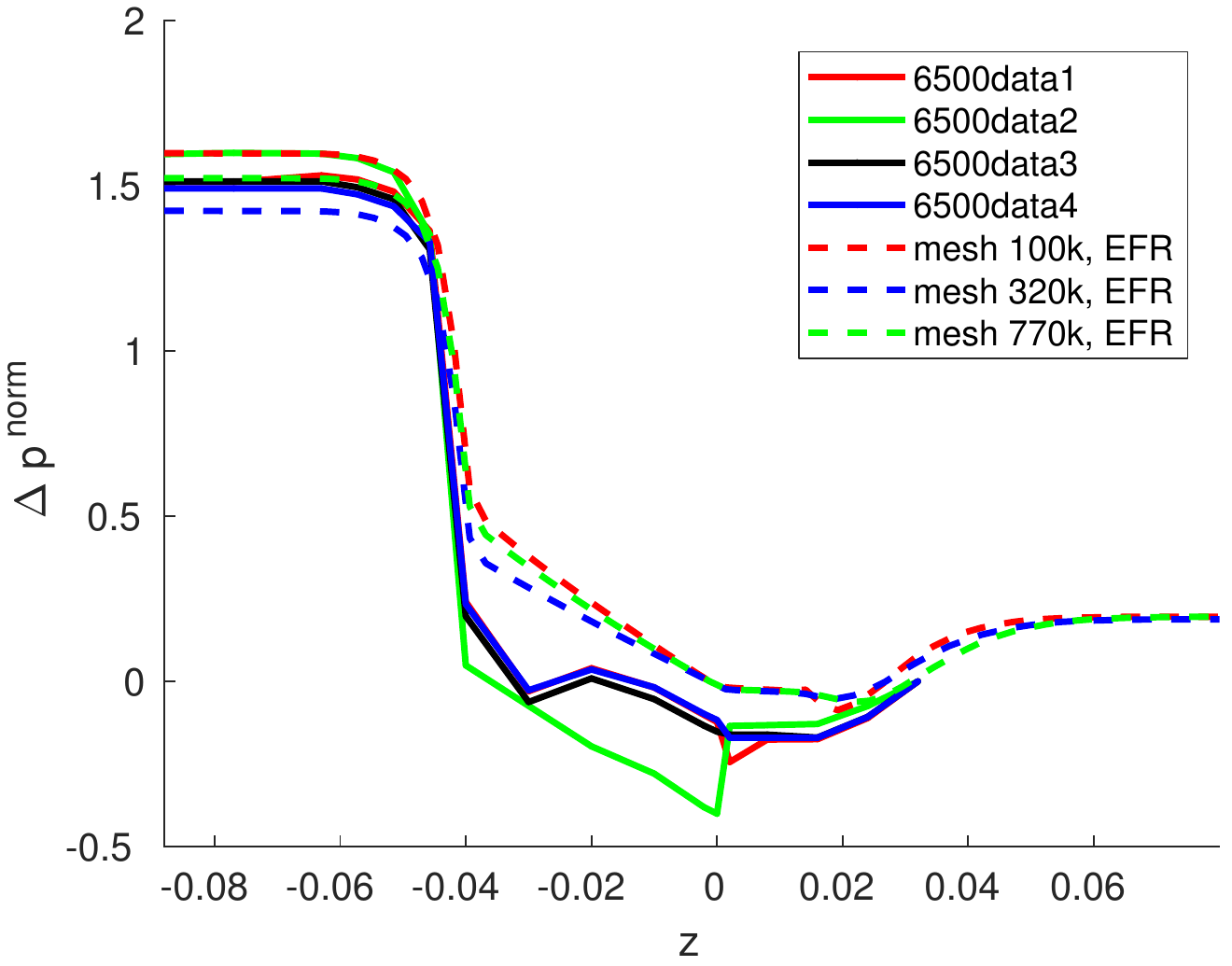}
       \put(13,78){\small{normalized pressure difference along $z$}}
      \end{overpic}
\caption{
Case $Re_t = 6500$, EFR: comparison between experimental data (solid lines) 
and numerical results (dashed lines) for the normalized axial velocity \eqref{eq:uz} (left)
and the normalized pressure drop \eqref{eq:pz} (right) along the $z$. The numerical 
results were obtained for all the meshes in Table \ref{tab:mesh_FDA}
coarser than mesh $1600k$.}
\label{fig:FDA_EFR_6500}
\end{figure}

Given the good EFR results shown in Sec.~\ref{sec:3500} and \ref{sec:5000}, for $Re_t = 6500$
we limit the investigation to the EFR algorithm with $\chi = \chi_2$ and all the meshes in Table \ref{tab:mesh_FDA}
coarser than mesh $1600k$. We report the comparison between computed 
and measured normalized axial velocity \eqref{eq:uz} and pressure drop \eqref{eq:pz} 
in Fig.~\ref{fig:FDA_EFR_6500}. From Fig.~\ref{fig:FDA_EFR_6500} (left)
we see that the axial velocity computed with meshes $330k$ and $770k$ is 
in very good agreement with the measurements all along the $z$-axis. 
As for mesh $100k$, we observe the same trend as for the $Re_t = 3500, 5000$ cases:
the computed jet starts to break down closer to the sudden expansion, while 
still in good agreement with the experimental data in the rest of the domain.
As for the pressure, from Fig.\ref{fig:FDA_EFR_6500} (right) we observe a significant 
overestimation in the conical convergent respect to the experimental data for all the meshes. 
A relaxation of the pressure (see Remark \ref{rem:p_relax}) might lead to a better agreement. 

\subsubsection{Case $Re_t = 2000$}\label{sec:2000}

The transitional case ($Re_t = 2000$) proved to be challenging from both the experimental and
computational sides, thus we present it last. 
From the experimental side, the interlaboratory velocity data in \cite{hariharang} 
showed significant differences downstream of the sudden expansion, in particular
concerning the jet breakdown point.
This was attributed mainly to a 10\% higher flow rate (and consequently higher $Re_t$), 
which caused premature jet breakdown in two experiments out of five \cite{hariharang}. 
However, minor differences in the fabricated geometrical models and inlet perturbation levels played a role also. 
From the numerical point of view, we found the results to be very sensitive to the mesh 
refinement level.

\begin{figure}[h]
\centering
 \begin{overpic}[width=0.45\textwidth]{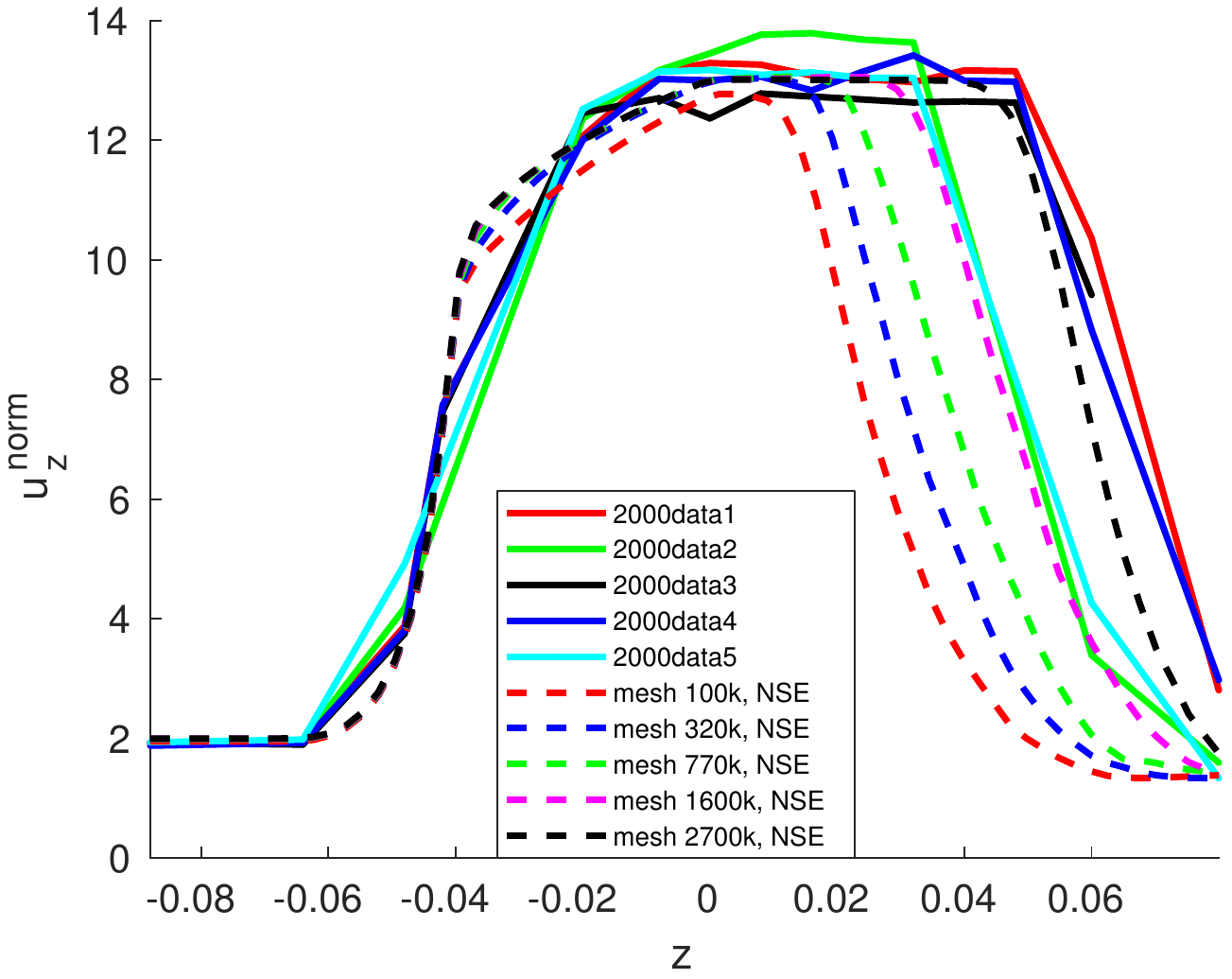}
    \put(17,78){\small{normalized axial velocity along $z$}}
      \end{overpic}
 \begin{overpic}[width=0.45\textwidth]{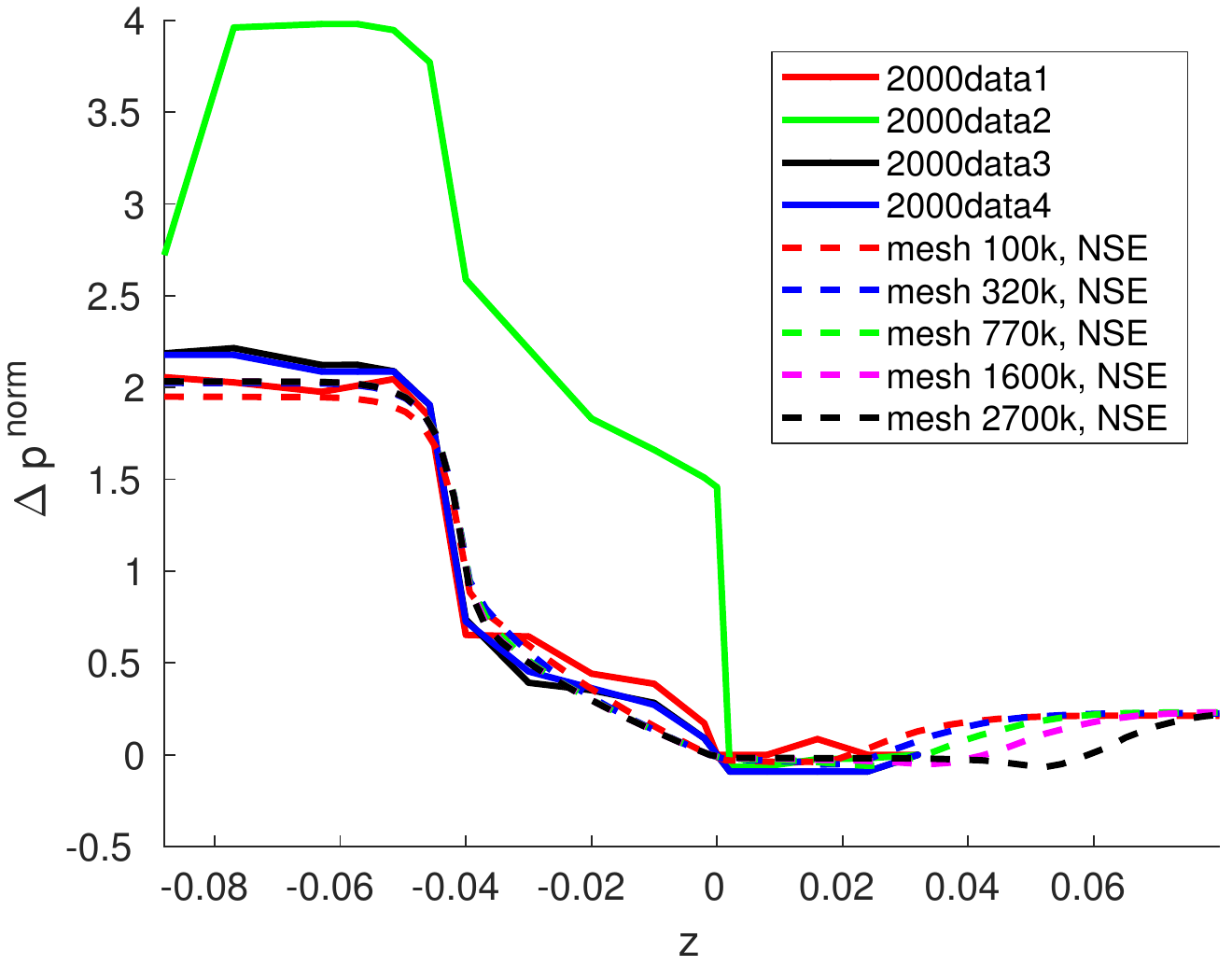}
       \put(13,78){\small{normalized pressure difference along $z$}}
      \end{overpic}
\caption{
Case $Re_t = 2000$, NSE: comparison between experimental data (solid lines) 
and numerical results (dashed lines) for the normalized axial velocity \eqref{eq:uz} (left)
and the normalized pressure drop \eqref{eq:pz} (right) along the $z$. The numerical 
results were obtained for mesh $2700k$ and coarser (see Table \ref{tab:mesh_FDA}).
}
\label{fig:FDA_NSE_2000}
\end{figure}

For this test, we look first at the results obtained with the NSE algorithm.
Fig.~\ref{fig:FDA_NSE_2000} reports the normalized axial velocity \eqref{eq:uz}
and the normalized pressure drop \eqref{eq:pz} along the $z$
for mesh $2700k$ and coarser (see Table \ref{tab:mesh_FDA}).
From Fig.~\ref{fig:FDA_NSE_2000} (left), we notice that the axial velocity
computed with mesh $1600k$ (resp., $2700k$) match quite well the measurements 
giving a shorter (resp., longer) jet. 
This agrees with what reported in \cite{Passerini2013}, where it is shown that the NSE
results obtained with a mesh with 2.5 millions tetrahedra compare well
with the experimental data giving a longer jet. However, the numerical
axial velocity computed with meshes $1600k$ and $2700k$ are much further apart
for $Re_t = 2000$ than for $Re_t = 3500$. Compare Fig.~\ref{fig:FDA_NSE_2000} (left)
with Fig.~\ref{fig:FDA_NSE_3500} (left). Thus, we cannot argue that either mesh is
refined enough for this test case. 
From Fig.~\ref{fig:FDA_NSE_2000} (right), we see that the simulated pressure drop 
is in very good agreement with the experimental data all along the $z$-axis
and regardless of the mesh.

\begin{figure}[h]
\centering
 \begin{overpic}[width=0.45\textwidth]{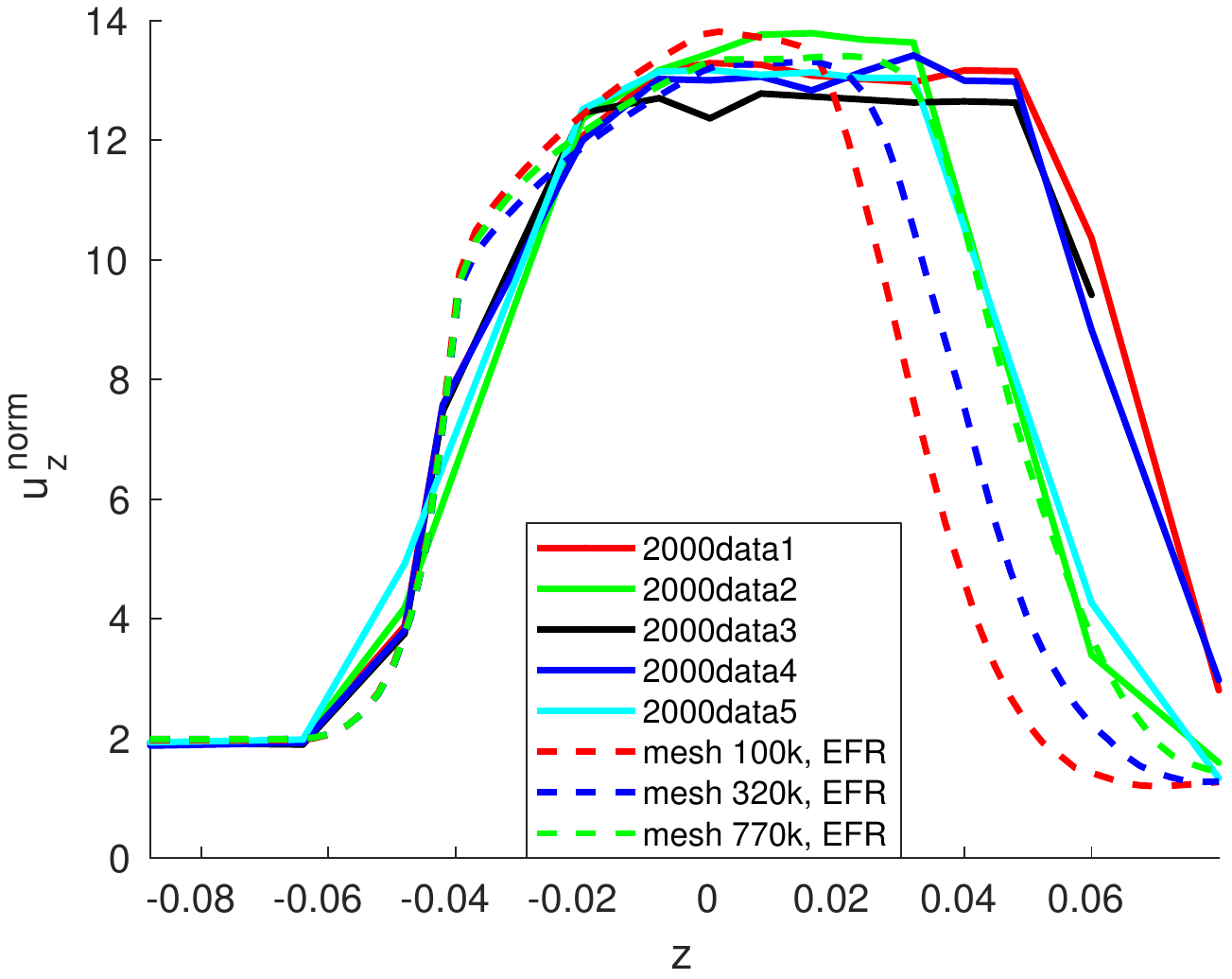}
      \put(17,78){\small{normalized axial velocity along $z$}}
      \end{overpic}
 \begin{overpic}[width=0.45\textwidth]{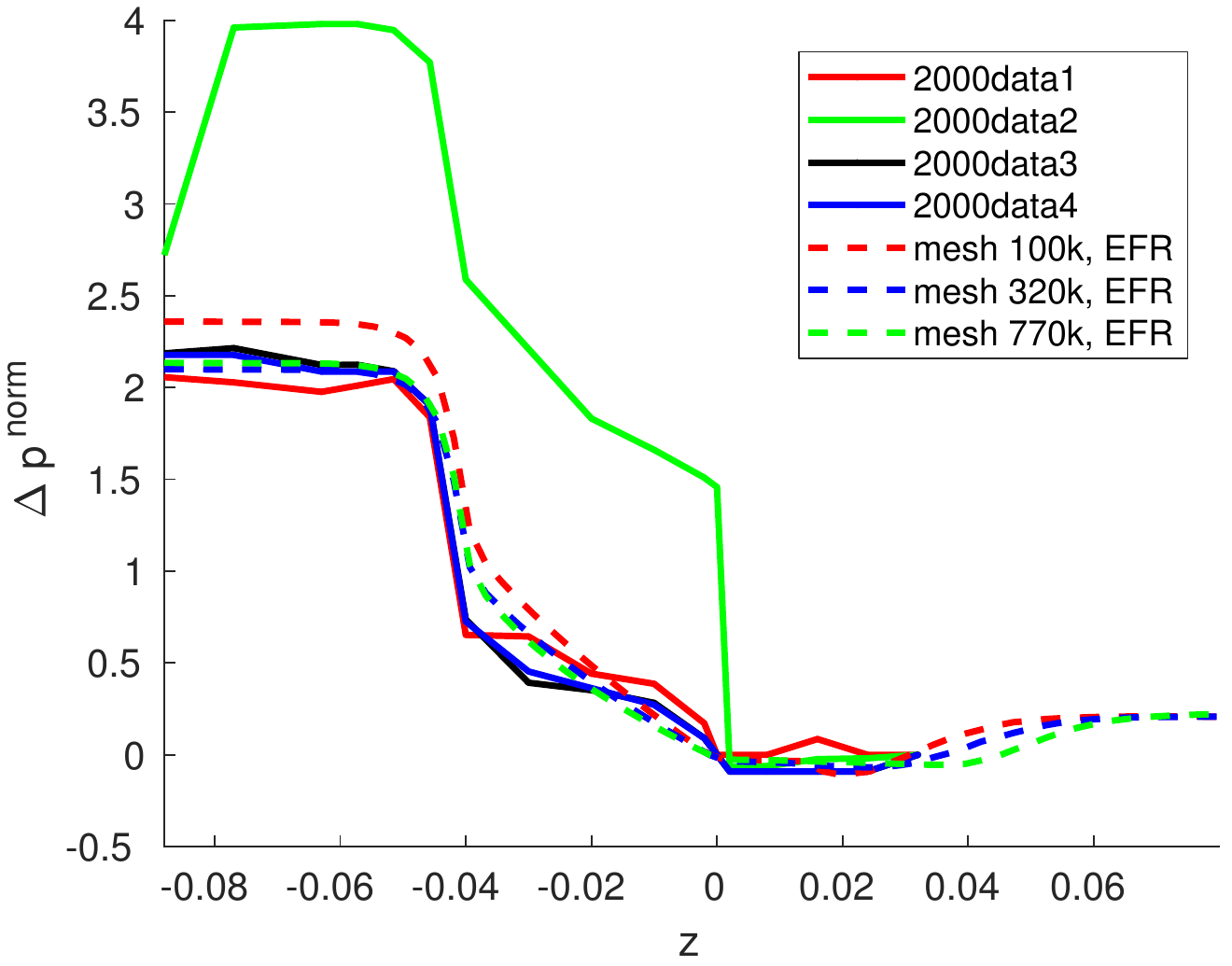}
       \put(13,78){\small{normalized pressure difference along $z$}}
      \end{overpic}
\caption{
Case $Re_t = 2000$, EFR: comparison between experimental data (solid lines) 
and numerical results (dashed lines) for the normalized axial velocity \eqref{eq:uz} (left)
and the normalized pressure drop \eqref{eq:pz} (right) along the $z$. The numerical 
results were obtained for all the meshes in Table \ref{tab:mesh_FDA}
coarser than mesh $1600k$.}
\label{fig:FDA_EFR_2000}
\end{figure}

Next, we consider the EFR algorithm. We report the comparison between computed 
and measured normalized axial velocity \eqref{eq:uz} and pressure drop \eqref{eq:pz} 
in Fig.~\ref{fig:FDA_EFR_2000}. We see that the axial velocity computed with mesh $770k$
is in very good agreement with the measurements giving a shorter jet.
To reach such an agreement with the NSE algorithm we had to use mesh $1600k$. 
However, just like for the NSE algorithm, the numerical axial velocity computed with 
the different meshes are further apart for $Re_t = 2000$ than for all the Reynolds
numbers considered in the previous sections. Compare Fig.~\ref{fig:FDA_EFR_2000} (left)
with Fig.~\ref{fig:FDA_EFR_3500} (left), \ref{fig:FDA_EFR_5000} (left), and \ref{fig:FDA_EFR_6500} (left).
These results confirm the difficulty in the numerical simulation of the transitional regime.

\section{Conclusions and Perspectives}\label{sec:conclusions}

We showed the effectiveness of a FV-based EFR algorithm in simulating flow problems at moderately large Reynolds numbers. Nonlinear filtering stabilizes marginally resolved scales without over-diffusing, thereby allowing to use less degrees of freedom than required by a DNS. To select the regions of the domain where filtering is needed, we employ a nonlinear differential low-pass filter.
The interest in the Finite Volume approximation is due to the fact that it has been widely used in the LES context. 
However, the application of the Leray model in a FV framework has been unexplored.

In order to showcase the features of our approach, we presented a computational study related to two benchmarks: 2D flow past a cylinder and a 3D benchmark from the FDA. 
With reference to the FDA benchmark, we performed a complete characterization of the flow at Reynolds numbers starting from 2000 up to 6500. We proposed a new formula based on physical and numerical arguments to tune a value of $\chi$ that leads to very good agreement with the experimental measurements. Several meshes were considered to understand how under-refined the mesh can be while still capturing the physical average quantities. Furthermore, we investigated the impact of the mesh features as well as the discretization of the convective term on the results obtained.

As a follow-up of the present work, we are going to develop a Leray Reduced Order Model (ROM) within a Finite Volume framework. We are also interested in coupling the Leray model with an elasticity model
to simulate fluid-structure interaction problems which are ubiquitous 
in science and engineering.

\section{Acknowledgements}\label{sec:acknowledgements}
We acknowledge the support provided by the European Research Council Executive Agency by the Consolidator Grant project AROMA-CFD ``Advanced Reduced Order Methods with Applications in Computational Fluid Dynamics" - GA 681447, H2020-ERC CoG 2015 AROMA-CFD and INdAM-GNCS projects.
This work was also partially supported by US National Science Foundation through grant  DMS-1620384.

\bibliographystyle{plain}
\bibliography{latexbi,references,lifev} 
\end{document}